%% file: main.tex
\documentclass[oneside]{sn-jnl} 

\usepackage{amsfonts, mathtools, graphicx} 
\usepackage{amsthm}
\usepackage{booktabs}

\usepackage{caption}
\usepackage{xcolor}

\usepackage[
backend=biber,
bibstyle=ieee,
citestyle=numeric,
natbib=true,
sorting=nyt,
eprint=false,
isbn = false,
url = false,
]{biblatex}
\AtEveryBibitem{\clearlist{language}} 

\usepackage{tikz}
\usetikzlibrary{math}
\usetikzlibrary{arrows.meta}
\usepackage{pgfplots}
\usetikzlibrary{plotmarks}

\definecolor{cyanEuler}{rgb}{0,1,1}
\definecolor{greenGS}{rgb}{0,1,0}
\definecolor{blueTalbot}{rgb}{0,0,1}
\definecolor{bordeauxCME}{rgb}{0.703,0.211,0.098}
\definecolor{orangeZakian}{rgb}{0.965,0.473,0.039}
\definecolor{redTAME}{rgb}{0.977,0.086,0.281}
\definecolor{greenTAME}{rgb}{0,0.5,0.18}
\definecolor{purpleTAME}{rgb}{0.39, 0, 066}
\definecolor{purpleCME}{rgb}{0.3, 0, 05}

\def\markEuler{triangle*}
\def\marksizeEuler{2.5pt}
\def\markGaver{diamond}
\def\marksizeGaver{2.5pt}
\def\markTalbot{o}
\def\marksizeTalbot{2.5pt}
\def\markCME{asterisk}
\def\marksizeCME{2.5pt}
\def\markZakian{pentagon}
\def\marksizeZakian{2.5pt}
\def\markTAME{*}
\def\marksizeTAME{1.8pt}

\usepackage[algoruled]{algorithm2e}

\usepackage{hyperref} 

\definecolor{orange_sidenote}{RGB}{255, 188, 24}   

\usepackage[normalem]{ulem}     

\newtheorem{prop}{Proposition}[section]
\newtheorem{thm}[prop]{Theorem}
\newtheorem{lem}[prop]{Lemma}
\newtheorem{defin}[prop]{Definition}
\newtheorem{cor}[prop]{Corollary}
\theoremstyle{remark}
\newtheorem{rem}[prop]{Remark}

\DeclarePairedDelimiter{\norm}{\lVert}{\rVert}
\DeclarePairedDelimiter{\abs}{\lvert}{\rvert}

\let\Re\relax
\DeclareMathOperator{\Re}{\mathfrak{Re}}		
\let\Im\relax
\DeclareMathOperator{\Im}{\mathfrak{Im}}		

\newcommand{\R}{\mathbb{R}}			
\newcommand{\C}{\mathbb{C}}				

\newcommand{\dd}{\,\mathrm{d}}    
\newcommand{\deltaxx}{\rho}     
\newcommand{\LT}[1]{\widehat{#1}}  
\newcommand{\ones}{\mathbf{1}}

\title{Error analysis of Abate--Whitt methods for Inverse Laplace Transforms and a new algorithm for queuing theory applications}

\author*[1]{\fnm{Nikita} \sur{Deniskin}}\email{nikita.deniskin@sns.it}
\author*[2]{\fnm{Federico} \sur{Poloni}}\email{federico.poloni@unipi.it}

\affil[1]{\orgdiv{Classe di Scienze}, \orgname{Scuola Normale Superiore}, \orgaddress{\city{Pisa}, \country{Italy}}}

\affil[2]{\orgdiv{Dipartimento di Informatica}, \orgname{Università di Pisa}, \orgaddress{\city{Pisa}, \country{Italy}}}

\date{}

\begin{document}

\abstract{
    We study the accuracy of a class of methods to compute the Inverse Laplace Transform, the so-called \emph{Abate--Whitt methods} [Abate, Whitt 2006], which are based on a linear combination of evaluations of $\widehat{f}$ in a few points. We provide error bounds which relate the accuracy of a method to the rational approximation of the exponential function. We specialize our analysis to applications in queuing theory, a field in which Abate--Whitt methods are often used; in particular, we study phase-type distributions and Markov-modulated fluid models (or \emph{fluid queues}). 
    
    We use a recently developed algorithm for rational approximation, the AAA algorithm [Nakatsukasa, Sète, Trefethen 2018], to produce a new family of methods, which we call TAME. The parameters of these methods are constructed depending on a function-specific domain $\Omega$; we provide a quasi-optimal choice for certain families of functions. We discuss numerical issues related to floating-point computation, and we validate our results through numerical experiments which show that the new methods require significantly fewer function evaluations to achieve an accuracy that is comparable (or better) to that of the classical methods.
}



\maketitle
    

\section{Introduction}

The Laplace Transform of a function $f: (0,\infty) \to \mathbb{C}$ is defined as
\begin{equation}
    \LT{f}(s) = \int\limits_0^{\infty} e^{-st} f(t) \dd t.
    \label{eq: Laplace Transform definition}
\end{equation}
We study methods to compute its inverse, i.e., reconstruct the value $f(t)$ in a given point $t$, from evaluations of $\LT{f}$ in several points $\beta_1,\dots,\beta_N\in\mathbb{C}$.

Unlike the better-behaved Fourier transform, inverting the Laplace transform is an ill-posed problem: for instance, there are examples of functions such that $\abs{f(1) - g(1)} = 1$ but $\norm{\LT{f} - \LT{g}}_{\infty} \leq \varepsilon$ for arbitrarily small $\varepsilon$. This property makes numerical computation challenging, since small errors in the computation of $\LT{f}$ can turn into large errors on $f(t)$. Because of its ill-posedness, the Inverse Laplace Transform (ILT) problem has attracted ample attention in the mathematical and numerical community. It would be impossible to cite all relevant works, but we refer the reader to the book~\cite{Schiff1999} for theoretical properties and to the book~\cite{Coh2007} for numerical methods.

In this work, we make two specific assumptions to restrict our interest to special cases:
\begin{itemize}
    \item We are interested in algorithms that compute $f(t)$ in one given point $t$, based on evaluating $\LT{f}$ in a small number of points that we are allowed to choose; these algorithms usually take the form 
    \[
        f(t) \approx \sum_{n=1}^N w_n \LT{f}(\beta_n),
    \] which resembles the general form of a quadrature formula. This setup is known as the \emph{Abate--Whitt framework} in certain fields, after~\cite{AbaWhi2006}. It excludes several algorithms, such as those based on approximating $f$ with a truncated series of functions~\cite[Chapter~3]{Coh2007} or a rational function with known coefficients~\cite[Chapter~5]{Coh2007}, those based on the derivatives of $f$ (e.g., \cite{CuoDRM08} or the Post--Widder formula~\cite[Chapter~7]{Coh2007}), methods where $f(t)$ is evaluated in multiple points simultaneously (e.g., \cite{CuoDMR07,deniseger2006}). We describe this framework and some of the classical methods in Section~\ref{sec:AW}.
    \item We are interested only in the case in which $f$ has a very special (``tame'') form: it belongs to one of certain classes of functions that are constructed based on weighted sums of exponentials. Namely:
    \begin{description}
        \item[SE class (sum of exponentials)] functions of the form
        \[
        f(t) = \sum_{m=1}^M c_m e^{\alpha_m t},
        \]
        where $c_1,\dots,c_m,\alpha_1,\dots,\alpha_m\in\mathbb{C}$ are given constants.
        \item[ME class (matrix exponential)] functions of the form
        \begin{equation} \label{class ME function}
        f(t) = v^* \exp(tQ)u, \quad Q\in\mathbb{C}^{d\times d}, \quad u,v\in\mathbb{C}^d,            
        \end{equation}
        where $\exp(A) = I + A + \frac{1}{2!}A^2 + \frac{1}{3!}A^3\dots$ is the matrix exponential. Using the Jordan form, one sees that these functions can be equivalently written as 
        \[
        f(t) = \sum_{m=1}^M c_m e^{\alpha_m t} t^{b_m},
        \]
        with $b_1,\dots,b_m \in \mathbb{Z}_+$. It is easy to see that this class is the one of all functions $f$ whose transform $\LT{f}$ is a rational function. We chose the acronym ME after the one used for \emph{matrix exponential distributions}~\cite{encystat}, but we note that some authors use the same acronym for other classes named \emph{mixture of exponentials}, e.g.,~\cite[Definition~9.4]{SchSonVon2012}. 
        
        The density of a Matrix Exponential distribution~\cite{encystat} has the form~\eqref{class ME function}, but with additional conditions: $f(t)$ has to be non-negative with total mass 1. In our definition of a ME function, we do not require such conditions. 
        \item[LS class (Laplace--Stieltjes)] functions of the form
        \begin{equation} \label{LS integral}
            f(t) = \int e^{-xt} d\mu(x),    
        \end{equation}
        where $\mu$ is a non-negative measure on $\mathbb{R}^+ = (0,\infty)$. In other words, $f$ in~\eqref{LS integral} is itself the Laplace transform of a measure $\mu$.
        A SE function with negative real weights $\alpha_i$ is in this class: it corresponds to a discrete measure $\mu = \sum_{m=1}^M c_m \delta_{|\alpha_m|}$. The LS class is a generalization with an integral rather than a discrete sum.
        
        Bernstein's Theorem \cite[Theorem 1.4]{SchSonVon2012} states that the class of Laplace-Stieltjes functions is equivalent to the class of \emph{completely monotone} functions, i.e. functions $f(t) \in \mathcal{C}^\infty((0,\infty))$ such that 
\[ (-1)^n f^{(n)}(t) \geq 0 \hspace{5mm} \forall n\geq 0, \, \forall t > 0.\]
    \end{description}
    These three classes are very favorable cases: in most other literature the inverse Laplace transform the main focus is on functions with jump discontinuities (e.g., the Heaviside step function, the square wave) or jumps in the derivative (e.g., formulas involving the absolute-value function $\abs{\cdot}$).
\end{itemize}
Under these two assumptions, we relate the accuracy of an Abate--Whitt method to a rational approximation problem, and prove several quantitative bounds on its error, revealing the role of quantities such as the evaluation point $t$, the range of the $\alpha_i$, and the \emph{field of values} of the matrix $Q$. This is done in Section~\ref{sec:AWaccuracy}. Another approach is based on the moments of a function related to the approximation problem. We present some bounds and an estimate in Section~\ref{sec:moments}.

Despite the strong assumptions that we make, this restricted setup is useful in practice in several applications in queuing theory, which we describe in more detail in Section~\ref{sec:queues}. In particular, we deal with distributions whose density function belongs to a subclass of the ME class, the so-called \emph{phase-type distributions}. They are obtained when $Q$ in~\eqref{class ME function} is the subgenerator of a continuous-time Markov chain and $u= -Q \ones$. Phase-Type distributions predate ME distributions, as they arise naturally in queuing theory as the probability distribution of the hitting time of an absorbing Markov chain. See \cite{BlaNie2017} for an exposition of the theory of ME and Phase-Type distributions. 

Another application in queuing theory is the computation of the so-called \emph{time-dependent first-return matrix} of \emph{Markov-modulated fluid models}, or \emph{fluid queues}, a class of stochastic processes that is used to model continuous-time queues and buffers~\cite{asmussen95,kk95,roger94,ram99}. Indeed, there are several algorithms to compute directly the Laplace transform of this matrix~\cite{BeaOReTay2008}, while the direct computation of the time-domain function is more involved. 

After reducing the problem to rational approximation, we describe the AAA algorithm~\cite{AAA2018}, a recently developed algorithm that produces high-quality rational approximants. With a few modifications, we can use this algorithm to produce the weights and poles $(w_n,\beta_n)_{n=1}^N$ of a family of Abate--Whitt methods, which we dub the TAME method (Triple-A for the Matrix Exponential). We describe the AAA algorithm and its modifications in Section~\ref{sec:AAA}. 

From the computational point of view, one needs to ensure not only that the parameters $(w_n,\beta_n)$ of an Abate--Whitt method produce an accurate rational approximation, but also that the magnitude of weights $w_n$ is moderate.
Weights of large magnitude may lead to precision loss in floating-point arithmetic. We discuss these issues, as well as other numerical aspects of our method, in Section~\ref{sec:TAME}, and
we conclude our paper with numerical experiments that prove the effectiveness of our method, in Section~\ref{sec:experiments}. The experiments show that the TAME method achieves accuracy comparable to the best Abate--Whitt methods, requires significantly fewer evaluations of the function $\LT{f}$, and has better stability properties.

\subsection{Notation}
In the following, $\norm{A}_2$ and $\norm{A}_{\infty}$ denote the Euclidean and $\infty$ operator norm of a matrix $A$. The spectrum of the matrix $A$ is $\Lambda(A)$. For a complex-valued function $f$, we define $\norm{f}_{\infty, \Omega} = \max\{\abs{f(z)} \colon z \in \Omega\}$. We denote with $B(c, R) = \{z\in\mathbb{C} \colon \abs{z-c}\leq R\}$ the ball with center $c$ and radius $R$ in the complex plane. 

\section{Inverse Laplace Transforms and the Abate--Whitt framework} \label{sec:AW}

We assume that the function $f$ is defined (at least) on $\R^+ = (0,\infty)$ and that it has real values. Most of our results are also applicable to functions $f:\R^+ \to \C$, but we shall see that in the real case we can use conjugate poles to speed up the computation of the Inverse Laplace Transform. 

For Equation~\eqref{eq: Laplace Transform definition} to be well-defined, the integral must converge. For the function $f(t) = e^{\alpha t}$, this means that $\LT{f}(s)$ is defined for all $s\in\mathbb{C}$ with $\Re(s-\alpha) > 0$. However, $\LT{f}(s) = \frac{1}{s-\alpha}$ is defined algebraically for all $s\neq \alpha$, irrespective of the convergence of the integral~\eqref{eq: Laplace Transform definition}. This is a form of analytical continuation. When we deal with linear combinations and integral averages of exponential functions, we can use this extended definition of $\LT{f}$: we only need to ensure that $\LT{f}$ (which in most of our paper is a rational function) can be computed in the required points $\frac{\beta_n}{t}$, even if those points do not lie in the domain of convergence of~\eqref{eq: Laplace Transform definition}. In particular, we note that also one of the classical Abate--Whitt methods has poles $\beta_n$ in the left half-plane, namely, the Talbot method (see e.g. Figure~\ref{fig:nodes}); so we deal with points where~\eqref{eq: Laplace Transform definition} does not converge even when using the Talbot method for functions of the form $f(t) = e^{\alpha t}$ with $\alpha < 0$.

A widely used formula for the Inverse Laplace Transform is the Bromwich Integral, which is a contour integral on the vertical line in the complex plane $\gamma(u) = b + iu$ with $u$ ranging from $-\infty$ to $\infty$, namely

\begin{equation}
    f(t) = \frac{1}{2 \pi i} \int\limits_{\gamma} e^{st} \LT{f}(s) \dd s.
    \label{eq: Bromwich Integral}
\end{equation}
As discussed in the Introduction, the Inverse Laplace Transform is an ill-posed problem. We can gain insight on this phenomenon from the above formula~\eqref{eq: Bromwich Integral}: the multiplying weight $e^{st}$ has constant modulus on $\gamma$, so we have to integrate on an unbounded domain a function that may not vanish when $\Im s \to \infty$. This difficulty does not appear in the Direct Laplace Transform, where the multiplying weight is decreasing, nor in the Fourier Transform, where the integral is on a finite domain. 

We focus on a class of methods that have been put in a common structure by Abate and Whitt in a series of works \cite{AbaChoWhi2000, AbaVal2004, AbaWhi2006}. This class includes some well-known classical methods, as well as some more recent ones.
\begin{defin}\label{def: Abate Whitt method}    
	An \emph{Abate--Whitt method} with $N$ weights $(w_n)_{n=1}^{N}$ and $N$ distinct nodes $(\beta_n)_{n=1}^{N}$
    is the formula
    \begin{equation}
    	f_N(t) = \sum\limits_{n=1}^N \frac{w_n}{t} \; \LT{f}\left( \frac{\beta_n}{t}\right).
        \label{eq: Abate Whitt definition}    
    \end{equation}
\end{defin}
\noindent This formula allows one to approximate $f(t) \approx f_N(t)$ with $N$ evaluations of the function $\LT{f}(s)$. In particular, for $t=1$ this formula becomes
\begin{equation}\label{eq: Abate Whitt definition for T = 1}
	f_N(1) = \sum\limits_{n=1}^{N} w_n \; \LT{f}\left( \beta_n \right).
\end{equation} 
Abate--Whitt methods are useful when $\LT{f}$ can be computed (by a black-box function) on a desired set of points $\{\beta_1, \ldots, \beta_n\}$ and we are interested in recovering $f(t)$ for one value of $t$. One typically is interested in the case when computation of $\LT{f}$ is expensive (as in the fluid queue case, Section~\ref{sec:queues}), hence we want to keep $N$ small. Indeed, we will see that in many common cases the TAME method works by evaluating $\LT{f}$ in just 4 points. 

\begin{rem}
Suppose that $f$ satisfies the property that $\LT{f}(\overline{s}) = \overline{\LT{f}(s)}$, where $\overline{s}$ is the complex conjugate of $s$: this form of symmetry with respect to the real axis is common for functions that map $\R^+$ to $\R$. If there is a pair of conjugate nodes and weights, i.e., $(w_n, \beta_n)$ and $(w_{k}, \beta_{k})$ with $w_{k} = \overline{w_n}$, $\beta_{k} = \overline{\beta_n}$, then
\[w_n \LT{f}(\beta_n) + w_{k} \LT{f}(\beta_{k}) = w_n \LT{f}(\beta_n) + \overline{(w_{n} \LT{f}(\beta_{n}))} = 2\Re(w_n \LT{f}(\beta_{n})).\]
Therefore, we can compute the sum with just one evaluation of $\LT{f}$ instead of two. More generally, for a real-valued function $f$, an Abate--Whitt method can be computed in the \emph{reduced form}
\begin{equation}\label{eq: Abate Whitt reduced form}
    f_N(t) = \sum\limits_{n=1}^{N'} \Re \left( \frac{w_n'}{t} \; \LT{f}\left( \frac{\beta_n}{t}\right) \right),  
\end{equation}
where we have reordered the indices such that each weight and node with index in $\{N'+1,\dots,N\}$ is the conjugate of one with index in $\{1,\dots,N'\}$. To keep the sum unchanged, we set $w_n'=2w_n$ if $n$ is part of a conjugate pair, and $w_n'=w_n$ otherwise.
\end{rem}

The Abate--Whitt methods that we present below use this speed-up trick and are already in reduced form. All their weights and nodes are in conjugate pairs, plus possibly an unpaired one when $N$ is odd, so $N'$ is either $N/2$ or $(N+1)/2$. Choosing conjugate pairs speeds up the computation by a factor 2, and ensures that the computed $f_N(t)$ is real. 

There are many ways to construct an Abate--Whitt method. For example, a quadrature of the Bromwich Integral \eqref{eq: Bromwich Integral} is a weighted sum of evaluations of $\LT{f}$, so it has the same form as Equation~\eqref{eq: Abate Whitt definition for T = 1}. The Euler method \cite{AbaChoWhi2000, AbaWhi2006} uses equispaced nodes $\beta_n$ on the vertical line $\gamma$. Abate and Whitt \cite[Section 5]{AbaWhi2006} propose the following choice for the parameters, with odd $N'$ and for $1 \leq n \leq N'$
\begin{equation}
	\beta_n =  \frac{\ln(10)}{6} (N'-1) + i \pi (n-1) \; \; \text{ and } \; \;  w_n' =  10^\frac{N'-1}{6} (-1)^n \xi_n, 
\end{equation}
where
\begin{equation*}
	\xi_1 = \frac{1}{2}, \; \;
	\xi_n = 1 \;\; \forall\, n : \; 2 \leq n \leq \frac{N'+1}{2},
\end{equation*}
\begin{equation*}
	\xi_{\frac{N'+3}{2}+j} = 1 -  2^{-\frac{N'-1}{2}} \left(\sum\limits_{k=0}^j \binom{\frac{N'-1}{2}}{k} \right) \;\;\; \forall\, j : \; 0 \leq j \leq \frac{N'-3}{2}.
\end{equation*}
Talbot \cite{Tal1979} proposed to deform the contour, to avoid the oscillation of $e^{st}$ due to its imaginary part. The new contour is the curve 
$\gamma(\theta) = r \theta \left( \cot(\theta) +a i\right)$,  with parameters $r>0, a>0$ and $-\pi < \theta < \pi$,
which has its imaginary part bounded in the interval $[-\pi ra, \pi ra]$. Abate and Valk\'o \cite{AbaVal2004} describe a way to fit the Talbot method in the Abate--Whitt framework using this contour. The parameters are \cite[Section 6]{AbaWhi2006} 
\begin{equation*}
	\begin{aligned}
	\beta_1 &= \frac{2 N'}{5},  
	 && \beta_n= \frac{2 (n-1) \pi }{5} \left( \cot(\theta_n)  + i   \right)  &&\text{ for } 2 \leq n \leq N', \\
	w_1' &= \frac{1}{5} e^{\beta_1},
	&& w_n' = \frac{2 }{5} \left(  1+ i \theta_n \left( 1 + \cot\left( \theta_n \right)^2 \right) 
	- i \cot\left( \theta_n \right) \right) e^{\beta_n}   &&\text{ for } 2 \leq n \leq N'.
	\end{aligned}
\end{equation*}
Further optimization of the contour was also investigated in \cite{Weideman2006}. 

The Gaver--Stehfest method \cite{Gav1996,Stehfest1970} is based instead on the Post-Widder formula \cite[Theorem 2.4]{Coh2007}. Its parameters are (\cite[Section 4]{AbaWhi2006}, \cite{HHAT2020})
\begin{equation*}
	\beta_n = n \ln(2), \hspace{2mm} 1 \leq n \leq N',
\end{equation*}
\begin{equation*}
	w_n' = (-1)^{N'/2+n} \ln(2) \sum\limits_{j=\lfloor (n+1)/2\rfloor }^{\min(n,N'/2)} \frac{j^{N'/2+1}}{(N'/2) !} \binom{N'/2}{j} \binom{2j}{j} \binom{j}{n-j}, \hspace{2mm} 1 \leq n \leq N'.
\end{equation*}
In this method all poles and weights are real, but the drawback is that it is prone to numerical cancellation, since it  computes implicitly higher-order derivatives through finite differences. 

%

These three methods are described extensively  by Abate and Whitt in \cite{AbaWhi2006}. They mention also the work of Zakian \cite{Zak1969, Zak1970} and that his method can be included in the framework, but they do not focus on it. 

The \emph{Concentrated Matrix Exponential} (CME) method has been introduced by Telek and collaborators in a series of papers \cite{HHAT2020, HorHorTel2020, HSTZ2016}. Both Zakian's and the CME method are based on the following idea. Let $(w_n,\beta_n)_{n=1}^N$ be the parameters of an Abate--Whitt method. Substitute the definition of $\LT{f}$ into the defining formula~\eqref{eq: Abate Whitt definition} of the Abate--Whitt method
\begin{equation}
    \begin{aligned}
    f_N(t) & = \sum\limits_{n=1}^N \frac{w_n}{t} \LT{f}\left( \frac{\beta_n}{t} \right) 
    = \sum\limits_{n=1}^N \frac{w_n}{t}  \int\limits_0^\infty e^{-\frac{\beta_n}{t} x} f(x) \dd x \\
    & = \int\limits_0^\infty f(x) \left( \sum\limits_{n=1}^N w_n e^{-\beta_n \frac{x}{t}}  \right) \, \frac{1}{t} \dd x \\
     & = \int\limits_0^\infty f(ty) \left( \sum\limits_{n=1}^N w_n e^{-\beta_n y}  \right) \dd y.    
\end{aligned}
\label{eq: Abate--Whitt through Dirac approximant}
\end{equation}
\begin{defin}
    The \emph{Dirac approximant} of the Abate--Whitt method $(w_n,\beta_n)_{n=1}^N$ is the function
    \[
        \deltaxx_N(y) = \sum\limits_{n=1}^N w_n e^{-\beta_n y}.
    \]
    \label{def: Dirac approximant}
\end{defin}
If in~\eqref{eq: Abate--Whitt through Dirac approximant} we replace $\deltaxx_N(y)$ with the Dirac delta $\delta_1(y)$ at point $y=1$, we recover $f(t)$ exactly:
\begin{equation*}
    f(t) = \int\limits_0^\infty f(ty) \, \delta_1(y) \dd y.   
\end{equation*}
We can create an effective Abate--Whitt method by choosing the parameters $(w_n,\beta_n)_{n=1}^N$ in such a way that the Dirac approximant $\deltaxx_{N}(y)$ is ``close" to the Dirac delta $\delta_1(y)$ in a suitable sense. The two methods use different approaches.

For the CME method \cite{HHAT2020}, the authors construct a bell-shaped Dirac approximant; that is, a (continuous) function $\deltaxx_N(y)$ that assumes a large value at $y=1$,  with most of the mass concentrated near $y=1$ (hence the name \emph{Concentrated Matrix Exponential}), and that tends rapidly to zero away from $1$. The authors choose the nodes $\beta_n$ equispaced on the vertical line (just as in the Euler method), and find the appropriate weights $w_n$ by an optimization procedure, aiming to minimize the normalized variance of the function $\deltaxx_N(y)$ (see Definition~\ref{def:SCV}). 


Zakian \cite{Zak1969, Zak1970} instead takes Laplace Transforms
\[\LT{\delta_1}(s) = \int\limits_{0}^\infty e^{-sy} \delta_1(y) \dd y = e^{-s},\]
and
\[\LT{\deltaxx_N}(s) = \int\limits_{0}^\infty e^{-sy} \left(\sum\limits_{n=1}^N w_n e^{-\beta_n y} \right)\dd y = \sum\limits_{n=1}^N \int\limits_{0}^\infty w_n e^{-(s+\beta_n) y} \dd y = \sum\limits_{n=1}^N \frac{w_n}{\beta_n+s}\]
 and  chooses the parameters to make the two transformed functions close. 

If we set $z=-s$, this becomes a classical problem: approximating the exponential $e^{z}$ with a rational function $\LT{\deltaxx_N}(-z)$. 
\begin{defin}
The \emph{rational approximant} of the Abate--Whitt method $(w_n,\beta_n)_{n=1}^N$ is 
\[\LT{\deltaxx_N}(-z) = \sum\limits_{n=1}^N \frac{w_n}{\beta_n-z}.\]
\end{defin}
Zakian suggests choosing $\LT{\deltaxx_N}(-z)$ to be an accurate rational approximation of the exponential $e^z$, to produce an accurate Abate--Whitt method. However, recall that the Inverse Laplace Transform is ill-posed: small perturbations to $\LT{g}$ can result in big perturbations to $g$. Thus, even if $\norm{\LT{\deltaxx_N}(-z) - \LT{\delta_1}(-z)}$ is small in a suitable norm, this does not mean that $\deltaxx_N(y) - \delta_1(y)$ is small as well. Indeed, we can not even use a classical norm to measure this error, since $\delta_1(y)$ is a distribution and not a function in the classical sense. In the next section, we make this approach more rigorous.

Zakian suggests using the $(N-1, N)$th Padé approximant of the exponential $e^{z}$ as the rational approximant $\LT{\deltaxx_N}(-z)$. With this choice, $\LT{\deltaxx_N}(-z)$ is an excellent approximation of $e^{z}$ in a neighbourhood of $z=0$, but gets progressively worse as $|z|$ grows. As noted in \cite{Zak1970}, the method is exact when $f$ is a polynomial of degree at most $2N-1$. In our experiments we found that Zakian method exhibits fast convergence, although it suffers from numerical instability, evident from $N'=5$ onward.

The use of Padé approximants is also discussed by Wellekens \cite{Wellekens1970}, refining an approach due to Vlach \cite{Vlach1969} (see also \cite{AbaWhi2006}). The idea is to approximate the exponential in \eqref{eq: Bromwich Integral} with a rational function, whose poles and residues are the nodes and weights of the Abate--Whitt method. While this structure is similar to the method we propose later, they differ in a key aspect: Wellekens approximates the exponential on the contour of the Bromwich integral, while the approximation domain $\Omega$ of the TAME method is tailored to the function $f$. 

\section{Accuracy of AW methods through rational approximants} \label{sec:AWaccuracy}
In this section we provide theoretical background for Zakian's approach, proving that an accurate rational approximation of $\exp(z)$ gives an accurate Abate--Whitt method. The relation between rational approximation and accuracy of a method for inverse Laplace transform is not new: a discussion appears for instance in~\cite{TreWeiSch}, though formulated in terms of integrals.
\begin{defin}
    We say that an Abate--Whitt method $(w_n,\beta_n)_{n=1}^N$ is \emph{$\varepsilon$-accurate on $\Omega\subseteq \mathbb{C}$} if 
    \begin{equation} \label{epsilon-accurateness condition}
    \norm*{\exp(z) - \LT{\deltaxx_N}(-z)}_{\infty,\Omega} = \norm*{\exp(z) - \sum_{n=1}^N \frac{w_n}{\beta_n - z}}_{\infty,\Omega} \leq \varepsilon.
    \end{equation}
\end{defin}

Based on this definition, we obtain accuracy results for functions of class SE and ME; we start with the SE case since the proof is simpler.
\begin{thm}\label{thm: approximation error SE}
    Let $f(t)$ be a SE function
    \[f(t) = \sum\limits_{m=1}^M c_m e^{\alpha_m t}.\]
    Let $(w_n, \beta_n)_{n=1}^N$ be an Abate--Whitt method that is $\varepsilon$-accurate on a region $\Omega$ which contains $\alpha_1t, \dots, \alpha_Mt$. 
    
    Then, the error of the method at point $t$ is bounded by
    \[\left| f(t) - f_N(t) \right| \leq \left(  \sum\limits_{m=1}^M |c_m|\right) \cdot \varepsilon.\]
\end{thm}
\begin{proof}
    Since
    \[
    \LT{f}(s) = \sum_{m=1}^M \frac{c_m}{s-\alpha_m},
    \]
    we have
    \begin{align*}
        \abs{f(t) - f_N(t)} &= \abs*{\sum_{m=1}^M c_m e^{\alpha_m t} - \sum_{n=1}^N \frac{w_n}{t} \left(\sum_{m=1}^M\frac{c_m}{\frac{\beta_n}{t}-\alpha_m} \right)}\\
        &= \abs*{\sum_{m=1}^M c_m \left(e^{\alpha_m t} - \sum_{n=1}^N \frac{w_n}{t \left(\frac{\beta_n}{t} - \alpha_m\right)}\right)}\\
        &= \abs*{\sum_{m=1}^M c_m \left(e^{\alpha_m t} - \sum_{n=1}^N \frac{w_n}{\beta_n - \alpha_m t}\right)}\\
        & \leq \sum_{m=1}^M \abs{c_m}\, \varepsilon.
    \end{align*}
    In the last line, we have used the $\varepsilon$-accurateness condition~\eqref{epsilon-accurateness condition} in the points $\alpha_m t$.
\end{proof}
\begin{rem} \label{rem:negative result}
    When $M=1$, the inequality in this theorem becomes an equality. So we have also a converse negative result: let $z_*\in\mathbb{C}$ be a point for which $\abs{\exp(z_*) - \LT{\deltaxx_N}(-z_*)} = C > \varepsilon$; then, if we choose $\alpha,t_*$ such that $\alpha t_* = z_*$ we have that $|f(t_*) - f_N(t_*)| = C > \varepsilon$ for the function $f(t) = e^{\alpha t}$. Or, informally: given a point $z_*$ in which the rational approximant of an Abate--Whitt method is inaccurate, we can find an exponential function $f$ and a point $t_*$ for which the method is inaccurate. Since no rational approximation of the exponential can be accurate on the whole complex plane, this means that no Abate--Whitt method can be accurate for all exponential functions and all points $t$.
\end{rem}

Therefore, the quality of the Inverse Laplace Transform with an Abate--Whitt method depends on the approximation error of $e^{z}$ with the rational function $\LT{\deltaxx_N}(-z)$ at points $\{\alpha_m t\}_{m=1}^M$. Given $f$, the coefficients $c_m$ are fixed, but we can increase the order $N$ and choose the parameters $(w_n,\beta_n)_{n=1}^N$ opportunely to obtain rational approximations of $e^{z}$ which get better on the points $\{\alpha_m t\}_{m=1}^M$ as $N$ increases. This observation ensures that a suitable family of Abate--Whitt approximations can achieve convergence when $N\to\infty$ (at least in exact arithmetic).

To extend the result to the class of ME functions, we need a few technical tools.
\begin{defin}
    The \emph{field of values}, also called \emph{numerical range}, of a matrix $A\in\mathbb{C}^{d\times d}$ is the complex set
    \[
    \mathbb{W}(A) = \{\mathbf{x}^*A\mathbf{x} \colon \mathbf{x}\in\mathbb{C}^d,\norm{\mathbf{x}}=1\}.
    \]
\end{defin}
The field of values is a well-known tool in linear algebra; here we recall only a few classical results (see, e.g., \cite[Section~I.2]{BhatiaMA} or \cite{FOVbook}).
\begin{lem} \label{lem:Wproperties}The following properties hold.
    \begin{enumerate}
        \item We have the inclusions
        \[
        \operatorname{hull}(\Lambda(A)) \subseteq \mathbb{W}(A) \subseteq B(0, \norm{A}_2),
        \]
        where $\operatorname{hull}(\Lambda(A))$ is the convex hull of the eigenvalues of $A$. The left inclusion is an equality when $A$ is a normal matrix.
        \item Translation and rescaling of a matrix changes its field of values in the same way: 
        \[\mathbb{W}(\alpha A + \beta I) = \{\alpha z + \beta \colon z \in \mathbb{W}(A)\}.\]
        \item If $B$ is a principal submatrix of $A$, then $\mathbb{W}(B) \subseteq \mathbb{W}(A)$.
        \item Let $J_0\in\mathbb{R}^{d\times d}$ be the size-$d$ Jordan block with eigenvalue $0$. We have $\mathbb{W}(J) = B(0, \cos\frac{\pi}{d+1})\subseteq B(0,1)$.
    \end{enumerate}
\end{lem}
An important result relating the field of values and matrix functions is the following.
\begin{thm}[\protect{Crouzeix-Palencia, \cite{CroPal17}}] \label{thm:CrouzeixPalencia}
    Let $A$ be a square matrix, and let $\phi(z)$ be a holomorphic function on $\mathbb{W}(A)$. Then,
    \[
    \norm{\phi(A)}_2 \leq (1+\sqrt{2}) \, \norm{\phi}_{\infty,\mathbb{W}(A)},
    \]
    where $\phi(A)$ denotes the extension of $\phi$ to square matrices.
\end{thm}
It is conjectured that the constant $1+\sqrt{2}$ can be replaced by $2$ (Crouzeix's conjecture).

With these tools, we extend our error bound to matrices. 
\begin{thm} \label{thm:class2approximation}
    Let $Q$ be a $d \times d$ matrix and $f(t) = \exp(tQ)$. 
    Let $(w_n, \beta_n)_{n=1}^N$ be an Abate--Whitt method that is $\varepsilon$-accurate on a region $\Omega$ which contains $\mathbb{W}(tQ)$.
    
    Then, the error of the ILT at point $t$ is bounded by
    \[
        \norm{f(t)-f_N(t)}_2 \leq (1+\sqrt{2}) \varepsilon.
    \]    
\end{thm}
\begin{proof}
    The Laplace Transform of $f$ is
    \[
        \LT{f}(s) = \left( sI - Q\right)^{-1}.
    \]
    This can be proven by using the spectral representation of a matrix function \cite[Section 7.9]{Meyer2000}. The Laplace Transform is well-defined when $\Re(s) > \Re(\Lambda(Q))$, but as discussed in Section \ref{sec:AW}, it can be extended to $s \not\in\Lambda(Q)$. The Abate--Whitt approximant of $f$ is 
    \[
        f_N(t) = \sum\limits_{n=1}^N \frac{w_n}{t}\, \LT{f}\left(\frac{\beta_n}{t}\right) = \sum\limits_{n=1}^N \frac{w_n}{t} \, \left(\frac{\beta_n}{t}I - Q\right)^{-1} =\, \sum\limits_{n=1}^N w_n \, \left(\beta_n I - tQ \right)^{-1} . 
    \]
    We apply the Crouzeix-Palencia Theorem~\ref{thm:CrouzeixPalencia} to the function
    \[
    \phi(z) = e^{z} - \sum_{n=1}^N \frac{w_n}{\beta_n - z},
    \]
    obtaining
    \[ 
    \norm*{\exp(tQ) - \sum_{n=1}^N w_n\, (\beta_n - tQ)^{-1}}_{2} \leq (1+\sqrt{2})  \varepsilon. \qedhere
    \]
\end{proof}
Now we can obtain any ME function by choosing an appropriate matrix $Q$ and doing a left and right vector product.  
\begin{cor} \label{cor:classME}
    Let $v,u\in\mathbb{C}^d$ and $Q\in\mathbb{C}^{d\times d}$. Let $f(t) = v^* \exp(tQ)u$. Let $(w_n, \beta_n)_{n=1}^N$ be an Abate--Whitt method that is $\varepsilon$-accurate on a region $\Omega$ which contains $\mathbb{W}(tQ)$.    
    Then,
    \[
    \abs{f(t) - f_N(t)} \leq (1+\sqrt{2})\varepsilon \norm{v}_2\norm{u}_2.
    \]
\end{cor}
\begin{cor}
    Let $f(t) = \frac{t^b}{b!} e^{\alpha t}$, with $b\in\mathbb{N}$ and $\alpha\in\mathbb{C}$. 
    Let $(w_n, \beta_n)_{n=1}^N$ be an Abate--Whitt method that is $\varepsilon$-accurate on a region $\Omega$ which contains $B(\alpha t,t)$.
    Then,
    \[
    \abs{f(t) - f_N(t)} \leq (1+\sqrt{2}) \varepsilon.
    \]
\end{cor}
\begin{proof}
    Let $J$ be the $(b+1)\times(b+1)$ Jordan block with eigenvalue $\alpha$; then we have
    \begin{gather*}
    J = \begin{bmatrix}
        \alpha & 1\\
        & \alpha & 1 \\
        & & \ddots & \ddots\\
        & & & \alpha & 1\\
        & & & & \alpha
    \end{bmatrix},\quad
    \exp(tJ) = e^{\alpha t}\begin{bmatrix}
         1 & t & \frac{t^2}{2} & \dots & \frac{t^b}{b!}\\
        & 1 & t & \ddots & \vdots\\
        & & 1 & \ddots & \frac{t^2}{2}\\
        & & & \ddots & t\\
        & & & & 1
    \end{bmatrix}
    \end{gather*}
    and $\mathbb{W}(tJ) \subset B(\alpha t,t)$ thanks to the properties in Lemma~\ref{lem:Wproperties}. The result now follows from Corollary~\ref{cor:classME}, by setting $v=e_1, u=e_{b+1}$, and noting that $f(t) = v^* \exp(tJ) u$.
\end{proof}

If a function $f$ is close to a function $g$ for which we know the error of an Abate--Whitt method (for example, if $g$ is in the SE class), then we can bound the error of the ILT with the Abate--Whitt method for $f$. 

\begin{thm} \label{thm: f approximated by g bound}
    Let $f(t)$ be a function, and suppose that $g(t)$ is an approximation of $f$ with error $\eta$ on $(0,\infty)$. 
    \[
        \left\Vert f(t) -  g(t) \right\Vert_{\infty,\, \R^+} \leq \eta.
    \]
    Let $(w_n, \beta_n)_{n=1}^N$ be an Abate--Whitt method.  Let $\deltaxx_N$ be the Dirac approximant of the method.
    
    Then, the error of the method at point $t$ is bounded by
    \[
     \abs*{f(t)- f_N(t)} \leq (1+\norm{\deltaxx_N}_1) \, \eta  +   \abs*{g(t)- g_N(t)}
    \]
\end{thm}
\begin{proof}
\[
\begin{aligned}
    &\phantom{=} f(t) - f_N(t) = f(t) - \int\limits_0^\infty  f(ty) \, \deltaxx_N(y)  \dd y \\
    &= \int\limits_0^\infty  (g(ty)-f(ty)) \, \deltaxx_N(y) \dd y  + g(t) - \int\limits_0^\infty  g(ty) \, \deltaxx_N(y) \dd x + (f(t)-g(t)).
\end{aligned}
\]
We estimate the first term as 
\[ 
\begin{aligned}
    \abs*{\int\limits_0^\infty  (g(ty)-f(ty)) \, \deltaxx_N(y) \dd y } &\leq  
    \int\limits_0^\infty  \abs*{g(ty)-f(ty)} \, \abs*{\deltaxx_N(y)} \dd y \\
    & \leq \max\limits_{y\in(0,\infty)} \abs*{g(ty)-f(ty)} \cdot \int\limits_0^\infty  \abs*{\deltaxx_N(y)} \dd y \\
    & = \norm*{g-f}_{{\infty},\R^+} \cdot  \norm*{\deltaxx_N}_1 \\
    & \leq \eta \, \norm*{\deltaxx_N}_1 .
\end{aligned}
\]
The second term is the Abate--Whitt approximation error $g(t) - g_N(t)$, while the third term is bounded as $\abs*{f(t) - g(t)} \leq \eta$.
Putting everything together, we obtain 
\[
    \abs*{f(t)- f_N(t)} \leq \eta \, (1+\norm{\deltaxx_N}_1)  + \abs*{g(t)-g_N(t)}. \qedhere
\]
\end{proof}
Combining Theorem~\ref{thm: f approximated by g bound} and Theorem~\ref{thm: approximation error SE} we get the following result. 
\begin{thm} \label{thm: f approximated by SE}
    Let $f(t)$ be a function, and suppose that $g(t) = \sum\limits_{m=1}^M c_m e^{\alpha_m t}$ is an approximation of $f$ with error $\eta$. 
    \[
        \left\Vert f(t) -  \sum\limits_{m=1}^M c_m e^{\alpha_m t} \right\Vert_{{\infty}, \R^+} \leq \eta
    \]
    Let $(w_n, \beta_n)_{n=1}^N$ be an Abate--Whitt method that is $\varepsilon$-accurate on a region $\Omega$ which contains $\alpha_1t, \dots, \alpha_Mt$. Let $\deltaxx_N$ be the Dirac approximant of the method. Then, the error of the method at point $t$ is bounded by
    \begin{equation} \label{eq:AW bound for SE approximation}
     \abs*{f(t)- f_N(t)} \leq (1+\norm{\deltaxx_N}_1) \, \eta  +  \left(\sum\limits_{m=1}^M |c_m|\right) \, \varepsilon. 
    \end{equation}
\end{thm}
Let us comment on this result. It is not immediate how to choose parameters that make the right-hand side of~\eqref{eq:AW bound for SE approximation} small. To approximate $f(t)$ with a small error $\eta$, we may have to choose a SE function with large weights $c_m$. Symmetrically, to approximate $\LT{\deltaxx_N}(s)$ with a small error $\varepsilon$, we may have to choose an Abate--Whitt method with large $\norm*{\deltaxx_N}_1$. But if $\deltaxx_N \approx \delta_1$, we can expect $\norm*{\deltaxx_N}_1$ not to be much larger than $\norm{\delta_1}_1 = 1$; say, $\norm*{\deltaxx_N}_1 \leq 10$. If this bound holds uniformly for a family of Abate--Whitt methods, we can first choose $\eta$ to make the first summand in~\eqref{eq:AW bound for SE approximation} small, and then select an Abate--Whitt method in the family to make $\varepsilon$ (and hence the second summand in~\eqref{eq:AW bound for SE approximation}) small.

We use Theorem~\ref{thm: f approximated by SE} to obtain a bound for a LS-class function with \emph{finite} measure $\mu$, i.e., one such that $\mu(\mathbb{R}^+) < \infty$. Note that in this case then $f(0) = \int_0^\infty \dd \mu(x) = \mu(\mathbb{R}^+)$ is also defined using the formula~\eqref{LS integral} and finite; hence an alternative characterization is LS functions that do not diverge in $0$.
\begin{prop} \label{prop: LS accuracy bound}
Let $f$ be a LS function with a finite measure $\mu$. Then, for each $\varepsilon>0$ there exists a SE-class function $g(t) = \sum_{m=1}^M c_m e^{-x_m t}$ such that $\norm{f-g}_{\infty,\R^+} \leq \varepsilon$. We can take $g(t)$ to have $c_m > 0$, $\sum_{m=1}^M c_m = \mu(\mathbb{R}^+)$, and $x_1,\dots,x_M \in (0,L]$, where $L>0$ is such that
\[
\frac{\mu((L,\infty))}{\mu(\R^+)} \leq \varepsilon.
\]
\end{prop}
\begin{proof}
    We assume, up to scaling, that $\mu$ is a probability measure, i.e., $\mu(\mathbb{R}^+)=1$. Let $F(x)$ be its cumulative distribution function (CDF). If $X$ is a random variable with distribution $F(x)$, then the distribution of the clamped random variable $\min(X,L)$ is
    \[
    F_L(x) = \begin{cases}
        F(x) & x < L,\\
        1 & x \geq L.
    \end{cases}
    \]
    By the Glivenko-Cantelli theorem (uniform strong law of large numbers, \cite{encystat}), if the reals $x_1, x_2, \dots$ are taken sampled from the distribution $F_L$, then
    the CDF $G_M(t)$ of the measure $\frac{1}{M}\sum_{m=1}^M \delta_{x_m}$ satisfies almost surely $\lim_{M\to\infty }\|F_L - G_M\|_{\infty,\R^+} = 0$. In particular, this implies that for each $\varepsilon$ there exist choices of $x_1,\dots, x_M \in (0,L]$ such that $\|F_L - G_M\|_{\infty,\R^+} \leq \varepsilon$. Then $\|F - G_M\|_{\infty,\R^+}\leq \varepsilon$ holds too, since for any $x>L$ we have $G_M(x)=1$ and $F(x) \geq 1 - \varepsilon$.

    We set $g(t) = \int_0^{\infty} e^{-xt} \dd G(x) = \frac{1}{M}\sum_{m=1}^M e^{-tx_m}$, a SE function. Integrating by parts, we have
    \begin{align*}
    \abs{f(t) - g(t)} &= \abs*{\int_{0}^{\infty} e^{-xt}(\dd F(x)- \dd G_M(x))}\\
    & = \abs*{\int_{0}^{\infty} \left(\frac{\dd }{\dd x}e^{-xt}\right)(F(x)- G_M(x)) \dd x} \\
    &\leq \int \abs*{\frac{\dd }{\dd x}e^{-xt}} \dd x \cdot \varepsilon = \varepsilon.
    \end{align*}
    To justify rigorously the use of integration by parts even when the measures are not Lebesgue-continuous, we can use for instance~\cite[Theorem~18.4]{billingsley}.
\end{proof}

 The convergence speed (as $N$ grows) of the approximation of LS functions with SE functions is also studied in detail in~\cite{Hackbush}.

Combining Theorem~\ref{thm: f approximated by SE} and Proposition~\ref{prop: LS accuracy bound}, we obtain the following result.
\begin{thm} \label{thm:error of LS class}
    Let $f$ be a LS function with a finite measure $\mu$. Let $L>0$ and $\eta>0$ be such that $\mu((L,\infty)) \leq \eta \mu(\mathbb{R}^+)$. Let $(w_n, \beta_n)_{n=1}^N$ be an Abate--Whitt method that is $\varepsilon$-accurate on $\Omega=[-L,0)$. Let $\deltaxx_N$ be the Dirac approximant of the method. Then, the error of the method at point $t$ is bounded by
    \[
     \abs*{f(t)- f_N(t)} \leq (1+\norm{\deltaxx_N}_1) \, \eta  +  \mu(\mathbb{R}^+) \, \varepsilon.
    \] 
    Moreover, if the Abate--Whitt method is $\varepsilon$-accurate on the whole half-line $\mathbb{R}^+$, we can let $L\to\infty$ and $\eta \to 0$, obtaining
    \[
     \abs*{f(t)- f_N(t)} \leq \mu(\mathbb{R}^+) \, \varepsilon.
    \]
\end{thm}

\begin{rem}
    Our previous theorems are valid for a function $f(t)$ in the SE or ME class that can take, in the most general case, complex values. Functions of the LS class are positive (since they are the integral of a positive function with a positive measure). However, the Laplace Transform, its Inverse, and the Abate--Whitt methods are linear. If we can recover both functions $g_1$ and $g_2$ with a small error, we can also recover any linear combination $c_1g_1+c_2g_2$. Therefore, if a function $f$ can be written as $f(t) = g_1(t) - g_2(t)$ where $g_1, g_2$ are LS functions with measures $\mu_1$ and $\mu_2$, then Theorem~\ref{thm:error of LS class} is valid also for $f$, with $\mu_1(\mathbb{R}^+)+\mu_2(\mathbb{R}^+)$ in place of $\mu(\mathbb{R}^+)$.
\end{rem}

\section{Bounds based on moments}
\label{sec:moments}
We now present bounds (and an estimate) based on the moments of a Dirac approximant $\deltaxx_N$, which are valid under small assumptions on the regularity of $f$.
\begin{defin} \label{def:moments}
	Let $\rho: \R^+ \to \R$ be a function. The \emph{moments} of $\rho$ are
	\begin{equation*}
		\mu_0 = \int\limits_0^\infty \rho(y) \dd y, \hspace{5mm} 
		\mu_1 = \int\limits_0^\infty y \rho(y) \dd y, \hspace{5mm} 
		\mu_2 = \int\limits_0^\infty y^2 \rho(y) \dd y.\\
	\end{equation*}
    The \emph{shifted moments} are
	\begin{equation*}
		\nu_2 = \int\limits_0^\infty (y-1)^2 \rho(y) \dd y, \hspace{5mm}
		\widetilde{\nu}_2 = \int\limits_0^\infty (y-1)^2 \, |\rho(y)| \, \dd y.		
	\end{equation*}    
\end{defin}
The shifted moment $\nu_2$ can be expressed through $\mu_0, \mu_1, \mu_2$ as
\[\nu_2  = \int\limits_0^\infty (y^2-2y+1) \rho(y) \dd y = \mu_2 - 2 \mu_1 + \mu_0.\]
If $\rho$ is non-negative, then $\nu_2(\rho) = \widetilde{\nu}_2(\rho)$.

It is a classical result that the moments of a function (when they exist) can be expressed through the derivatives of its Laplace transform in $0$:
\begin{equation} \label{moments and Laplace transforms}
\mu_0 = \LT{\rho}(0), \quad -\mu_1 = \frac{\dd}{\dd s} \LT{\rho}(s)\bigg|_{s=0}, \quad \mu_2 = \frac{\dd^2}{\dd s^2} \LT{\rho}(s)\bigg|_{s=0};    
\end{equation}
see, e.g.,\cite[Section~21]{billingsley}, which shows this fact more generally for probability distributions.

Another quantity related to moments appear prominently in \cite{HorHorTel2020}.
\begin{defin} \label{def:SCV}
    Let $\rho:\R^+ \to \R^+$ be a function with finite moments $\mu_0,\mu_1,\mu_2$. The \emph{Squared Coefficient of Variation (SCV)} is
    \begin{equation*}
        \mathrm{SCV}(\rho) = \frac{\mu_0 \mu_2}{\mu_1^2}-1.
    \end{equation*}
\end{defin}
If $\rho$ is the pdf of a random variable $X$, then the SCV is the normalized variance of $X$: $\mathrm{SCV}(\rho) = \frac{\mathrm{Var}(X)}{\mathbb{E}(X)^2}$. In general, the relation between the SCV and the second shifted moment is
\begin{equation*}
	\mathrm{SCV}(\rho) =  \frac{\mu_0 \mu_2 - \mu_1^2}{\mu_1^2}  =  \frac{\mu_0 (\nu_2 +2\mu_1 - \mu_0)- \mu_1^2}{\mu_1^2}   = \nu_2 \, \frac{\mu_0}{\mu_1^2} - \left(1-\frac{\mu_0}{\mu_1} \right)^2.
\end{equation*}
Note that if $\mu_0(\rho)=\mu_1(\rho) = 1$, then $\mathrm{SCV}(\rho) = \nu_2$. In \cite{HHAT2020}, the CME method is computed with an optimization procedure that minimizes $\mathrm{SCV}(\deltaxx_N)$, and the following bound is obtained. That article assumes $\mu_0=1$ for the definition of the SCV and for Theorem~\ref{thm: HHAT2020 thm 4, Lipschitz}. However, the result can be easily extended to a generic non-negative function $\rho$ by rescaling it as $\frac{1}{\mu_0(\rho)} \rho$.
\begin{thm}[{\cite[Theorem 4]{HHAT2020}}] \label{thm: HHAT2020 thm 4, Lipschitz}
	Let $\deltaxx(x): \R^+ \to \R^+$ be a non-negative function with finite moments $\mu_0,\mu_1,\mu_2$, and assume that $\mu_0 = \mu_1 = 1$. Let $f: \R^+ \to \R$ be bounded by a constant $H$, and Lipschitz-continuous with constant $L$ at point $t$, i.e.
	\begin{equation*}
		|f(t)| \leq H 
		\hspace{3mm} \text{and} \hspace{3mm} 
		|f(t) - f(t_1) | \leq L |t-t_1| 
		\hspace{3mm} \text{for }\,  t_1 \geq 0.
	\end{equation*}
	If $f_N(t) = \int\limits_0^\infty f(tx) \deltaxx(x) \dd x$, then the error of this approximation is bounded by
	\vspace*{-3mm}
	\begin{equation*}
		|f_N(t) - f(t) | \, \leq\, 3\left(2 H L^2 t^2  \right)^{\frac13} \, (\mathrm{SCV}(\deltaxx))^{\frac13}.
	\end{equation*}
\end{thm}

We are interested in giving similar bounds for functions $\deltaxx_N$ which are approximations of the Dirac delta distribution in 1. As the Dirac delta is a probability distribution with mean 1, we expect to have $\mu_0(\deltaxx_N) \approx 1$ and $\mu_1(\deltaxx_N) \approx 1$. For the Dirac approximants of the Abate--Whitt methods these are not exact equalities, so we give bounds in which $\mu_0(\deltaxx_N)$ and $\mu_1(\deltaxx_N)$ can differ from 1. 
\begin{thm}
	Let $\rho(y): \R^+ \to \R$ be a function, with moments as in Definition~\ref{def:moments}. Let $f: \R^+\to \R$ be a $\mathcal{C}^2$ function and $f_N(t) = \int_0^\infty f(ty) \rho(y) \dd y$. Then
		\begin{equation*}
			|f_N(t) - f(t)|  \leq \left|\mu_0 -1 \right| \, |f(t)| +  t \, |f'(t)| \, |\mu_1 - \mu_0| + \frac{1}{2} \, t^2 \; \Vert f'' \Vert_{\infty,\R^+} \; \widetilde{\nu}_2 .
		\end{equation*}
\end{thm}
\begin{proof}
We have $\int_0^\infty \rho(y) \dd y = \mu_0$, therefore $ \int_0^\infty f(t) \rho(y) \dd y  = \mu_0 f(t)$. Then
\begin{equation*}
	\begin{aligned}
		 f_N(t) - f(t) & =    f_N(t) - \mu_0  f(t) + \left(\mu_0 -1 \right)\, f(t) \\			
	     & =   \int\limits_0^\infty f(ty) \rho(y) \dd y  -  \int\limits_0^\infty f(t) \rho(y) \dd y + \left(\mu_0 -1 \right)\, f(t) \\
		 & =   \left({\mu_0} -1 \right)\, f(t) +  \int\limits_0^\infty (f(ty) - f(t) ) \, \rho(y) \dd y \\ 
		 & =  \left(\mu_0 -1 \right)\, f(t) +  \int\limits_0^\infty (f'(t) (ty-t) + \frac{1}{2} f''(\zeta_y) (ty-t)^2  ) \,  \rho(y) \dd y \\		 
		  &=  \left({\mu_0} -1 \right)\, f(t) + t f'(t)  \int\limits_0^\infty (y-1) \rho(y) \dd y +  \frac{1}{2} t^2 \int\limits_0^\infty f''(\zeta_y) (y-1)^2 \rho(y) \dd y.
		\end{aligned}
\end{equation*}
The above expression is an equality, where between the fourth and fifth lines we used the Taylor series expansion of $f$ centered on $t$, computed at point $ty$, with $\zeta_y$ being the point of the Lagrange remainder. To obtain an upper bound, we take absolute values and get
    \begin{align*}
        & \phantom{\leq}  \,\,\,  | f_N(t) - f(t) | \\
        &\leq \left|\left({\mu_0} -1 \right)\, f(t)\right| + \left| t f'(t)  \int\limits_0^\infty (y-1) \rho(y) \dd y \right| \,+  \,\left| \frac{1}{2} t^2 \int\limits_0^\infty f''(\zeta_y) (y-1)^2 \rho(y) \dd y \right| \span \span \\         
        & \leq \left|{\mu_0} -1 \right| \, |f(t)|+ t \, |f'(t)| \, \left| \int\limits_0^\infty y \rho(y) \dd y - \int\limits_0^\infty \rho(y) \dd y \right| + \frac{1}{2} \, t^2 \, \Vert f'' \Vert_{\infty,\R^+} \int\limits_0^\infty (y-1)^2   \left| \rho(y) \right| \dd y  \\        
        & =   \left|{\mu_0} -1 \right| \, |f(t)| + t \, |f'(t)| \, |\mu_1 - \mu_0| + \frac{1}{2} \, t^2 \, \Vert f'' \Vert_{\infty,\R^+} \, \widetilde{\nu}_2. \span \span \qedhere
    \end{align*}
\end{proof}

We apply this theorem to an Abate--Whitt method, with $\rho=\deltaxx_N$. For the CME method, $\deltaxx_N$ is positive and concentrated near $x=1$, so $\widetilde{\nu}_2(\deltaxx_N) = \nu_2(\deltaxx_N)$ is small and the error bound is small also. For the other methods, $\deltaxx_N$ has both positive and negative components, and unfortunately this makes $\widetilde{\nu}_2(\deltaxx_N)$ orders of magnitude larger than $\nu_2(\deltaxx_N)$ in practical cases. Thus, this error bound is too large to be useful.  However, we can truncate the Taylor series expansion at the first-order term, and compute just an approximation instead of an upper bound. We obtain the following estimate.
\begin{defin}\label{def:est1moment}
    Let $\deltaxx_N$ be the Dirac approximant of an Abate--Whitt method. Let $\mu_0, \mu_1$ be the moments of $\deltaxx_N$ as in Definition~\ref{def:moments}. Let $f:\R^+ \to \R$ be a $\mathcal{C}^1$ function. The first-order moment estimate of the Abate--Whitt approximation error is 
    \begin{equation}
		|f_N(t) - f(t)|  \approx   \left| {\mu_0} -1 \right| \, |f(t)| + t \, |f'(t)| \, |\mu_1 - \mu_0|.
	\end{equation}
\end{defin}
We compare this estimate to the actual error in Figure~\ref{fig:expA_bounds}. 

Note that the coefficients $\abs{\mu_0-1}$ and $\abs{\mu_1-\mu_0}$ can be related to the rational approximation of $e^z$ with $\LT{\deltaxx_N}(-z)$. Recalling~\eqref{moments and Laplace transforms}, we have
$\mu_0(\deltaxx_N) = \int_{0}^{\infty} \deltaxx_N(y) \dd y = \LT{\deltaxx_N}(0)$ and $e^0 = 1$, so 
\[1 - \mu_0 = (e^z - \LT{\deltaxx_N}(-z))\rvert_{z=0}. \]
Hence if $0 \in \Omega$ and the Abate--Whitt method is $\varepsilon$-accurate on $\Omega$, then $|\mu_0-1| \leq \varepsilon$. Similarly, 
\[\mu_1-\mu_0=\mu_1-1+1-\mu_0 =\left( \frac{\dd}{\dd z}(\LT{\deltaxx_N}(-z)-e^z) \right)\biggr\rvert_{z=0} + (e^z - \LT{\deltaxx_N}(-z))\biggr\rvert_{z=0}.\]
Hence $\abs{\mu_1-\mu_0}$ is small when both the value and the first derivative in zero of the rational approximation $\widehat{\deltaxx_N}(-z)$ are close to those of $e^z$.

\section{Laplace transforms in queuing theory} \label{sec:queues}
In this section, we specialize our general bounds to the functions appearing in some applications in queuing theory.
\subsection{Continuous-time Markov chains}
The \emph{generator matrix} (or \emph{rate matrix}) of a continuous-time Markov chain (CTMC) is a matrix $Q\in\mathbb{R}^{d\times d}$ such that $Q_{ij} \geq 0$ whenever $i\neq j$ and $Q\ones = 0$, where $\ones$ is the vector with all entries equal to 1. In particular, $Q$ is a singular $-M$-matrix.

The time-dependent distribution of a CTMC with initial probability distribution $\boldsymbol{\pi}_0 \in \mathbb{R}^d$ is given by
\begin{equation} \label{pictmc}
    \boldsymbol{\pi}(t) = \boldsymbol{\pi}_0 \exp(Qt).
\end{equation}
On the other hand, its Laplace transform
\[
\LT{\boldsymbol{\pi}}(s) = \boldsymbol{\pi}_0 (sI-Q)^{-1}
\]
has a simpler algebraic expression, which is more convenient to work with. One of the reasons for this convenience is an appealing probabilistic interpretation for $s>0$: $\LT{\boldsymbol{\pi}}(s)$ is the expected state of the Markov chain at time $\tau$, where $\tau \sim \operatorname{Exp}(s)$ is an exponentially distributed random variable \cite{encystat}. Often one can resort to algorithms and arguments that exploit the connection between this time $\tau$ and the many other exponentially-distributed random variables that appear in the theory of CTMCs; see for instance~\cite{Tel2022}.

This connection is useful also when working with other quantities defined in terms of a CTMC. An important example are \emph{phase-type distributions}, which model the time it takes for a CTMC to reach a specified set of states. A phase-type distribution is a probability distribution on $[0,\infty)$ with probability density function (pdf) 
\begin{equation}\label{phpdf}
    f(t) = \boldsymbol{\alpha}^T \exp(tQ)\mathbf{q}, \quad \mathbf{q} =  -Q\ones \geq 0,
\end{equation}
where $\alpha\in\mathbb{R}^{d}$ is a stochastic vector, $\ones\in\mathbb{R}^{d}$ is the vector of all ones, and $Q$ is a \emph{sub}generator matrix, i.e., a matrix such that $Q_{ij}\geq 0$ for all $i\neq j$ and $-Q\ones \geq \mathbf{0}$.

For a probability distribution, one usually deals with the Laplace transform of its pdf, which is known as the \emph{Laplace-Stieltjes transform}. Phase-type distributions appear together with their Laplace-Stieltjes transforms in various contexts in queuing theory; see for instance the examples in~\cite{AbaChoWhi2000}.

Clearly both~\eqref{pictmc} and~\eqref{phpdf} are ME functions, hence the bounds in Section~\ref{sec:AWaccuracy} can be applied. With some manipulations, we obtain the following result.

\begin{thm} \label{thm:phasetype}
    Let $f(t)$ be the pdf of a phase-type distribution~\eqref{phpdf}, and $F(t) = \int_0^t f(t)\,dt$ be its corresponding cumulative distribution function (CDF). Let $(w_n, \beta_n)_{n=1}^N$ be an Abate--Whitt method that is $\varepsilon$-accurate on a region $\Omega$ which contains $\mathbb{W}(tQ)$.

    Then, the following bound holds for the inverse Laplace transform of $f(t)$:
    \[
    \abs{f(t) - f_N(t) } \leq (1+\sqrt{2})\varepsilon \norm{\mathbf{q}}_1,
    \]
    and the following two bounds for that of $F(t)$, assuming $0 \in \Omega$:
    \begin{align*}
    \abs{F(t) - F_N(t) } &\leq \varepsilon + (1+\sqrt{2})\varepsilon   \sqrt{d},\\
    \abs{F(t) - F_N(t) } &\leq \varepsilon + (1+\sqrt{2})\varepsilon  (-\boldsymbol{\alpha}^T Q^{-1} \ones) \norm{\mathbf{q}}_1.
    \end{align*}
\end{thm}
Note that $-\boldsymbol{\alpha}^T Q^{-1}\ones$ is the first moment of the phase-type distribution~\cite[Page~6105]{encystat}.
\begin{proof}
    For the first bound, it is enough to use Corollary~\ref{cor:classME} and note that $\norm{\boldsymbol{\alpha}}_2 \leq \norm{\boldsymbol{\alpha}}_1 = \boldsymbol{\alpha}^T \ones=1$, and $\norm{\mathbf{q}}_2 \leq \norm{\mathbf{q}}_1$.

    For the first bound, we integrate $f(t)$ to get the well-known expression for the CDF of a phase-type distribution
    \[
    F(t) = 1 - \boldsymbol{\alpha}^T \exp(tQ) \ones.
    \]
    The function $F(t)$ is the sum of the constant function $F^{(1)}(t) =1$ and $F^{(2)}(t) = - \boldsymbol{\alpha}^T \exp(tQ) \ones$; hence we have
    \[
    \abs{F(t) - F_N(t) } \leq \abs{F^{(1)}(t) - F^{(1)}_N(t) } + \abs{F^{(2)}(t) - F^{(2)}_N(t) } \leq \varepsilon + \abs{F^{(2)}(t) - F^{(2)}_N(t) },
    \]
    where we can bound the first term with $\varepsilon$ since the constant function $1$ is SE with $\alpha=0$. The second term is a ME function; to bound it, we can use Corollary~\ref{cor:classME}, leading to
    \[
    \abs{F^{(2)}(t) - F^{(2)}_N(t) } \leq \varepsilon (1+\sqrt{2}) \norm{\boldsymbol{\alpha}}_2 \norm{\ones}_2 
    \leq \varepsilon (1+\sqrt{2}) \norm{\boldsymbol{\alpha}}_1 \norm{\mathbf{1}}_2 
    = \varepsilon (1+\sqrt{2}) \sqrt{d}.
    \]
    Alternatively, since $Q, Q^{-1}$ and $\exp(tQ)$ commute, we can write
    \[
    F^{(2)}(t) = -\boldsymbol{\alpha}^T \exp(tQ)\ones = -\boldsymbol{\alpha}^T Q^{-1}\exp(tQ)Q\ones.
    \]
    This expression leads to the second bound, since
    \[
    \norm{(-\boldsymbol{\alpha}^T Q^{-1})^T}_2 \leq \norm{(-\boldsymbol{\alpha}^T Q^{-1})^T}_1 = -\boldsymbol{\alpha}^T Q^{-1}\ones.
    \]
    Indeed, $Q^{-1}\leq 0$, so the vector $-\boldsymbol{\alpha}^T Q^{-1}$ has non-negative entries and its 1-norm reduces to the sum of its entries.
\end{proof}

Unfortunately, there is no simple expression for the field of values $\mathbb{W}(Q)$ of a generator or subgenerator matrix, but the following results give inclusions. Let us recall that a \emph{uniformization rate} for a (sub)generator matrix $Q\in\mathbb{R}^{d\times d}$ is any $\lambda\in\mathbb{R}$ such that $\lambda \geq \max_{i=1,\dots,d} \abs{Q_{ii}}$; the definition comes from the concept of \emph{uniformization}, a popular tool to discretize CTMCs~\cite[Section~6.7]{Rossbook}.

\begin{thm}\label{thm:Q_circle}
    Let $Q\in\mathbb{R}^{d\times d}$ be a generator or sub-generator matrix, and $\lambda$ be a uniformization rate for it. Then, $\mathbb{W}(Q) \subseteq B(-\lambda, \lambda\sqrt{d})$.
\end{thm}
\begin{proof}
A sub-generator matrix can always be written as the principal submatrix of a $(d+1)\times (d+1)$ generator matrix, so in view of Lemma~\ref{lem:Wproperties} we reduce to the case in which $Q$ is a generator. If $Q$ is a generator matrix and $\lambda$ is a uniformization rate, then $P=I + \frac{1}{\lambda} Q$ is a stochastic matrix. In particular, $\norm{P}_\infty=1$, where $\norm{\cdot}_{\infty}$ is the operator norm induced by the max-norm on vectors. It is a classical inequality~\cite[Page 50-5]{handbook} that $\norm{P}_{2} \leq \sqrt{d} \norm{P}_{\infty}$; in particular we have
\[
\mathbb{W}(P) \subseteq B(0,\norm{P}_2) \subseteq B(0,\sqrt{d}).
\]
This inclusion implies the result, thanks again to Lemma~\ref{lem:Wproperties}.
\end{proof}
A stronger bound is the following:
\begin{thm}\label{thm:Q_rectangle_large}
    The smallest rectangle in the complex plane  (with its sides parallel to the axes) such that $\mathbb{W}(Q) \subset \Omega_d$ for each generator or subgenerator matrix $Q\in\mathbb{R}^{d\times d}$ with uniformization rate $\lambda$ is
    \[\Omega_d =  \begin{cases}
        [-2\lambda, \frac{\lambda}{2}(\sqrt{2}-1)] + i [-\frac12\lambda,\frac12\lambda], & d=2;\\
        [-(1+\frac{\sqrt{5}}{2})\lambda, \frac{\lambda}{2}\sqrt{3}-1] + i [-\frac{\sqrt{3}}{2}\lambda,\frac{\sqrt{3}}{2}\lambda], & d=3;\\
        [-(1+\frac{\sqrt{6}}{2})\lambda, \frac{\lambda}{2}] + i [-\lambda,\lambda], & d=4;\\
        \begin{multlined}
        \left[-(1+\frac{\sqrt{d+2}}{2})\lambda, \frac{\lambda}{2}(\sqrt{d}-1)\right] + \\
        i\left[-\frac{\sqrt{2}\lambda}{4} \sqrt{d+\sqrt{d^2-4d+12}}, \frac{\sqrt{2}\lambda}{4} \sqrt{d+\sqrt{d^2-4d+12}}\right],
        \end{multlined}
          & d\geq 5.\end{cases}
    \]
\end{thm}
\begin{proof}
    This result follows from the argument in Theorem~\ref{thm:Q_circle}, paired with the bound on the field of values of a stochastic matrix $\mathbb{W}(P)$ given in~\cite[Theorem~6.9]{GauWW2016}.
\end{proof}

These bounds are valid for a generic generator matrix of given size. If we have information on the eigenvalues of the real and imaginary parts of a matrix, sharper bounds can be obtained for its field of values. The next theorem follows from the results in \cite[Section~5.6]{FOVbook}, taking a very coarse discretization $\theta \in \{0,\frac{\pi}{2}, \pi, \frac{3\pi}{2}\}$.
\begin{thm} \label{thm:Q_rectangle_small}
Let $A\in \C^{d \times d}$ be a matrix. Let $X = \frac{A+A^*}{2}$ be the real part of $A$ and $Y=\frac{A-A^*}{2i}$ be the imaginary part of $A$. $X$ and $Y$ are real matrices and thus have real eigenvalues. We have 
\[\mathbb{W}(A) \subseteq \left[\lambda_{\min}\left(X\right),\lambda_{\max}\left(X\right)\right] + i \left[\lambda_{\min}\left(Y\right),\lambda_{\max}\left(Y\right)\right].\]
\end{thm}

\subsection{Fluid queues}
A setting where we can obtain useful results is the one of \emph{Markov-modulated fluid models}, or \emph{fluid queues}~\cite{asmussen95,kk95,roger94,ram99}. A fluid queue with transition matrix $Q\in\mathbb{R}^d$ and \emph{rates} $r_1,\dots,r_d$ is an infinite-dimensional continuous-time Markov process $(\varphi(t), \ell(t))$, where the random variable $\varphi(t) \in \{1,\dots,d\}$ (\emph{phase}) is a CTMC with transition matrix $Q$, and the random variable $\ell(t)\in\mathbb{R}$ (\emph{level}) is a continuous function of $t$ that evolves according to $\frac{d}{dt} \ell(t) = r_{\varphi(t)}$. This model is often paired with boundary conditions, for instance to enforce $\ell(t) \geq 0$. Fluid queues model buffers, telecommunication queues, or performance measures associated to being in a state $\varphi(t)$ for a certain period of time. A fundamental quantity to analyze their transient distribution is the so-called \emph{first-return matrix}. Let $\varphi(0)=i$ be a phase with $r_i>0$, so that $\ell(t)$ is increasing for $t=0$. Set
\[
\tau = \min \{t>0 \colon \ell(t) = \ell(0)\},
\]
the \emph{first-return time to the initial level}, with the convention that $\tau=\infty$ if the level never returns to $\ell(0)$. The (time-dependent) first-return matrix $\Psi(t)$ is then defined by
\[
[\Psi(t)]_{ij} = \mathbb{P}[\tau \leq t, \varphi(\tau)=j \mid \varphi(0)=i],
\]
i.e., the probability that the first return happens before time $t$, and at the same time the phase changes from $i$ to $j$. Usually, one includes in $\Psi(t)$ only phases $i$ and $j$ with $r_i>0$ and $r_j<0$, since these conditions must hold for a return to level $0$ to be possible; hence $\Psi(t) \in \mathbb{R}^{d_+ \times d_-}$, where $d_+, d_- \leq d$ are the number of positive and negative rates, respectively.

The matrix-valued function $\Psi(t)$ is the cumulative density function (CDF) of the return time: $\Psi(t)$ is increasing in $t$, and converges for $t\to\infty$ to a finite substochastic matrix $\Psi(\infty)$. Its derivative $\frac{d}{dt}\Psi(t) = \psi(t)$ is the corresponding probability density function (pdf), which is non-negative for each $t$ and converges to $0$ for $t\to\infty$. One can compute the Laplace-Stieltjes transform $\LT{\psi}(s)$ (and with it $\LT{\Psi}(s) = \frac{1}{s}\LT{\psi}(s)$) as the solution of a certain nonsymmetric algebraic Riccati equation; see \cite{BeaOReTay2005,BeaOReTay2008} for more detail and several algorithms. On the other hand, algorithms to compute $\Psi(t)$ directly \cite{BarSerTel01,BeaNguPol2016} are less common and more complicated. Once again, the convenience of working with Laplace transforms is related to the physical interpretation (for $s>0$) of $\LT{\Psi}(s)$ as the probability that the first-return time is smaller than a random variable with exponential distribution $\operatorname{Exp}(s)$.

In order to apply the techniques introduced earlier, we recall the following expression for $\Psi(t)$. 
\begin{thm}[\protect{\cite[Lemma~2]{BeaNguPol2016}}]
    Consider a fluid queue model, and let $\lambda$ be a uniformization rate for its generator matrix $Q$. Then, there exist nonnegative matrices $\Psi_k^{-} \in \mathbb{R}^{d_+\times d_-}$, for $k=1,2,\dots$ (dependent on $\lambda$), such that
    \begin{equation} \label{Psit}
    \Psi(t) = e^{-\lambda t}\sum_{h=1}^\infty \frac{(\lambda t)^h}{h!} \sum_{k=1}^h \Psi_k^-
    =
    e^{-\lambda t} \sum_{k=1}^{\infty} \Psi_k^- \sum_{h= k}^{\infty} \frac{(\lambda t)^h}{h!}.
    \end{equation}
    Moreover, $\lim_{t\to\infty }\Psi(t) = \Psi(\infty) = \sum_{k=1}^{\infty} \Psi_k^{-}$.
\end{thm}
The matrices $\Psi_k^{-}$ have a physical interpretation using uniformization: they represent the probability that the first return to level $\ell(0)$ happens between the $n$th and $n+1$st uniformization event; see~\cite{BeaNguPol2016} for more detail.

Differentiating~\eqref{Psit} term by term, we can obtain an expression for $\psi(t)$, i.e.,
\begin{equation} \label{psit}
    \psi(t) = \frac{d}{dt} \Psi(t) = \lambda e^{-\lambda t} \sum_{k=1}^\infty  \Psi_k^- \frac{(\lambda t)^{k-1}}{(k-1)!} .
\end{equation}

Following the same approach as Corollary~\ref{cor:classME}, we can get the following result.
\begin{thm} \label{thm:psibound}
    Consider a fluid queue with uniformization rate $\lambda$, and let $f(t) = [\psi(t)]_{ij}$.  Let $(w_n, \beta_n)_{n=1}^N$ be an Abate--Whitt method that is $\varepsilon$-accurate on a region $\Omega$ which contains $B(-t\lambda,t\lambda)$. 
    
    Then, the error of the method at point $t$ is bounded by
    \[
     \left| f(t) - f_N(t) \right|< \varepsilon \lambda (1+\sqrt{2}) [\Psi(\infty)]_{ij}.
    \]
\end{thm}
\begin{proof}

Let
    \begin{equation} \label{JexptJ}
    Q = \begin{bmatrix}
        -\lambda & \lambda\\
        & -\lambda & \lambda \\
        & & \ddots & \ddots\\
        & & & -\lambda & \lambda\\
        & & & & -\lambda
    \end{bmatrix},\quad
    \exp(tQ) = e^{-t\lambda}\begin{bmatrix}
         1 & t\lambda & \frac{(t\lambda)^2}{2} & \dots & \frac{(t\lambda)^{(K-1)}}{(K-1)!}\\
        & 1 & t\lambda & \ddots & \vdots\\
        & & 1 & \ddots & \frac{(t\lambda)^2}{2}\\
        & & & \ddots & t\lambda\\
        & & & & 1
    \end{bmatrix}
    \end{equation}
be a scaled Jordan block of size $K\times K$ and its matrix exponential. Note that $\mathbb{W}(tQ) \subseteq B(-t\lambda,t\lambda)$, thanks to the properties in Lemma~\ref{lem:Wproperties}. We truncate~\eqref{psit} after the first $K$ terms; thanks to the expression above of $\exp(tQ)$, we see that
\[
    f^{(K)}(t) = \lambda \sum_{k=1}^K \exp(-\lambda t) \frac{(\lambda t)^{(k-1)}}{(k-1)!} [\Psi^-_k]_{ij} =  \boldsymbol{\alpha}^T \exp(tQ) \mathbf{q},
\]
with
\begin{equation} \label{vw}
    \boldsymbol{\alpha} = \begin{bmatrix}
        [\Psi^-_K]_{ij} \\
        [\Psi^-_{K-1}]_{ij} \\
        \vdots \\
        [\Psi^-_2]_{ij} \\
        [\Psi^-_1]_{ij}
    \end{bmatrix}, \quad
    \mathbf{q} = -Q\mathbf{1} = \begin{bmatrix}
        0\\0\\
        \vdots\\
        0\\
        \lambda
    \end{bmatrix}.
\end{equation}
This is, essentially, a phase-type distribution, but up to a scaling factor, because in general
\[
\boldsymbol{\alpha}^T \mathbf{1} = \sum_{k=1}^K [\Psi_k^-]_{ij} \neq 1.
\]
Now Corollary~\ref{cor:classME} gives
\begin{align*}
\abs{f^{(K)}(t) - f^{(K)}_N(t)} &\leq (1+\sqrt{2})\varepsilon \,\norm{\boldsymbol{\alpha}}_2 \norm{\mathbf{q}}_2 \\
&\leq  (1+\sqrt{2})\varepsilon \,\norm{\boldsymbol{\alpha}}_1 \norm{\mathbf{q}}_1 \\
&=  (1+\sqrt{2})\varepsilon \lambda \sum_{k=1}^K [\Psi^-_k]_{ij} \leq  (1+\sqrt{2})\varepsilon \lambda \, [\Psi(\infty)]_{ij}.
\end{align*}
Since this bound is uniform in $K$, we can pass to the limit and obtain a bound for $f(t)-f_N(t)$.
\end{proof}
With a little more work, one can obtain an accuracy bound for $\Psi(t)$ as well.
\begin{thm}
    Consider a fluid queue with uniformization rate $\lambda$, and let $F(t) = [\Psi(t)]_{ij}$. Then, the error of the method at point $t$ is bounded by
    \[
     \left| F(t) - F_N(t) \right|< \varepsilon F(\infty) + (1+\sqrt{2}) \varepsilon \lambda\,\mathbb{E}[\psi(t)]_{ij}.
    \]
    Here, $\mathbb{E}[\psi(t)] = \int_0^\infty t\psi(t) \dd t$ is the first moment of the function $\psi(t)$.
\end{thm}
\begin{proof}
    We compute a formula for $\Psi(t)$ analogous to the expression of the CDF of a phase-type distribution. First, we make some algebraic manipulations in the summations to obtain
    \begin{align*}
    \Psi(t) &= e^{-\lambda t}\sum_{h=1}^\infty \frac{(\lambda t)^h}{h!} \sum_{k=1}^h \Psi_k^-\\
    &=e^{-\lambda t} \left(\sum_{h=0}^{\infty}\frac{(\lambda t)^h}{h!}\right)\left(\sum_{k=1}^{\infty} \Psi_k^-\right)
    -
    e^{-\lambda t}\sum_{h=0}^\infty \frac{(\lambda t)^h}{h!} \sum_{k=h+1}^\infty \Psi_k^-\\
    &= 
    \left(\sum_{k=1}^{\infty} \Psi_k^-\right)
    -
    e^{-\lambda t}\sum_{h=0}^\infty \frac{(\lambda t)^h}{h!} \sum_{k=h+1}^\infty \Psi_k^-.
    \end{align*}
    Then, we truncate the summations after the term $\Psi^-_K$, to get
    \[
    \Psi^{(K)}_{ij}(t) = \left(\sum_{k=1}^{K} [\Psi_k^-]_{ij}\right)
    -
    e^{-\lambda t}\sum_{h=0}^\infty \frac{(\lambda t)^h}{h!} \sum_{k=h+1}^K [\Psi_k^-]_{ij}
    = \boldsymbol{\alpha}^T \mathbf{1} - \boldsymbol{\alpha}^T \exp(tQ) \mathbf{1}.
    \]
    with $Q,\boldsymbol{\alpha}$ as in~\eqref{JexptJ} and~\eqref{vw} above. 
    
    Arguing as in the proof of Theorem~\ref{thm:phasetype},  with the only difference that $\boldsymbol{\alpha}^T \mathbf{1} \neq 1$, we have
    \begin{align*}
        \abs{\Psi^{(K)}(t) - \Psi^{(K)}_N(t)} \leq \varepsilon \alpha\mathbf{1} + (1+\sqrt{2})\varepsilon \lambda (-\boldsymbol{\alpha}^T Q^{-1}\mathbf{1}) \leq \varepsilon F(\infty) + 
        (1+\sqrt{2})\varepsilon \lambda (-\boldsymbol{\alpha}^T Q^{-1}\mathbf{1}).
    \end{align*}
    It remains to show that $(-\boldsymbol{\alpha}^T Q^{-1}\mathbf{1})$ converges to the first moment of $\Psi(t)_{ij}$, when $K\to\infty$. To this purpose, we compute the remainder
    \begin{align*}
    \mathbb{E}[f(t)] - \mathbb{E}[f^{(K)}(t)] 
    &= \int_0^{\infty} \left(tf(t) - tf^{(K)}(t)\right)\, \dd t\\ 
    &= \int_0^{\infty}  (\lambda t) \sum_{k=K+1}^\infty \exp(-\lambda t) \frac{(\lambda t)^{(k-1)}}{(k-1)!} [\Psi_k^-]_{ij}\, \dd t\\
    &= \sum_{k=K+1}^\infty [\Psi_k^-]_{ij} \int_0^{\infty}  (\lambda t)\exp(-\lambda t) \frac{(\lambda t)^{(k-1)}}{(k-1)!}\,\dd t\\
    &=\sum_{k=K+1}^\infty [\Psi_k^-]_{ij} \frac{k}{\lambda}.
    \end{align*}
    In particular, when $K=0$, this computation gives an expression for the first moment of $\psi$ as a non-negative series
    \[
    \mathbb{E}[f(t)] = \int_0^\infty tf(t)\, \dd t = \sum_{k=1}^\infty [\Psi_k^-]_{ij} \frac{k}{\lambda}
    \]
    If $\mathbb{E}[\psi(t)]_{ij} < \infty$, the series must be summable, and this implies that
    \[
    \lim_{K\to\infty} \sum_{k=K+1}^\infty [\Psi_k^-]_{ij} \frac{k}{\lambda} = 0. \qedhere
    \]
\end{proof}

\section{AAA approximation} \label{sec:AAA}
We have seen that the error of the Abate--Whitt approximant depends on the $\varepsilon$-accurate approximation of $e^{z}$ on $\Omega$ with the rational function $\LT{\deltaxx_N}(-z)$. Hence, we can construct accurate Abate--Whitt methods by choosing suitable approximations. Zakian used the Padé approximant of $e^{z}$ at the point $z=0$, which can be useful if $\Omega$ is close to $z=0$, but for generic $\Omega$ it provides a worse approximation. 

Therefore, we need an approach that can find a good $\varepsilon$-accurate rational approximation on an arbitrary region $\Omega$. There are many algorithms in the literature; the state of the art is the AAA algorithm, proposed by Nakatsukasa, Sète and Trefethen \cite{AAA2018}. We refer to the paper for more details and below we briefly present the algorithm. 

\subsection{AAA algorithm}
The AAA algorithm computes a rational function which approximates a prescribed $f(z)$ on a given finite set of points $Z \subseteq \mathbb{C}$.
One of the key tools is the so-called \emph{barycentric representation}
\begin{equation}\label{eq: barycentric representation r(z)}
r(z) = \frac{n(z)}{d(z)} = \frac{\displaystyle\sum\limits_{k=1}^K \dfrac{u_k f_k}{z-z_k}}{\displaystyle\sum\limits_{k=1}^K \dfrac{u_k}{z-z_k} }.
\end{equation}
Here $u_1,\dots,u_K, z_1,\dots,z_K,f_1,\dots,f_K \in \mathbb{C}$ are parameters that will be chosen appropriately. Both $n(z)$ and $d(z)$ are rational functions of degree $(K-1,K)$. Since $z-z_k$ appears in denominators, one may think that the $z_k$ must be poles of $r(z)$. However, clearing denominators we obtain the equivalent expression
\begin{equation} \label{AAA alternative form}
r(z) = \frac{\displaystyle\sum\limits_{k=1}^K u_k f_k \prod\limits_{i \neq k}(z-z_i) }{\displaystyle\sum\limits_{k=1}^K u_k \prod\limits_{i \neq k}(z-z_i) }.
\end{equation}
By evaluating~\eqref{AAA alternative form} for $z=z_k$, in both $n(z)$ and $d(z)$ all summands except one vanish, and we obtain $r(z_k) = f_k$. Therefore we see that $r(z)$ is a rational function of degree at most $(K-1,K-1)$ that takes the values $f_1, \ldots, f_K$ in the \emph{support points} $z_1, \ldots, z_K$. 

Evaluating a rational function in the barycentric representation~\eqref{eq: barycentric representation r(z)} has surprisingly good numerical stability properties: even though the denominators $z - z_k$ can become arbitrarily large, the ratio $\frac{n(z)}{d(z)}$ stays bounded, and the floating-point errors in computing $\frac{1}{z-z_k}$ for $z\approx z_k$ cancel out, since that sub-expression appears in both $n(z)$ and $d(z)$; see~\cite{AAA2018} for more discussion.

The AAA algorithm takes as its input a function $f(z)$ that we wish to approximate with a rational function, and a finite set of points  $Z \subseteq \C$, with $|Z|=L$, on which $f(z) \approx r(z)$ should hold.

The algorithm chooses iteratively inside $Z$ a set of points to use as support points, and adds one new point greedily at each iteration. We describe the iteration step $K+1$, assuming that points $\{z_1, \ldots, z_K\} \subseteq Z$ have already been selected, and set $f_k = f(z_k)$ for $k=1,\dots,K$. In this way, $f(z_k) = f_k = r(z_k) = \frac{n(z_k)}{d(z_k)}$, so the approximation is exact in the support points.

These interpolation conditions are not sufficient to determine uniquely the function $r(z)$ in~\eqref{eq: barycentric representation r(z)}: we also need to choose the weights $u_1,\dots,u_K$. We choose these weights so that $f(z) \approx r(z) = \frac{n(z)}{d(z)}$ on the points of $Z \setminus \{z_1,\dots,z_K\}$, which we shall label $Z_1,Z_2,\dots,Z_{L-K}$. To this purpose, we take
\begin{equation} \label{eq:AAA minimum problem}
(u_1,\dots,u_K) = \arg\min  \left\{\sum_{i=1}^{L-K} \abs*{f(Z_i)d(Z_i) - n(Z_i)}^2 \colon u \in \mathbb{C}^K, \norm{u}=1 \right\}.    
\end{equation}
(Note that we can restrict to $\norm{u}=1$, since $r(z)$ does not depend on the scaling of $u$.) The problem~\eqref{eq:AAA minimum problem} is, in effect, a linear least-squares problem in the weights $u$:
\[
(u_1,\dots,u_K) = \arg\min \{\norm{Au}_2^2 \colon u \in \mathbb{C}^K, \norm{u}=1 \}.
\]
With some computations, one sees that the associated matrix is $A = D_1 C - C D_2$, with
\[
C = \begin{pmatrix}
\dfrac{1}{Z_1-z_1} & \cdots & \dfrac{1}{Z_1-z_K} \\
\vdots & \ddots & \vdots \\
\dfrac{1}{Z_{L-K}-z_1} & \cdots & \dfrac{1}{Z_{L-K}-z_K} 
\end{pmatrix} \in \mathbb{C}^{(L-K)\times K},
\]
and
\begin{align*}
D_1 &= \operatorname{diag}(f(Z_1),\dots,f(Z_k)) \in \mathbb{C}^{(L-K)\times (L-K)},    \\
D_2 &= \operatorname{diag}(f(z_1),\dots,f(z_k)) \in \mathbb{C}^{K\times K}.
\end{align*}
The problem~\eqref{eq:AAA minimum problem} can be solved with some linear algebra tools: in the typical case where $L-K \geq K$, the optimal $u$ is the singular vector associated to the minimum singular value of $A$. By computing $u$, we fix all the parameters in the rational function $r(z)$. 

For the next step of the iteration, we add a new point $z_{K+1}$ to the set of support points: we use a greedy strategy, and select the point $Z_i$ on which the current approximation is worse:
\[z_{K+1} = \arg\max\{\abs*{f(z) - r(z)} : z \in \{Z_1,\dots,Z_{L-K}\} \}.\]
One continues adding support points (and increasing the degree $K$) until the desired accuracy is reached.

There are few theoretical guarantees for the optimality of the rational functions produced by AAA, since~\eqref{eq:AAA minimum problem} does not ensure that we minimize $\abs{f(z)-r(z)}$, but in practice the algorithm is very efficient and stable, producing results that are very close to the theoretical optima~\cite{AAA2018}.

The form~\eqref{eq: barycentric representation r(z)} does not reveal the poles of $r(z)$ immediately, but in~\cite{AAA2018} it is shown how to compute them. The poles of $r(z)$ are the solutions of the generalized eigenvalue problem
\begin{equation} \label{AAA generalized eigenvalue problem}
    \det \begin{pmatrix}
0 & u_1 & u_2 & \cdots & u_K\\
1 & z_1 -\lambda &     &        & \\
1 &     & z_2-\lambda &        & \\
\vdots & &    & \ddots & \\
1 &     &     &        & z_K-\lambda
\end{pmatrix} = 0.
\end{equation}
Indeed, adding a multiple of the $k$-th row with coefficient $u_k / (\lambda-z_k)$ to the first row, we obtain
\[\det \begin{pmatrix}
\sum\limits_{k=1}^{K} \frac{u_k}{\lambda-z_k} & 0 & 0 & \cdots & 0\\
1 & z_1 -\lambda &     &        & \\
1 &     & z_2-\lambda &        & \\
\vdots & &    & \ddots & \\
1 &     &     &        & z_K-\lambda
\end{pmatrix} = 0.
\]
The first diagonal entry is $d(\lambda)$; so the $K-1$ zeros of $d$ are the solutions of~\eqref{AAA generalized eigenvalue problem}.

Once the poles $\beta_k$ have been computed as the solutions of~\eqref{AAA generalized eigenvalue problem}, the corresponding residues $w_k$ can be obtained with the residue formula from complex analysis, $w_k = \frac{n(\beta_k)}{d'(\beta_k)}$.
While the rest of the AAA algorithm has good floating-point stability properties, solving~\eqref{AAA generalized eigenvalue problem} may be more problematic. We advise using higher precision in this final part of the algorithm.

\subsection{Modifications to AAA}
We have seen the basics of the AAA algorithm; now we apply it to the problem of computing $\varepsilon$-accurate Abate--Whitt methods. We want $\LT{\deltaxx_N}(-z)$ to be a rational approximation of $e^{z}$ on the set $\Omega$. However, $\LT{\deltaxx_N}(-z)$ has degree $(N-1,N)$, while the AAA algorithm produces a rational function of degree $(K-1,K-1)$. A modification of the algorithm to produce rational approximations of a general degree $(M,N)$ is described in~\cite{Fil18}; here, we adopt a simpler approach instead.
We can regard the condition $\deg n(z) < \deg d(z)$ as requiring that $r(\infty)=0$. This relation is similar to the interpolation conditions $r(z_k)=f_k$ that are imposed in the support points. We would like to impose this condition already from the first step, hence starting our iteration from a ``0th support point'' $z_0 = \infty$, $f_0 = f(z_0) = 0$. However, we cannot run the algorithm without modification with a support point equal to $\infty$, but we have to modify slightly the barycentric representation~\eqref{eq: barycentric representation r(z)}.
Namely, we set 
\begin{equation} \label{eq: modified barycentric representation with a pole at infinity}
r(z) = \frac{\displaystyle\sum\limits_{k=1}^K \dfrac{u_k f_k}{z-z_k}}{\displaystyle u_0 + \sum\limits_{k=1}^K \dfrac{u_k}{z-z_k}}.
\end{equation}
Clearing denominators, we see that $r(z)$ in~\eqref{eq: modified barycentric representation with a pole at infinity} has degree $(K-1,K)$. With this representation, the analogue of~\eqref{eq:AAA minimum problem} is
\begin{equation} \label{eq:modified AAA minimum problem}
(u_0,\dots,u_K) = \arg\min \{\norm{\widetilde{A}u}_2^2 \colon u \in \mathbb{C}^{K+1}, \norm{u}=1 \}.
\end{equation}
The associated matrix is $\widetilde{A} = D_1\widetilde{C} - \widetilde{C}\widetilde{D}_2 \in \mathbb{C}^{(L-K)\times (K+1)}$, with
\begin{equation}\label{eq: tildeC}
    \widetilde{C} = \begin{pmatrix}
1 & \dfrac{1}{Z_1-z_1} & \cdots & \dfrac{1}{Z_1-z_K} \\
\vdots & \vdots & \ddots & \vdots \\
1 & \dfrac{1}{Z_{L-K}-z_1} & \cdots & \dfrac{1}{Z_{L-K}-z_K} 
\end{pmatrix}, 
\end{equation}
and
\[\widetilde{D}_2 = \text{diag}(0,f(z_1),\ldots,f(z_{K})).\]
Hence, the only modification required is including an initial column of ones in the matrix $C$, and a corresponding zero in $D_2$; the rest of the iteration can proceed without modifications.

\begin{rem}
With the same strategy, for any given $f_0\in\mathbb{C}$ we can construct a rational approximant such that $r(\infty)=f_0$: it is sufficient to add the term $u_0 f_0$ to the numerator $n(z)$ of~\eqref{eq: modified barycentric representation with a pole at infinity}.
\end{rem}

We also make another modification to the algorithm to ensure that the non-real weights $w_n$ and nodes $\beta_n$ come in conjugate pairs: inside the main loop of the algorithm, whenever we add a non-real $z_k$ to the set of support points, we also add in the same iteration its conjugate $z_{k+1} = \overline{z_k}$.

The computed weights $u_k$ in AAA then also come in conjugate pairs, in exact arithmetic; to compensate for machine arithmetic errors, we replace the pair $(u_k, u_{k+1})$ with $\left(\frac{u_k+\overline{u_{k+1}}}{2}, \frac{\overline{u_{k}}+u_{k+1}}{2}\right)$.

We need a small modification also to the eigenvalue problem~\eqref{AAA generalized eigenvalue problem} to compute the poles; it becomes
\begin{equation} \label{eq: modified AAA eigenvalue problem}
    \det \begin{pmatrix}
0 & u_0 & u_1 & u_2 & \cdots & u_K\\
1 & -1 & \\
1 & & z_1 -\lambda &     &        & \\
1 & &     & z_2-\lambda &        & \\
\vdots&  & &    & \ddots & \\
1 & &     &     &        & z_K-\lambda
\end{pmatrix} = 0.
\end{equation}
Indeed, we can carry out row operations that preserve the determinant to transform this problem into the equivalent one
\begin{equation}
    \det \begin{pmatrix}
u_0 + \sum_{k=1}^K \frac{u_k}{\lambda-z_k} & 0 & 0 & 0 & \cdots & 0\\
1 & -1 & \\
1 & & z_1 -\lambda &     &        & \\
1 & &     & z_2-\lambda &        & \\
\vdots&  & &    & \ddots & \\
1 & &     &     &        & z_K-\lambda
\end{pmatrix} = 0.
\end{equation}
Solving this eigenvalue problem provides the parameters of an Abate--Whitt method: since
\[
r(z) = \LT{\deltaxx_K}(-z) = \sum\limits_{k=1}^K \frac{w_k}{\beta_k-z},
\]
we recover $\beta_k$ as the poles of $r(z)$, and $w_k$ as the corresponding residues.

The AAA algorithm with these modifications is described in Algorithm~\ref{algo:AAA}. We note that the number of computed support points $N$ may be either $N_{\max}$ or $N_{\max}-1$; the second case happens when we would like to add a pair of conjugate support points, but there is no room to do it. The number $N'$ is also not known a priori: the poles are computed only at the end of the loop, and we do not know in advance how many are real.
\begin{algorithm}
    \caption{The modified AAA algorithm} \label{algo:AAA}
    \KwData{Finite set of points $Z\subseteq \mathbb{C}$, function $f$, maximum order $N_{\max}$, tolerance $\varepsilon$}
    \KwResult{Poles and weights $(w_n, \beta_n)_{n=1}^N$, either real in or conjugate pairs, such that $\|f(z) - \sum_{n=1}^N \frac{w_n}{\beta_n-z}\|_{\infty,Z}\leq \varepsilon$ and $r(\infty)=0$ with $N \leq N_{\max}$.}
    $\widetilde{C} \gets \mathbf{1}_{|Z|}$\;
    $u_0 \gets 0$\;
    $K \gets 0$\;
    \While{$\|f-r\|_{\infty,Z} > \varepsilon$ and $K < N_{\max}$}{
        $z_{K+1} = \arg\max\{\abs*{f(z) - r(z)} : z \in Z \setminus \{z_1,\dots,z_K\}\}$\;
        Update $\widetilde{C}$ according to~\eqref{eq: tildeC}: remove row corresponding to $z_{K+1}$, add column corresponding to it\;
        \uIf{$z_{K+1}$ is real}{
            $K \gets K+1$\;
        }
        \Else{
            \If{$K+2 > N_{\max}$}{
                \tcp{No room to add a conjugate pair:}
                \tcp{ignore $z_{K+1}$ and terminate with $K$ support points}
                \textbf{break}\;
            }
            $z_{K+2} \gets \overline{z_{K+1}}$\;
            Update $\widetilde{C}$ according to~\eqref{eq: tildeC}: remove row corresponding to $z_{K+2}$, add column corresponding to it\;
            $K \gets K+2$\;
        }
        $(u_0,u_1,\dots,u_K) \gets \arg\min \{\norm{\widetilde{A}u}_2^2 \colon u \in \mathbb{C}^{K+1}, \norm{u}=1 \}$\;
    }
    Compute $\beta_1,\dots,\beta_K$ by solving~\eqref{eq: modified AAA eigenvalue problem} (in higher precision)\;
    Compute residues $w_k = \frac{n(\beta_k)}{d'(\beta_k)}$, $k=1,\dots,K$\;
\end{algorithm}
\section{Computing TAME parameters} 
\label{sec:TAME}
\subsection{Floating-point precision issues} \label{sec:floatingpoint}
While some previous works focus on arbitrary precision arithmetic~\cite{AbaVal2004,AbaWhi2006}, we deal with the case in which the computations in an Abate--Whitt method~\eqref{eq: Abate Whitt definition} are performed in the standard IEEE754 \texttt{binary64} (double precision), possibly using poles and weights $(w_n,\beta_n)$ precomputed in higher precision.

We have seen that, in exact arithmetic, the error of the ILT approximation is bounded as
$\abs{f(t)-f_N(t)} \leq C \varepsilon$,
where $\varepsilon$ is the rational approximation error and $C$ is a constant depending on the class of $f$. We now assess the impact of inaccuracies in the computation of $\LT{f}$, such as the ones due to floating-point arithmetic.

Let us assume for simplicity that $t=1$. Let $g(\beta_n)$ be the  value obtained by computing $\LT{f}(\beta_i)$ numerically, and suppose that it has relative precision $\delta$, i.e.,
\[\abs*{\frac{g(\beta_n) - \LT{f}(\beta_n)}{\LT{f}(\beta_n)} }\leq \delta.\]
In floating-point arithmetic, we can only ensure a $\delta$ at least as large as the machine precision, if not larger.
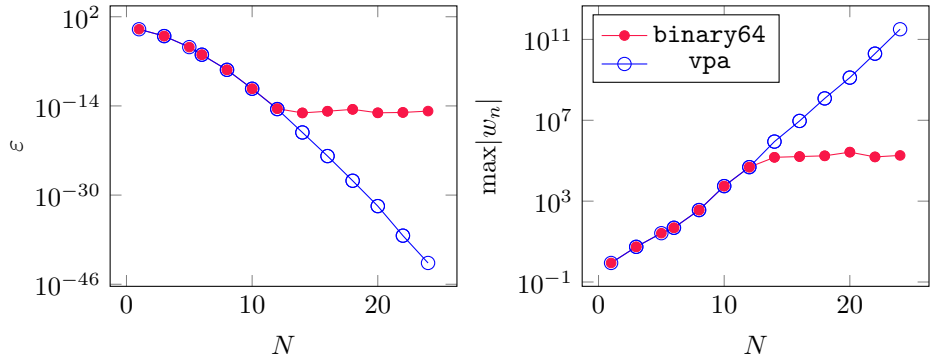
\begin{figure}
    \input{plot_magnitudeW}
    \caption{Approximation error $\varepsilon$ and magnitude of weights $\max\abs{w_n}$ when running AAA with $\Omega=B(-5,5)$ in either \texttt{binary64} or VPA with 100 significant digits.}
    \label{fig:magnitudeW}
\end{figure}

We can write the computed value of the Abate--Whitt approximant as 
\[ \widetilde{f}_N(1) = \sum\limits_{n=1}^N w_n g\left(\beta_n\right) = \sum\limits_{n=1}^N w_n \LT{f}\left(\beta_n\right) (1+\delta_n), \; \; \abs{\delta_n} \leq \delta .\]
Then we have the bound
\[ \abs{\widetilde{f}_N(1) - f_N(1)} \leq \sum\limits_{n=1}^N \abs{w_n \LT{f}\left(\beta_n\right)} \delta  \]
and hence
\begin{equation} \label{Numerical error of AW approximants in double precision}
   \abs{\widetilde{f}_N(1) - f(1)} \leq \sum\limits_{n=1}^N \abs{w_n \LT{f}\left(\beta_n\right)} \delta + C \varepsilon 
\end{equation}
This bound reveals that limiting the magnitude of the weights $w_n$ is important when we evaluate $\LT{f}$ in limited precision, as already noted in~\cite[Section~5.1]{HHAT2020}.
In practice, we observe that increasing $N$ leads to a smaller $\varepsilon$ but also to larger weights; hence increasing $n$ after a certain point (which depends on the machine precision) is no longer beneficial.

\begin{rem}\label{rem:AAA vpa vs double}
    In the TAME method, the weights are computed by optimizing a function that is itself computed numerically, by evaluating the rational approximant $\LT{\deltaxx_N}(-z)$. If the computations in the main AAA loop (Algorithm~\ref{algo:AAA}) are performed with sufficiently many significant digits (e.g., 100 decimal digits of precision), then we observe the approximation error $\varepsilon$ decreasing even below $10^{-16}$, but at the same time the magnitude of the weights increasing steadily. This is detrimental if we plan to compute the ILT in \texttt{binary64}, since the large weights cause a large numerical error irrespective of $\varepsilon$: in~\eqref{Numerical error of AW approximants in double precision} the summation becomes the dominant term, even if $\varepsilon$ is small. If we run the AAA algorithm in \texttt{binary64} instead, then the error stagnates around the machine precision $2.2 \times 10^{-16}$, and at that points increasing $N$ further only adds spurious poles with weights that are of the order of the machine precision; in particular, the magnitude of the weights does not increase anymore. So running the AAA algorithm in \texttt{binary64} acts as a safeguard against increasing weights. We observe this behavior in an example in Figure~\ref{fig:magnitudeW}, and also its consequences later in Section~\ref{sec:expB}.
\end{rem}

\subsection{Choice of $\Omega$}
We have seen in Sections~\ref{sec:AWaccuracy} and~\ref{sec:queues} that an Abate--Whitt method which is $\varepsilon$-accurate on $\Omega$ can recover the original function $f$ with an error proportional to $\varepsilon$ (Theorems~\ref{thm: approximation error SE}, \ref{thm:class2approximation}, \ref{thm:phasetype}, \ref{thm:psibound}, \ref{thm:error of LS class}). We can use the AAA algorithm to construct a TAME method with a small $\varepsilon$. The first step is to select the region $\Omega$, depending on the information we have on $f$.

\begin{itemize}
    \item \textbf{Fluid queues}. If $f$ arises from a fluid queue model  with uniformization rate $\lambda$ (i.e., $f(t) = \Psi(t)$ or $f(t) = \psi(t)$), then use $\Omega = B(-r, r)$, with $r= \lambda t$. 
    
    \item \textbf{ME class}. 
    If $f$ is in the ME class (e.g. $f$ is a phase type distribution), then $\Omega$ should contain $\mathbb{W}(Q)$. With no information on $Q$ apart from its dimension,  Theorem~\ref{thm:Q_rectangle_large} can be used. If the user has more information about $\mathbb{W}(Q)$, tighter bounds can be used, for instance the one in Theorem~\ref{thm:Q_rectangle_small}.

    \item \textbf{LS class}. If $f(t)$ is in the LS class, we use $\Omega = [-L,0]$ with $L$ chosen according to Theorem~\ref{thm:error of LS class}. 

    \item If $f(t)$ can be approximated by a function of the above classes, use the corresponding $\Omega$. Otherwise, if nothing is known about $f$, then we recommend trying three possible domains: the circle $B(-r,r)$, the segment on the real half-line $[-r,0]$, the segment on the imaginary line $i[-r,r]$. However, as we note in Remark~\ref{rem:negative result}, no Abate--Whitt method can give good results for all functions $f$.
\end{itemize}

\begin{rem}
    One may be tempted to use a large region $\Omega$ to cover as many functions as possible; however, usually the magnitude of the weights $\abs{w_n}$ is larger for a bigger region $\Omega$. This observation discourages using overly large sets $\Omega$: if we need to compute an ILT for which we know that the domain $\Omega = B(-1,1)$ is sufficient, then using $\Omega = B(-10,10)$ instead would cause loss of accuracy, at least when the computations are done in double-precision arithmetic.
\end{rem}


\subsection{Choice of $Z$}
\label{sec:choiceZ} 
Once $\Omega$ is chosen, we wish to construct a method that is $\varepsilon$-accurate on $\Omega$ using Algorithm~\ref{algo:AAA}. However, the domain for AAA is a discrete set of points $Z$, while in many of our earlier examples $\Omega$ is a bounded region such as a disc or a rectangle. Hence we need a way to approximate $\Omega$ with a finite set of points.

The simplest approach would be to use a regular grid of points inside $\Omega$, thus requiring a number of points that scales quadratically with the diameter of $\Omega$. We use another more efficient approach instead, following~\cite[Section~6]{AAA2018}. Assume that $\Omega$ is a simply connected domain with boundary $\partial \Omega$. Consider the function $g(z) = e^z-\LT{\deltaxx_N}(-z)$, which is holomorphic on $\C \setminus\{\beta_1, \ldots, \beta_n\}$. If $g$ has no poles inside $\Omega$, then by the maximum modulus principle~\cite[Theorem~10.24]{rudin} $\max_{z \in \Omega}|g(z)| = \max_{z \in \partial\Omega}|g(z)|$. That is, it is enough to require that $|g(z)|$ be small on the boundary of $\Omega$. Therefore we can take $Z$ as a discretization of $\partial \Omega$. In particular, when $\Omega$ is a disc, we take $Z$ to be a set of equispaced points on a circle. 

We then apply the modified AAA algorithm on the support set $Z$, obtaining a rational approximation $\LT{\deltaxx}(-z) = \sum_{n=1}^N \frac{w_n}{\beta_n - z}$ to $e^z$. This rational approximation is guaranteed to be $\varepsilon$-accurate only on $Z$, not on the whole $\Omega$, but in practice this is sufficient if the points in $Z$ are sufficiently tight. 

Finally, we recover the parameters $(w_n, \beta_n)_{n=1}^N$ and obtain the TAME method. We have to check that the poles $\beta_n$ are outside the domain $\Omega$, otherwise the maximum modulus principle may be violated; furthermore, we could not have a $\varepsilon$-accurate method on $\Omega$, since $\LT{\deltaxx_n}(-z)$ would have a pole inside $\Omega$. This condition is satisfied in the TAME methods we computed. 

\subsection{TAME with optimal $N'$ for a given $B(-r,r)$}

As we discussed above, increasing $N$ (or $N'$) beyond a certain value does not improve accuracy, if the transforms are evaluated in~\texttt{binary64}. In this section, we wish to give a strategy to select a quasi-optimal $N'$ to use, for the common case of a region of the form $\Omega=B(-r,r)$ and several values of $r$. First, we display in Figure~\ref{fig: errors vs r for many N} the errors $\varepsilon$ obtained for several values of $r$, and the corresponding values of $N$ and $N'$. 
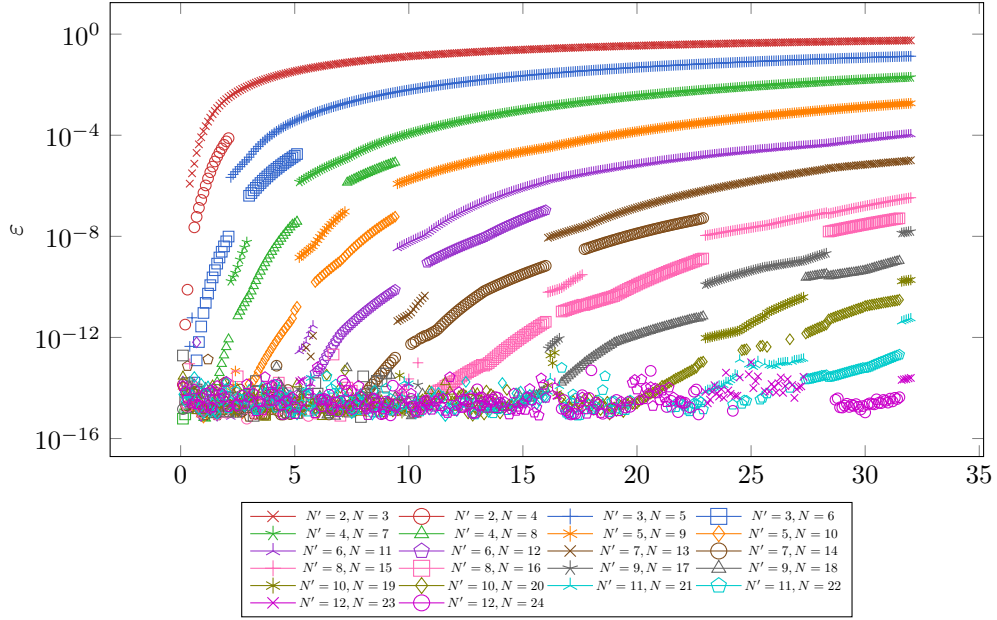
\begin{figure}
    \centering
    \def\currentylabel{$\varepsilon$}
    \def\currentycolumn{error}
    \input{plot_optimal_parameters}
    \caption{Errors $\varepsilon$ obtained by TAME on $\Omega=B(-r,r)$, with the corresponding values of $n$ and $n'$.}
    \label{fig: errors vs r for many N}
\end{figure}
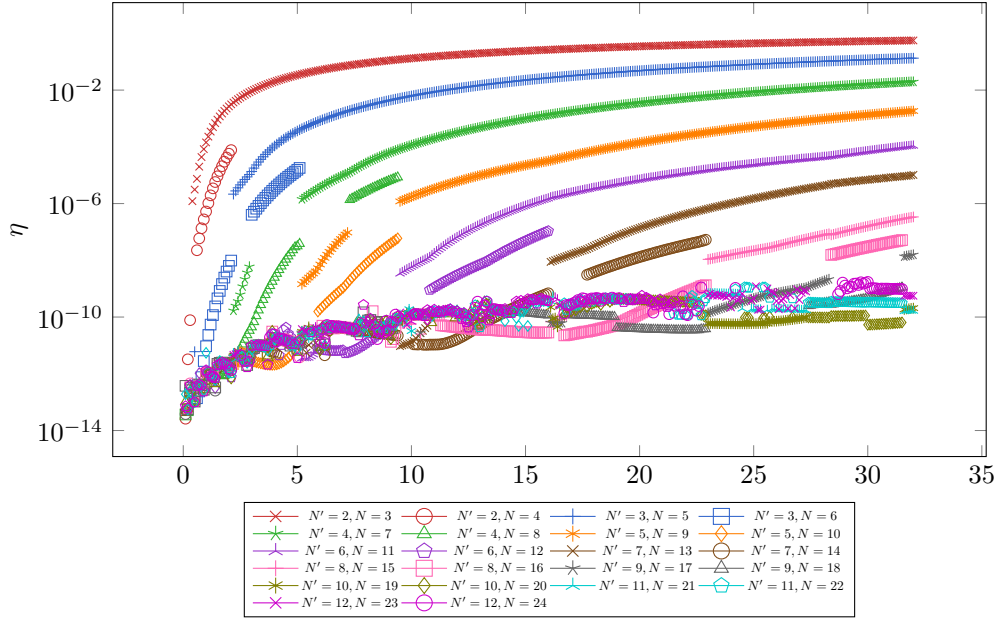
\begin{figure}
    \centering
    \def\currentylabel{$\eta$}
    \def\currentycolumn{estimator}
    \input{plot_optimal_parameters}
    \caption{Error estimator $\eta$ obtained by TAME on $\Omega=B(-r,r)$, with the corresponding values of $n$ and $n'$.}
    \label{fig: error estimator vs r for many N}
\end{figure}
Several comments on this figure are in order.

Since $Z$ is the discretization of a circle, all points in it apart from two are non-real. Hence, our AAA variant almost always adds support points in conjugate pairs. As a consequence, for any given $r$ one obtains only about half of the possible values of $N$. For instance, for $r=30$, TAME produces rational approximants with $N\in\{3,5,7,9,11,13,15,16,18,20,22,24\}$, since a real support point is added as the 16th point. The iterations at which real support points are added vary depending on $r$: this explains why the colored lines are sometimes broken in the figure.

We note that each value of $N'$ is reached with more than one value of $N$: typically, one has $N=2N'-1$ and $N=2N'$, i.e., all poles come in conjugate pairs apart from at most one real pole.

The AAA algorithm was run in \texttt{binary64}, as suggested in Remark~\ref{rem:AAA vpa vs double}, with only the eigenvalue problem~\eqref{eq: modified AAA eigenvalue problem} solved in higher precision. For each value of $N'$ (represented by a different color in the figure), the method reaches a value of $\varepsilon$ close to the machine precision $2.2\cdot 10^{-16}$ for all values of $r$ up to a certain threshold, but then the error increases steadily. For instance, for $N'=6$, the error starts drifting away from machine precision at $r\approx 5$. There is a considerable amount of numerical noise due to floating point errors, clearly visible at the bottom of the graph, so that ``close to machine precision'' may mean a value larger than $10^{-13}$, in certain cases.

It follows from the discussion in Section~\ref{sec:floatingpoint} that $\varepsilon$ is not the only factor that affects the final accuracy of the method when the ILT is evaluated in~\texttt{binary64}. We take as a proxy for the error~\eqref{Numerical error of AW approximants in double precision} the quantity $\eta = \varepsilon + \mathsf{u}\max \abs{w_n}$, where $\mathsf{u} \approx 2.2\cdot 10^{-16}$ is the machine precision. This choice is motivated by the heuristic that (a) $f(\beta_n)$ has the same order of magnitude as $C$ for our classes of functions and for typical values of $\beta_n$, and that (b) the maximum of $\abs{w_n}$ is a better estimate for the error than the worst-case bound $\sum \abs{w_n}$, because the errors $\delta_n$ do not have all the same phase, and the $w_n$ do not have all the same magnitude. We plot this error estimator in Figure~\ref{fig: error estimator vs r for many N}. 

Based on the plot, we selected an array of methods with different values of $N'$. Some relevant quantities for these methods are reported in Table~\ref{tab:optimal TAME}, while their poles and weights are available for download on~\url{https://github.com/numpi/tame-ilt}. When one needs a method for a certain $\Omega = B(-r,r)$, we suggest using the method with the smallest $r_{\max} \geq r$. In this way we have a method that is $\varepsilon$-accurate on $\Omega \subseteq B(-r_{\max},r_{\max})$, and at the same time it does not have unnecessarily large weights.
\begin{table}[]
     \centering
     \caption[0.8\textwidth]{Table of relevant quantities for a family of quasi-optimal precomputed TAME methods. The method in each row was generated with $r=r_{\max}$. 
     }\label{tab:optimal TAME}
\begin{tabular}{ccc}
\phantom{pippopippo}
&
     \begin{tabular}{cccccc}
     \toprule
     $N'$ & $N$ & $r_{\max}$ & $\varepsilon$ & $\max \abs{w_n}$ & Error proxy $\eta$\\
     \midrule
 3 & 6 & 0.6 & 1.637909e-14 & 5.462668e+02 & 1.376747e-13 \\
 4 & 8 & 1.8 & 8.229408e-14 & 3.402665e+03 & 8.378375e-13 \\
 5 & 10 & 4.0 & 3.973128e-13 & 7.669538e+03 & 2.100292e-12 \\
 6 & 12 & 7.0 & 1.420889e-12 & 1.684302e+04 & 5.160790e-12 \\
 7 & 14 & 11.2 & 1.983229e-12 & 3.981692e+04 & 1.082436e-11 \\
 8 & 16 & 16.8 & 1.140301e-11 & 5.302010e+04 & 2.317583e-11 \\
 9 & 18 & 22.7 & 6.075754e-12 & 1.296546e+05 & 3.486487e-11 \\
 10 & 20 & 31.6 & 3.192135e-11 & 1.495150e+05 & 6.512034e-11\\
 \bottomrule
     \end{tabular}
     &
     \phantom{pippopippo}
\end{tabular}
\end{table}

\section{Experiments} 
\label{sec:experiments}
\subsection{Comparison of different Abate--Whitt methods}
To gain additional insight on how the different methods relate to each other, we present a comparison of their parameters. 

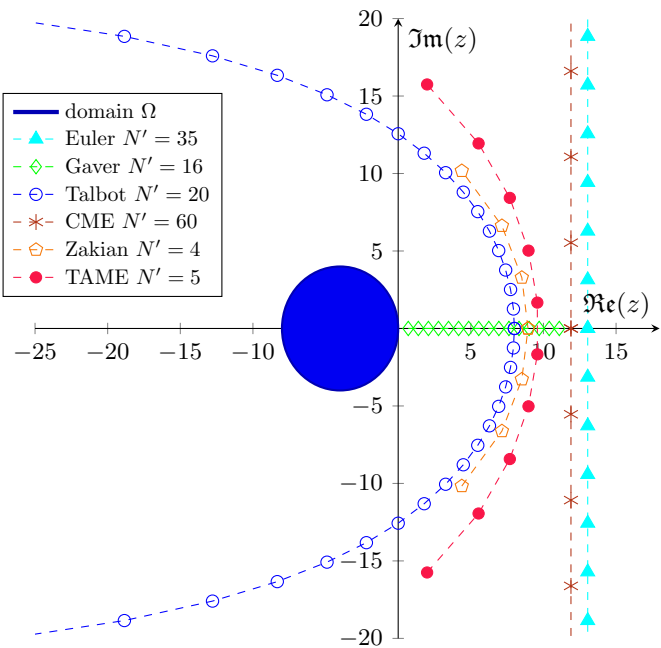
\begin{figure}
\centering
\input{plot_Z-beta.tex}
\caption{Plot of nodes $\beta_n$ for different methods. The TAME method uses $\Omega = B(-4,4)$.}
\label{fig:nodes}
\end{figure}

In Figure~\ref{fig:nodes} the nodes $\beta_n$ (along with their conjugates) are plotted in the complex plane. 
The TAME method uses $N'=5$ and $r=4$ chosen according to Table~\ref{tab:optimal TAME}. The domain $\Omega = B(-r,r)$ is shown as well; none of the nodes lie inside the domain, so the choice of the discretization $Z$ is justified (see Section~\ref{sec:choiceZ}). Some methods position the nodes according to a chosen geometry: the nodes of the CME and the Euler methods are aligned on vertical lines, while the nodes of the Talbot method follow the deformed integration contour. On the other hand, the position of TAME and Zakian nodes is determined by solving respectively a minimization problem and a system of equations. It is interesting to note that the obtained nodes are arranged on curves which are close to the Talbot contour; we stress, however, that having close nodes do not make two methods similar because the weights $w_n$ can be vastly different. 

The Dirac approximants $\deltaxx_N(y)$ are plotted in Figure~\ref{fig:deltas}. The CME method has a peaked distribution at $y=1$, and is close to zero elsewhere. The other methods have a smaller and wider peak, as well as oscillation around zero in other parts of the domain; this is clearly visible for the Euler method. 
The presence of oscillation in $\deltaxx_N$ does not imply necessarily that the method is ineffective. Indeed, for \eqref{eq: Abate--Whitt through Dirac approximant} we only need convergence of distributions in the weak sense, and this kind of convergence is possible even for wildly oscillating functions, for example $\sin(ny) \rightharpoonup 0$. The Talbot method is excluded because it has nodes with $\Re(\beta_n)<0$, which make the addends $w_n e^{\beta_ny}$ extraordinarily large for $y>1$; the interpretation of $\deltaxx_N(y)$ as a Dirac approximant is not clear for Talbot.

\begin{figure}
\centering
\def\pathtodir{dataset/deltaplot}
\input{plot_deltas.tex}
\caption{Plot of the Dirac approximants $\deltaxx_N(y)$ for different methods. The TAME method uses $\Omega = B(-4,4)$.}
\label{fig:deltas}
\end{figure}
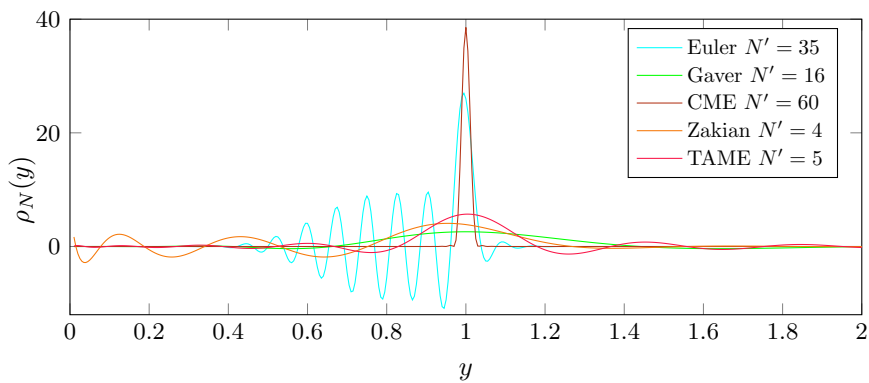
\begin{figure}
\centering
\begingroup
\setlength{\tabcolsep}{-0.2cm}
\begin{tabular}{cc}
\includegraphics[width=0.5\textwidth]{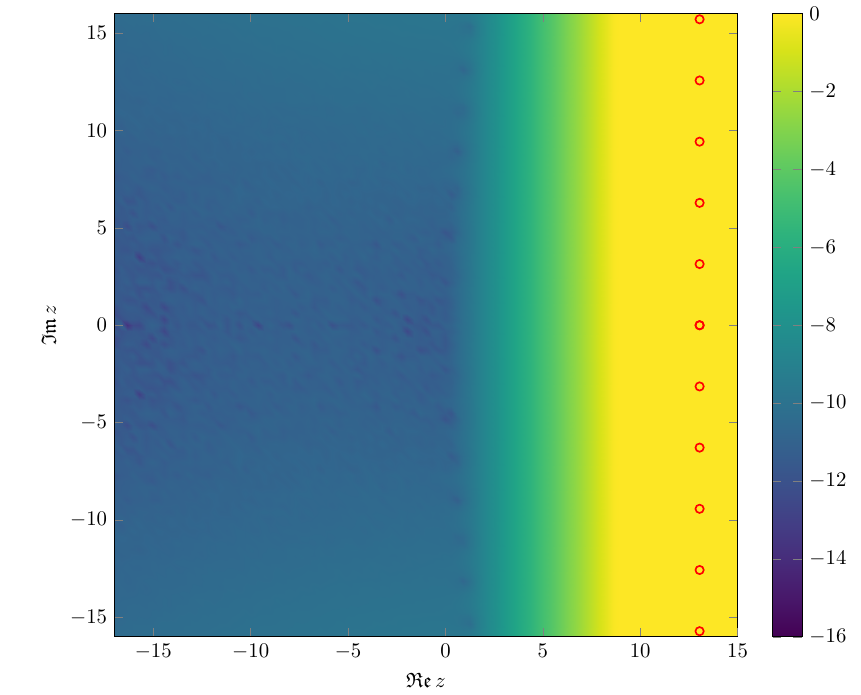}
&
\includegraphics[width=0.5\textwidth]{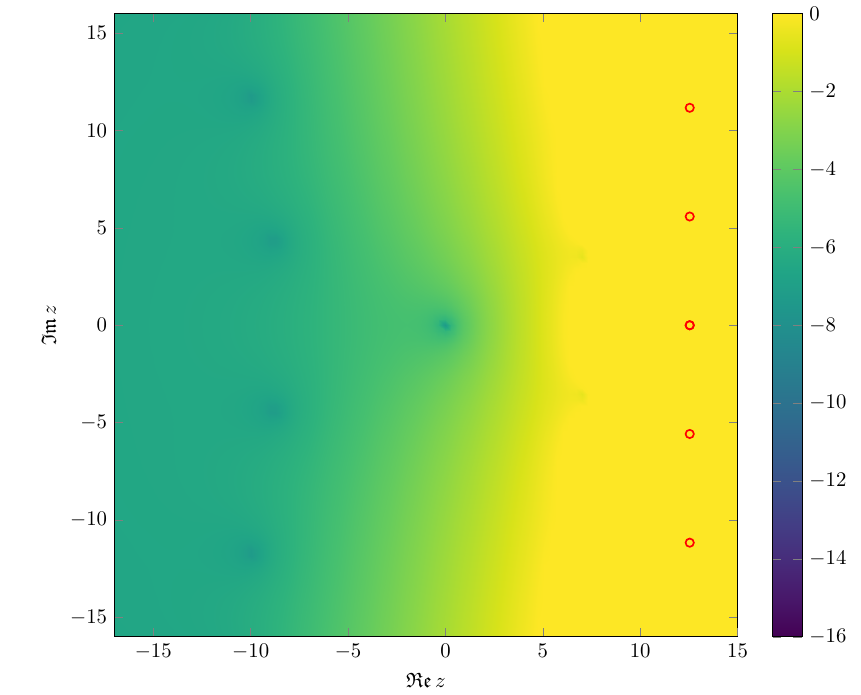}
\\
a) Euler $N'=35$. & d) CME $N'=75$.\\
\includegraphics[width=0.5\textwidth]{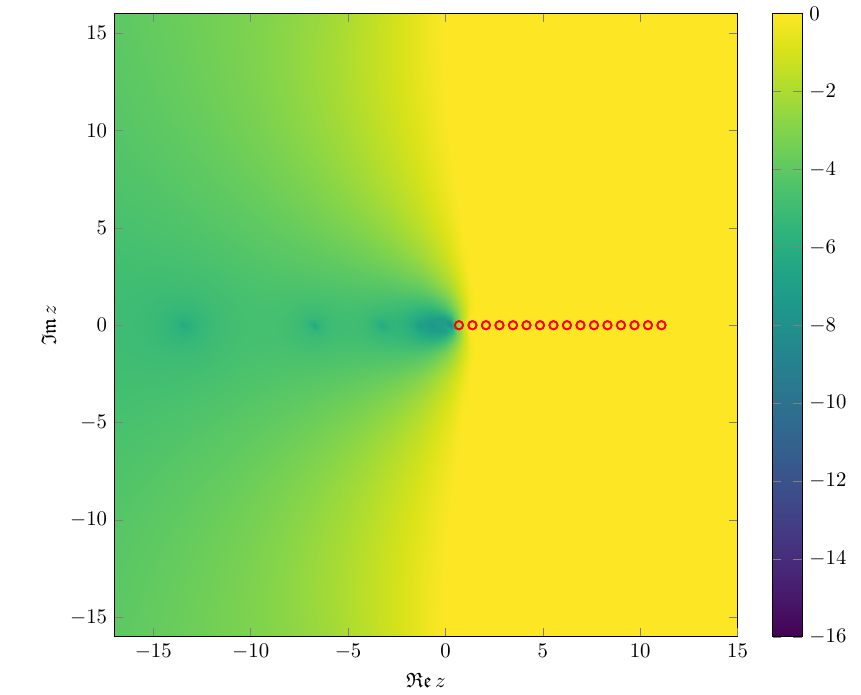}
&
\includegraphics[width=0.5\textwidth]{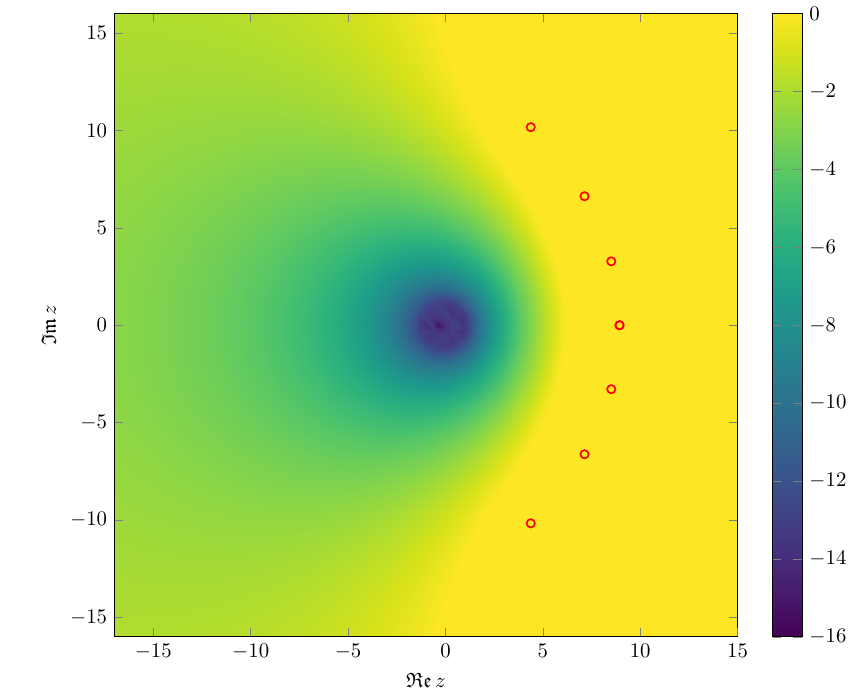}
\\
b) Gaver--Stehfest $N'=16$. & e) Zakian $N'=4$. \\
\includegraphics[width=0.5\textwidth]{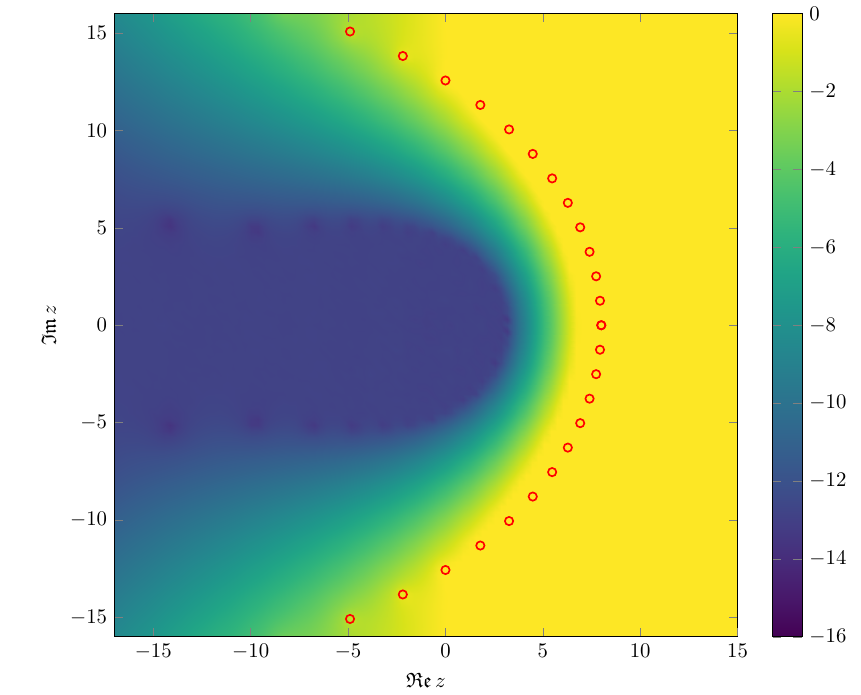}
&
\includegraphics[width=0.5\textwidth]{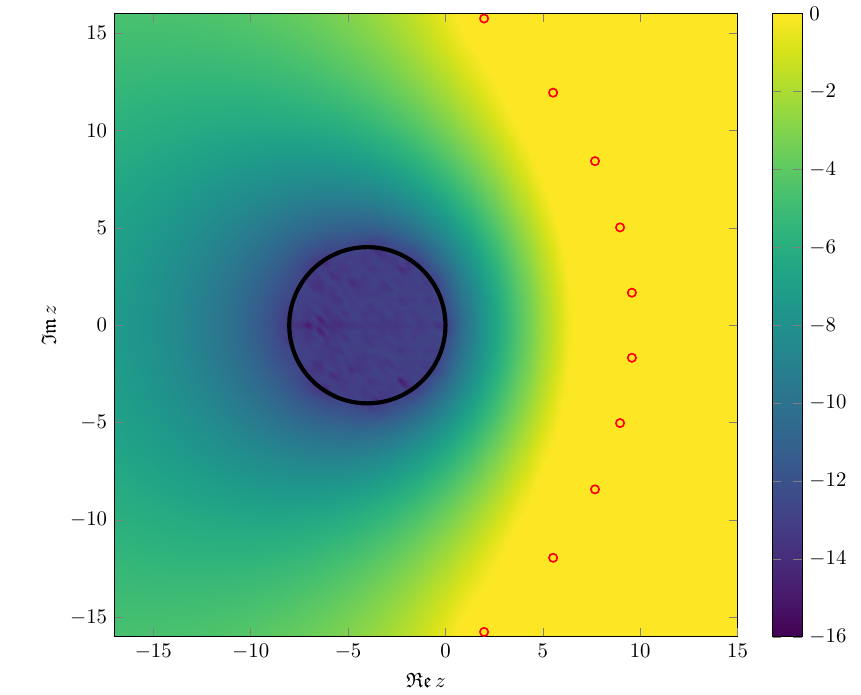}
\\
c) Talbot $N'=20$. & f) TAME $N'=5$,\\
& $\Omega=B(-4,4)$. \\
\end{tabular}
\endgroup
\caption{Base-10 logarithm of the approximation error $\log_{10}\abs*{\exp(z) - \sum_{n=1}^N \frac{w_n}{\beta_n - z}}$, and the poles $\beta$ in red.}
\label{fig:expapprox_comparison}
\end{figure}

In Figure~\ref{fig:expapprox_comparison}, we display the error of the rational approximant $\LT{\deltaxx}_N(-z)$ for different Abate--Whitt methods on a rectangular domain. The TAME method is constructed using $B(-r,r)$ with $r=4$, according to Table~\ref{tab:optimal TAME}; we display also a circle with radius $4$ on the picture. 

In the rest of this section, we present a numerical evaluation of the performance of the various Abate--Whitt methods on several problems; the first three are from the classes of functions studied in this paper, while the last two are from functions outside these classes, to determine whether the numerical properties observed and proved extend to more general problems. Our code is publicly available on~\url{https://github.com/numpi/tame-ilt}.

\subsection{Experiment A}
\begin{figure}[]    
    \centering
    \def\pathtodir{dataset/expA_PDF}
    \def\currentylabel{$\norm*{\psi(1)-\psi_N(1)}_\infty$}
    \input{plot_experA.tex}
    \def\pathtodir{dataset/expA_CDF}
    \def\currentylabel{$\norm*{\Psi(1)-\Psi_N(1)}_\infty$}
    \input{plot_experA.tex}
    \captionof{figure}{Experiment A. Plot of the error of the six Abate--Whitt methods. Above:~$\psi(t)$ (pdf), below:~$\Psi(t)$ (CDF). $\lambda=1$, $t=1$. The TAME method uses $\Omega=B(-1,1)$.}
    \label{fig:expA_PDF}
\end{figure}
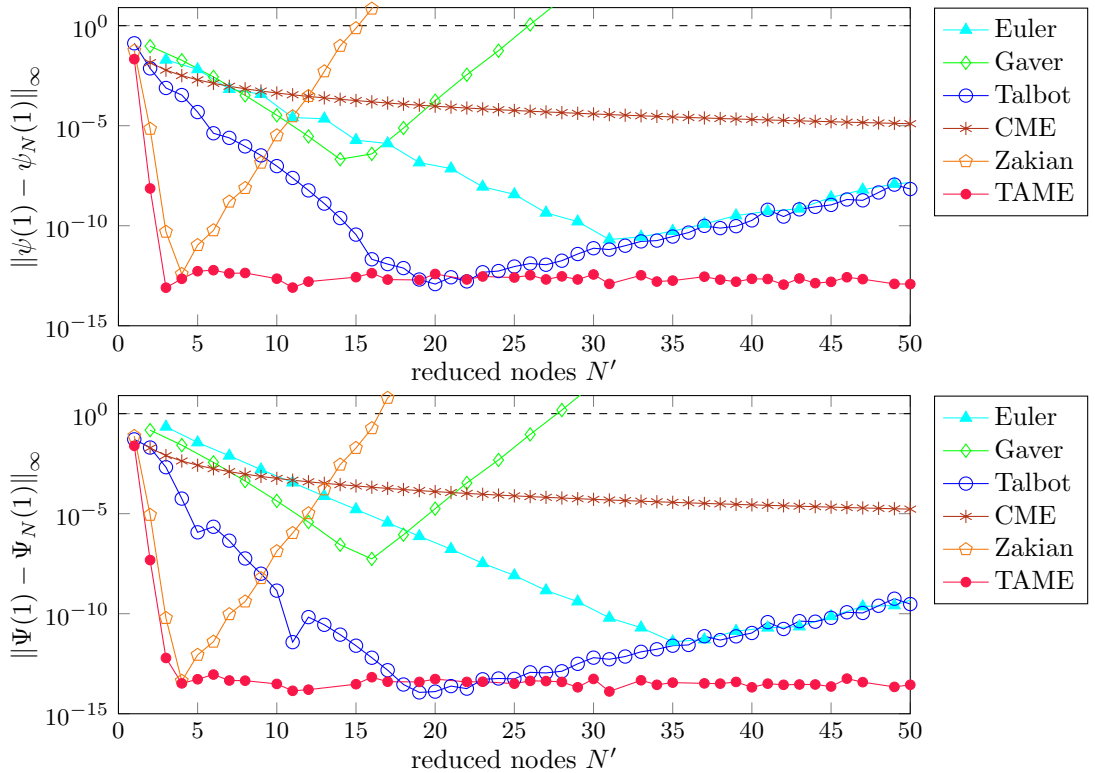

\begin{table}[]
    \centering
    \caption{Experiment A. Minimum of the approximation error for the given Abate--Whitt method and corresponding value of $N$. }\label{tab:exp01_minerror}
    \begin{tabular}{cc}
    $\Psi(t)$ \;  (CDF) &  $\psi(t)$ \;   (pdf) \\
        \hspace*{5mm}\begin{tabular}{l l r}
            \toprule
            Method  &  min error & $N$ \\
            \midrule
            Euler & $4.0 \cdot 10^{-12}$ & 35 \\
            Gaver & $5.5 \cdot 10^{-8}$ & 16 \\
            Talbot & $1.2 \cdot 10^{-14}$ & 18 \\
            CME & $1.7 \cdot 10^{-5}$ & 50 \\
            Zakian & $4.3 \cdot 10^{-14}$ & 4 \\
            TAME & $3.3 \cdot 10^{-14}$ & 4 \\
            \bottomrule
        \end{tabular}\hspace*{5mm}
    &
        \hspace*{5mm}\begin{tabular}{l l r}
            \toprule
            Method  &  min error & $N$ \\
            \midrule
            Euler & $2.0 \cdot 10^{-11}$ & 31 \\
            Gaver & $2.1 \cdot 10^{-7}$ & 14 \\
            Talbot & $1.2 \cdot 10^{-13}$ & 20 \\
            CME & $1.2 \cdot 10^{-5}$ & 50 \\
            Zakian & $3.8 \cdot 10^{-13}$ & 4 \\
            TAME & $8.0 \cdot 10^{-14}$ & 3 \\
            \bottomrule
        \end{tabular}\hspace*{5mm}
    \end{tabular}
\end{table}


\begin{figure}[]
\centering
\def\currentylabel{$\norm*{\psi(1)-\psi_N(1)}_\infty$}
\def\pathtodir{dataset/expA_PDF}
\begingroup
\setlength{\tabcolsep}{-6mm}
\hspace*{-8mm}\begin{tabular}{cc}
\def\methodName{Euler}
\def\methodColor{cyanEuler}
\def\methodMark{\markEuler}
\def\methodMarkSize{\marksizeEuler}
\def\methodError{euler.txt}
\def\methodRatAppr{RatAppr_euler.txt}
\def\methodMoments{Moments_euler.txt}
\def\legendpos{south west}
\input{plot_experA_bounds.tex}
&
\def\methodName{Gaver}
\def\methodColor{greenGS}
\def\methodMark{\markGaver}
\def\methodMarkSize{\marksizeGaver}
\def\methodError{gaver.txt}
\def\methodRatAppr{RatAppr_gaver.txt}
\def\methodMoments{Moments_gaver.txt}
\def\legendpos{south east}
\input{plot_experA_bounds.tex}
\\
\def\methodName{Talbot}
\def\methodColor{blueTalbot}
\def\methodMark{\markTalbot}
\def\methodMarkSize{\marksizeTalbot}
\def\methodError{talbot.txt}
\def\methodRatAppr{RatAppr_talbot.txt}
\def\methodMoments{Moments_talbot.txt}
\def\legendpos{north east}
\input{plot_experA_bounds.tex}
&
\def\methodName{Zakian}
\def\methodColor{orangeZakian}
\def\methodMark{\markZakian}
\def\methodMarkSize{\marksizeZakian}
\def\methodError{zakian.txt}
\def\methodRatAppr{RatAppr_zakian.txt}
\def\methodMoments{Moments_zakian.txt}
\def\legendpos{south east}
\input{plot_experA_bounds.tex}
\\
\def\methodName{CME}
\def\methodColor{bordeauxCME}
\def\methodMark{\markCME}
\def\methodMarkSize{\marksizeCME}
\def\methodError{cme.txt}
\def\methodRatAppr{RatAppr_cme.txt}
\def\methodMoments{Moments_cme.txt}
\def\legendpos{south west}
\input{plot_experA_bounds.tex}
& 
\def\methodName{TAME}
\def\methodColor{redTAME}
\def\methodMark{\markTAME}
\def\methodMarkSize{\marksizeTAME}
\def\methodError{tame1.txt}
\def\methodRatAppr{RatAppr_tame1.txt}
\def\methodMoments{Moments_tame1.txt}
\def\legendpos{north east}
\input{plot_experA_bounds.tex}\\
\end{tabular}
\endgroup
\caption{Experiment A. Upper bound, moment estimate, and  approximation error for the pdf $\psi(t)$.}
\label{fig:expA_bounds}
\end{figure}
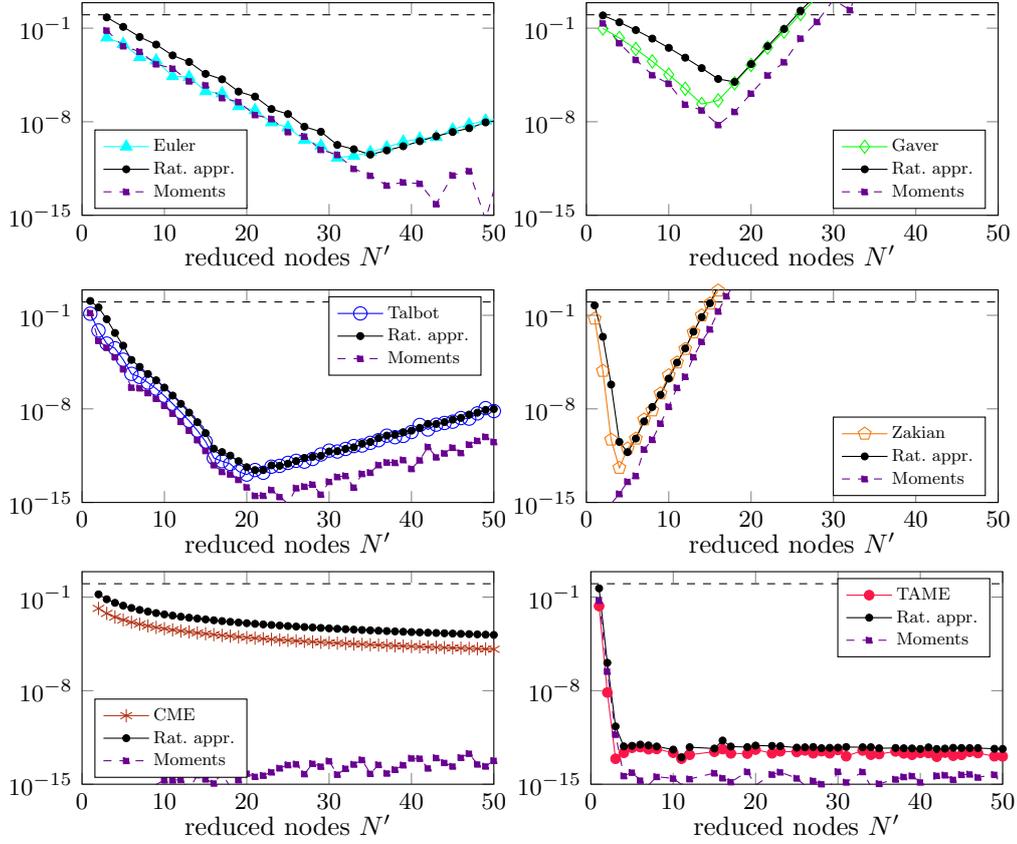

We consider a fluid queue model of size $d_+=5$, $d_-=10$ and uniformization rate $\lambda=1$. The matrix $Q$ of the underlying Markov Chain has size $d=15$. To construct $Q$, we generate a random matrix (with \texttt{abs(randn(d))} in Matlab) and then modify the diagonal entries to have zero row sums. The positive rates $r_1, \ldots, r_{d_+}$ are uniformly distributed in $[0,1]$, while the negative rates $r_{d_++1},\ldots, r_d$  in $[-1,0]$. To allow reproducibility, the seed of the Random Number Generator is set to \texttt{rng(0)}. 

The Transform $\LT{\psi}(s)$ of the pdf is calculated by solving a nonsymmetric algebraic Riccati equation (see \cite{BeaFacTay2008}). Using the properties of derivatives of the Laplace Transform (\cite[Section 1.3]{Coh2007}), the Transform of the CDF is calculated as $\LT{\Psi}(s) = \frac{1}{s} \LT{\psi}(s)$. We recover both $\Psi(t)$ and $\psi(t)$ by means of the Inverse Laplace Transform at time $t=1$. Following Theorem~\ref{thm:psibound} we choose the domain $\Omega = B(-1,1)$ for the TAME method and compare it to the classical Abate--Whitt methods.

We plot the errors $\norm*{\psi(1) - \psi_N(1)}_{\infty}$ and $\norm*{\Psi(1) - \Psi_N(1)}_{\infty}$ as $N$ increases in Figure~\ref{fig:expA_PDF}. All methods except CME exhibit a linear rate of convergence in the first part of the graph. The error of the CME method decreases approximately as $\sim (N')^{-2.21}$.  This result is consistent with the findings in \cite[Section 4.2]{HorHorTel2020}, where it was observed that the SCV of the CME method decreases approximately as $\sim (N')^{-2.14}$; therefore by Theorem~\ref{thm: HHAT2020 thm 4, Lipschitz} the error of the CME method decreases at least as fast as $\mathcal{O}((N')^{-0.71})$. In the CDF case, for $N'=50$ the error is $1.7 \cdot 10^{-5}$; the other methods reach this level of precision much sooner: the Talbot method at $N'=5$ and the TAME method at $N'=2$.

However, the classical methods soon encounter a problem: the error starts growing instead of decreasing. This is caused by an increase in the magnitude of the weights $w_n$, which causes numerical instability, as noted in Section~\ref{sec:floatingpoint}. 

Zakian's method is the fastest, but it is also the most prone to instability: after reaching the minimum at $N'=4$, the error starts rapidly rising again. For other functions the threshold $N'$ may be different, so one risks overshooting the value of $N'$ which gives a satisfying approximation. 
Talbot's method is the one that reaches the smallest error among classical methods and the growth of the error due to instability is not as fast as for Zakian.

The TAME method combines both positive aspects. It reaches an error similar to Talbot's error, and it is as fast as the Zakian method: just three evaluations of $\LT{f}(s)$ are sufficient to recover the original function $f$. Moreover, it does not suffer from numerical instability, because Algorithm~\ref{algo:AAA} constructs a good rational approximation already with $N'=4$ (i.e., with degree $N=8$), and increasing $N'$ further produces weights with an absolute value comparable to machine precision, as noted in Remark~\ref{rem:AAA vpa vs double}. While one would like to avoid using unnecessary nodes, the TAME method does not punish the user if they choose $N'$ too large.

In Figure~\ref{fig:expA_bounds} we show for each method, its error, the upper bound given by Theorem~\ref{thm:psibound}, and the moment estimate given by Definition~\ref{def:est1moment}. We see that up to the start of numerical instability, both the upper bound and the moment estimate are good approximations of the actual error for the Euler, Gaver, Talbot, and Zakian methods. For the CME method, the moment estimate is of order of machine precision, and the same happens for TAME for $N'\geq 6$. We highlight that the upper bound in Theorem~\ref{thm:psibound} depends on $N'$ only through $\varepsilon$, while the moment estimate in Definition~\ref{def:est1moment} depends on $N'$ through $|1-\mu_0|$ and $|\mu_1-\mu_0|$, which are quantities that depend only on the Abate--Whitt method and not on the function.

As further confirmed in Figure~\ref{fig:expapprox_comparison}, the upper bound is valid to measure the precision even of an Abate--Whitt method created with a different perspective in mind. Numerically also the moment estimates are sharp for all methods apart from CME; they are a surprisingly good approximation for the Talbot method. 

\subsection{Experiment B}
\label{sec:expB}
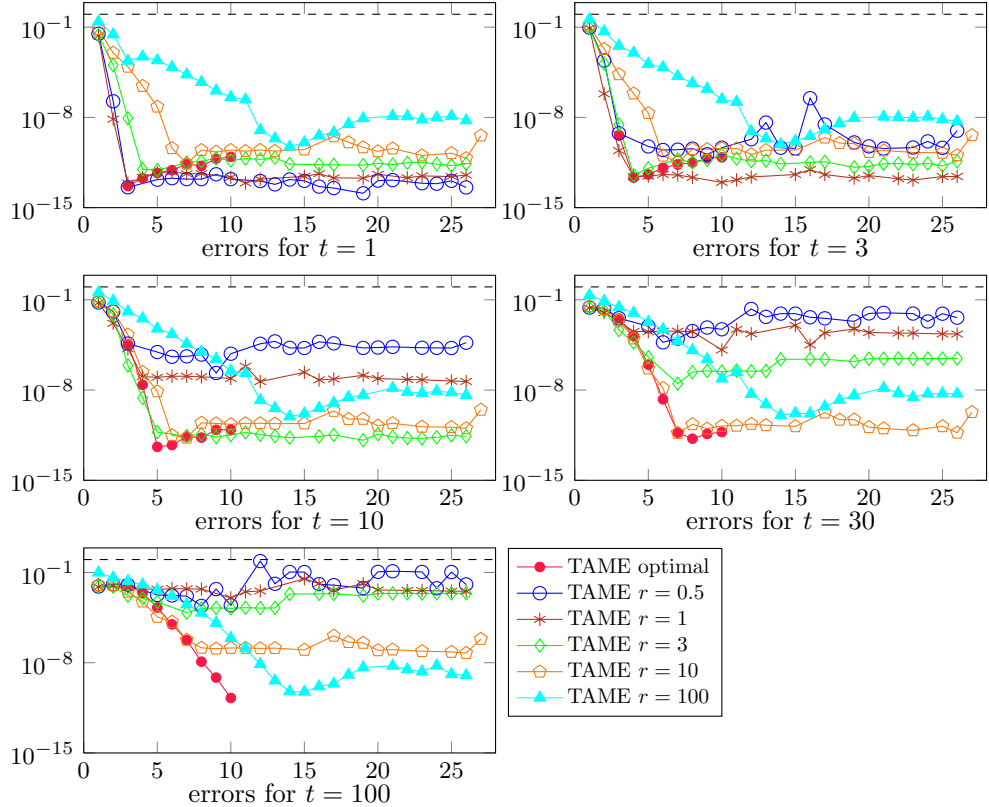
\begin{figure}[]
\centering
\def\currentylabel{}
\begingroup
\setlength{\tabcolsep}{-0.2cm}
\begin{tabular}{cc}
\def\pathtodir{dataset/expB_T=1}
\def\currentxlabel{errors for $t=1$}
\input{plot_experB.tex}
&
\def\pathtodir{dataset/expB_T=3}
\def\currentxlabel{errors for $t=3$}
\input{plot_experB.tex}
\\
\def\pathtodir{dataset/expB_T=10}
\def\currentxlabel{errors for $t=10$}
\input{plot_experB.tex}
&
\def\pathtodir{dataset/expB_T=30}
\def\currentxlabel{errors for $t=30$}
\input{plot_experB.tex}
\\
\multicolumn{2}{l}{
    \def\pathtodir{dataset/expB_T=100}
    \def\currentxlabel{errors for $t=100$}
    \input{plot_experB_withLegend.tex}
}
\end{tabular}
\endgroup
\caption{Experiment B. Error of the pdf $\psi(t)$, computed at different times $t$. Comparison of TAME methods on a circular domain $\Omega$ with different $r$.}\label{fig:expB}
\end{figure}
With the same fluid queue as in Experiment~A, we compute the pdf $\psi(t)$ at different times $t\in\{1,3,10,30,100\}$. By Theorem~\ref{thm:psibound}, the radius of the domain $\Omega = B(-r,r)$ is $r=t\lambda$ (here $\lambda=1$); we construct different TAME methods for each $r \in \{0.5,1,3,10,100\}$. We also consider the ``optimal TAME" method (see Table~\ref{tab:optimal TAME}), where the radius $r$ is not constant, but depends on the number of nodes. The comparison of these methods is shown in Figure~\ref{fig:expB}. 

In general, TAME methods with $r>t$ are quite accurate, while for $r < t$ they are not. For a fixed $t$, increasing $r$ makes the convergence slower in $N'$; this is due to the fact that larger $r$ leads to larger weights $w_n$ (see Section~\ref{sec:floatingpoint}). The TAME method with $r=100$ has the slowest convergence, but has almost the same convergence speed for any value of $t$. Conversely, methods with smaller $r$ have faster convergence, but they provide bad approximations when $t$ is much larger than $r$. 

\goodbreak

\subsection{Experiment C}
\begin{figure}[]
\centering
\def\pathtodir{dataset/expC}
\input{plot_experC.tex}
\caption{Experiment C. Target function is $f(t) = e^{tQ}$ at $t=1$.
Plot of the error of the six Abate--Whitt methods and four TAME methods as $N'$ increases.}
\label{fig:expC_errors}
\end{figure}
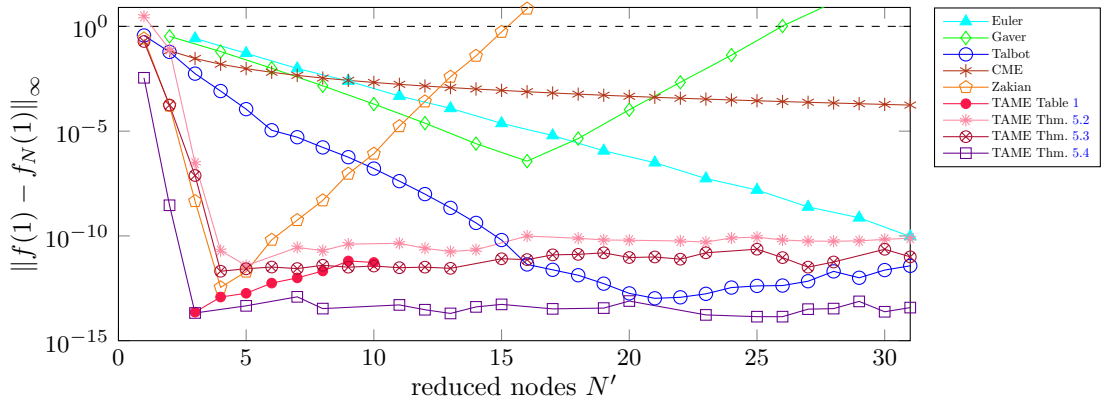
\begin{figure}
\centering
\input{plot_experC_fov} 
    \caption{Experiment C. Plot of the field of values and its bounds.} 
\label{fig:expC_fov}
\end{figure}
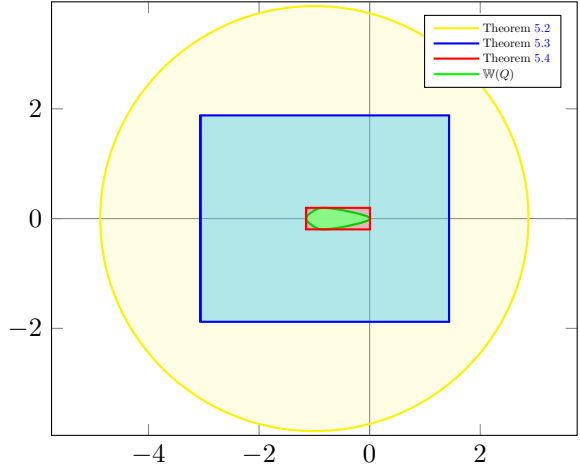
We consider the same generator matrix $Q$ of the previous two experiments, but instead of constructing the fluid queue we just analyze the underlying continuous-time Markov chain. The size of $Q$ is $d=15$ and the target function is $f(t) = e^{tQ}$, computed at time $t=1$. We use TAME methods constructed with the following domains $\Omega$:

\begin{itemize}
    \item the quasi-optimal ones in Table~\ref{tab:optimal TAME},
    \item the circle $\Omega = B(-1, \sqrt{15})$, based on Theorem~\ref{thm:Q_circle},
    \item the rectangle $\Omega = [-3.06, 1.43] + i[-1.88, 1.88]$, based on Theorem~\ref{thm:Q_rectangle_large},
    \item the smaller rectangle $\Omega = [-1.15, 0] + i[-0.19, 0.19]$, based on Theorem~\ref{thm:Q_rectangle_small}.
\end{itemize}
The three domain bounds for $\mathbb{W}(Q)$ are plotted in Figure~\ref{fig:expC_fov}.

The approximation errors are displayed in Figure~\ref{fig:expC_errors}. The four TAME methods exhibit the steepest convergence, along with the Zakian method. As the region $\Omega$ gets smaller, the error of the TAME method diminishes, consistently with the discussion in Section~\ref{sec:floatingpoint}. 


\subsection{Experiment D}
Consider two non-smooth signals: the square wave and the triangular wave. We rescale and shift them so that they assume values in the interval $[0,1]$. They are periodic functions, so they admit a Fourier series expansion on all $\R^+$. 

The triangular wave is defined as 
\[f_{\triangle}(t) = \begin{cases}
    t - \left\lfloor t\right\rfloor  & \text{ if } \mod(\left\lfloor t\right\rfloor,2) = 0,\\
    -t + \left\lfloor t\right\rfloor +1  & \text{ if } \mod(\left\lfloor t\right\rfloor,2) = 1.
\end{cases}\]
Its Laplace Transform is 
\[\LT{f_{\triangle}}(s) = \frac{1}{s^2}\frac{1-e^{-s}}{1+e^{-s}}.\]

\begin{figure}[]
\centering
\def\pathtodir{dataset/expDtriangwave}
\def\experDxmax{6}
\input{plot_experD_triangle}
\caption{Experiment D, triangular wave. Upper plot: functions $f(t)$ and $f_N(t)$. Lower plot: error in linear scale.}\label{fig:expD_triang}
\end{figure}
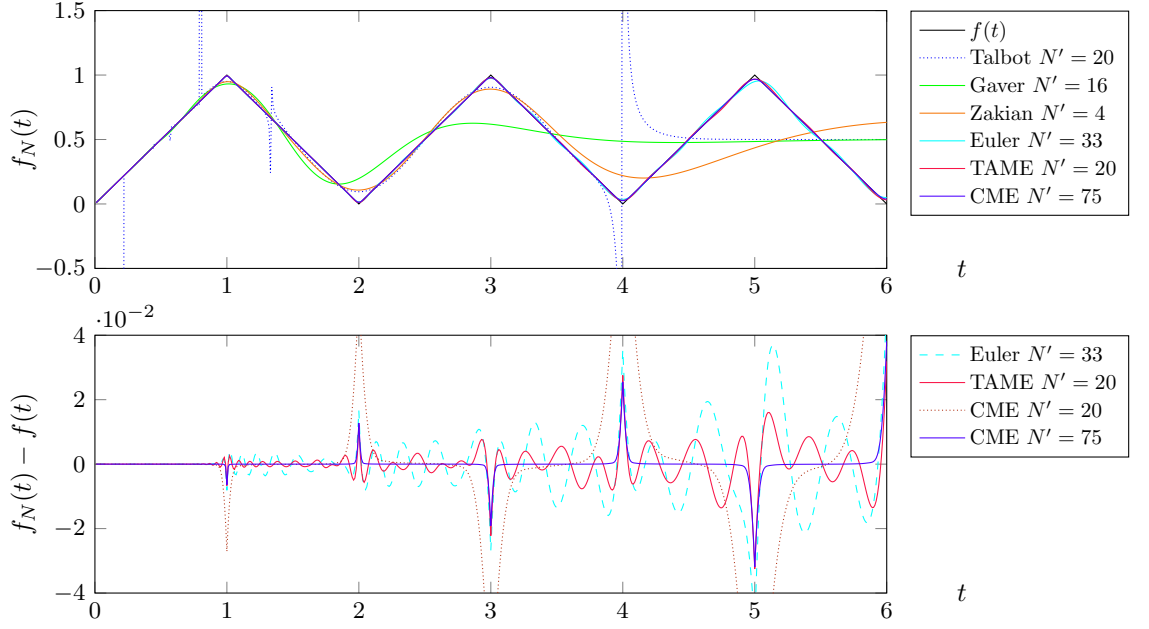

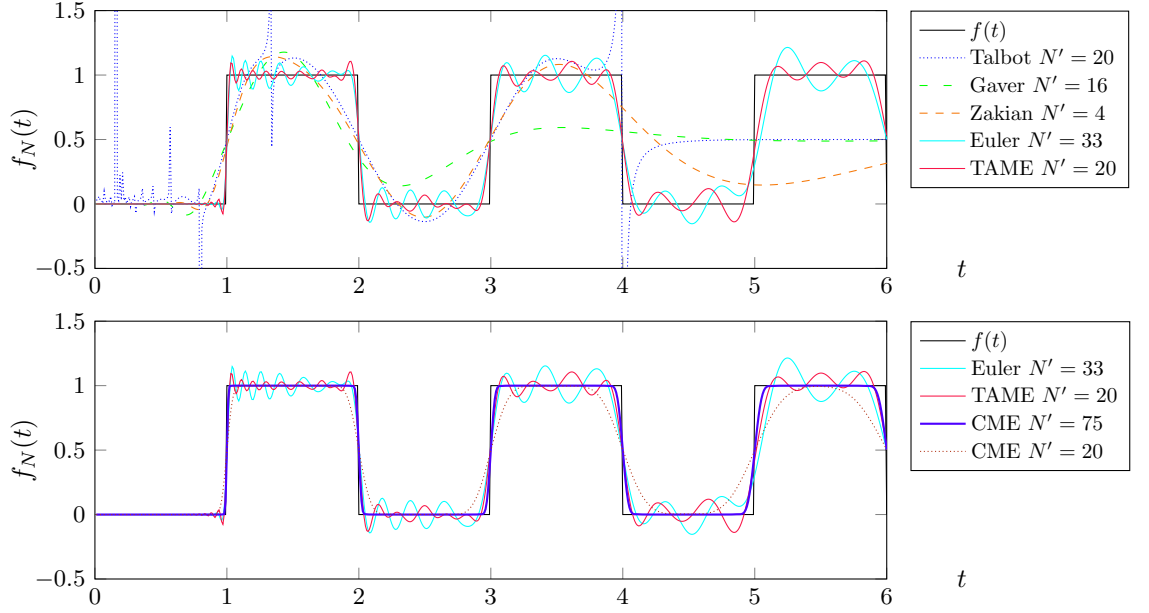
\begin{figure}[]
\centering
\def\pathtodir{dataset/expDsquarewave}
\def\experDxmax{6}
\input{plot_experD_square}
\caption{Experiment D, square wave. Function $f(t)$ and approximants are displayed $f_N(t)$. Upper plot: Euler and TAME compared to worse perfoming methods (Talbot, Gaver, Zakian). Lower plot: Euler and TAME compared to better performing CME.}\label{fig:expD_square}
\end{figure}

The triangular wave is a continuous function, but it is not differentiable everywhere. 
The $k$-th coefficient of its Fourier series is $c_k = \frac{4}{(2k+1)^2\pi^2}$. We note that the coefficients are summable and $\sum_{k=0}^\infty |c_k |= \frac{4}{\pi^2}\sum_{k=0}^\infty\frac{1}{(2k+1)^2} = \frac{1}{2}$. Therefore, the Fourier series converges absolutely to $f_{\triangle}$. The truncated Fourier series is in the SE class, so we have the necessary hypotheses to apply Theorem~\ref{thm: f approximated by SE}. The exponents of the summands of the truncated series are purely imaginary numbers, so we use $\Omega = i [-r,r]$ for a suitable chosen $r$. 

The square wave and its Laplace Transform are respectively
\[f_{\square}(t) = \bmod(\left\lfloor t\right\rfloor,2), \; \; \; \;  \LT{f_{\square}}(s) = \frac{1}{s (1+e^s)}.\]
The $k$-the coefficient of the Fourier series is $c_k = \frac{1}{2k+1}$. The square wave is discontinuous, so it cannot be uniformly approximated by any continuous function and we cannot apply Theorem~\ref{thm: f approximated by SE}. Nevertheless, the Fourier series converges at the points where $f_{\square}$ is continuous, so it is reasonable to use
a TAME method that acts accurately on the truncated Fourier series; therefore we use $\Omega = i [-r,r]$ also for the square wave. 
However, observe that the absolute series $\sum_{k=0}^\infty |c_k|$ diverges, so the approximation bound of Theorem~\ref{thm: approximation error SE} grows worse as $K$ grows.

For both problems, we have constructed a TAME method with $N'=20$ and $r=80$.

\noindent The triangular wave is shown in Figure~\ref{fig:expD_triang}. We can see that Talbot, Gaver, and Zakian do not provide a good approximation. The Talbot method has numerical issues at some values of $t$, reaching values of order of $10^{14}$. This happens because the contour curve of the Talbot integral intersects the domain $\Omega$, so the nodes $\frac{\beta}{t}$ may happen to be close to singularities of $\LT{f_{\triangle}}$. This leaves us the Euler, CME, and TAME methods. We see in the plot that the TAME method with $N'=20$ provides a better approximation that the Euler method with $N'=33$ and the CME method with $N'=20$. Unfortunately, increasing $N'$ further in the TAME method does not reduce this error, while CME with $N'=75$ is more accurate and avoids the oscillation that can be seen in the other plots. 

The square wave is shown in Figure~\ref{fig:expD_square}. As for the triangular wave, the results of Talbot, Gaver, and Zakian methods are not close to the target function. The Euler and TAME methods provide better approximations, but they suffer from Gibbs phenomena and oscillation around the correct value. In this case, CME method is the one that provides the best approximation by far, even with the same number of nodes ($N'=20$). Increasing it to $N'=75$ yields an even better approximation with the CME method.  

Indeed $f_{\square}$ is outside of our framework of ``tame" functions and, as noted, the Fourier coefficients $\frac{1}{2k+1}$ are not summable. The experiment shows that CME performs better than TAME on discontinuous functions.

\subsection{Experiment E}

We consider the European Call Option Pricing problem. For a detailed exposition of the topic, see the books \cite{Bjork,Ross_MathFin}. The goal of the Option Pricing problem is to determine the value of a contract (i.e., a \emph{call}) $C(t,Q)$ depending on time $t$ and the price of a given asset $Q$. Under standard hypotheses, $C(t,Q)$ satisfies the Black-Scholes PDE; one of the methods for its computation is through $\LT{C}(s)$ and the inverse Laplace transform. 

For the vanilla European Call Option, both $C(t)$ and $\LT{C}(s)$ are known explicitly, allowing us to use them to test the TAME method. For some exotic options, however, only $\LT{C}(s)$ is known, see e.g. \cite{Kimura2011}. The analytical solution of the European Call Option is \cite{TomShu2007}
%
%
\begin{equation}\label{eq: cost C(t)}
    C(t) = Q \, \Phi(d_+) - K e^{-R t} \, \Phi(d_-),
\end{equation}
where $Q=Q(t)$ is the asset price at current time, $R$ and $\sigma$ are parameters: $R$ is the rate of interest, $\sigma$ is the volatility. $\Phi(x)$ is the CDF of the normal distribution and 
$d_{\pm}(t) = \frac{1}{\sigma \sqrt{t}} \left(\log\left(\frac{Q}{K}\right) + \left(R\pm\frac{1}{2} \sigma^2\right)t \right).$ The Laplace Transform of $C(t)$ is (see \cite{TomShu2007})
\[
\LT{C}(s) = \begin{cases}
    \frac{K}{\gamma_+-\gamma_-} \left(\frac{S}{K}\right)^{\gamma_-}\left( \frac{\gamma_+}{R+s} - \frac{\gamma_+-1}{s}\right)  + \frac{Q}{s} - \frac{K}{R+s} & \;\;\text{ if } \;Q \geq K, \\
    \frac{K}{\gamma_+-\gamma_-} \left(\frac{S}{K}\right)^{\gamma_+}\left( \frac{\gamma_-}{R+s} - \frac{\gamma_--1}{s}\right)  & \;\;\text{ if }\; Q < K,
\end{cases}
\]
where $\gamma_+(s)\geq\gamma_-(s)$ are $\gamma_{\pm}(s) = \frac{1}{\sigma^2} \left(-\left(R-\frac{1}{2} \sigma^2\right) \pm \sqrt{\left(R-\frac{1}{2} \sigma^2\right)^2+2\sigma^2(R+s)}\right).$
The function $C(t)$ is in none of the three classes we described. Nevertheless, the summands $\Phi(d_{\pm}(t))$ of \eqref{eq: cost C(t)} are integrals of $e^{-x^2}$ on a domain determined by $d_{\pm}(t)$, which is similar to the structure of a LS function. For this reason, we construct  TAME methods using $\Omega = [-r,0]$. 

\begin{figure}[]
\def\pathtodir{dataset/expE}
\def\experDxmax{50}
\centering
\input{plot_experE}

\caption{Experiment E. Upper plot: function $C(t)$. Lower plot: absolute value of the error in log scale.}
\label{fig:expE}
\end{figure}
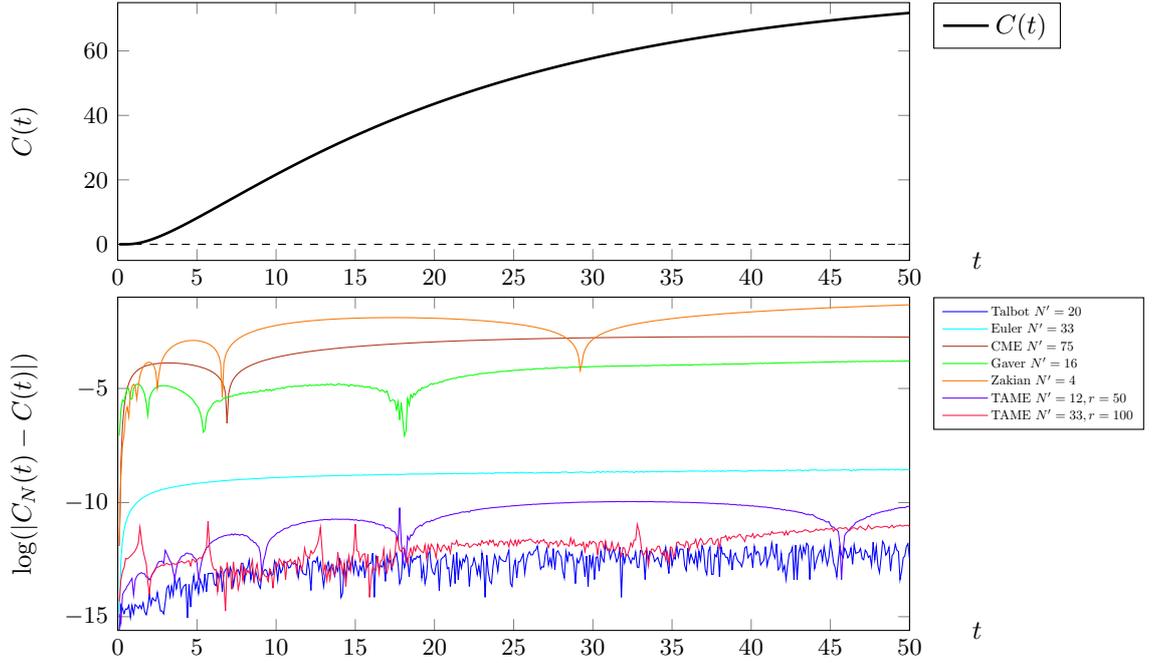

In Figure~\ref{fig:expE} we show the function $C(t)$ and the approximation errors for $K=100$, $Q=80$, $\sigma = 0.1$, $R=0.05$. The function is regular and all methods manage to provide a good enough approximation, with error smaller than $5\cdot 10^{-2}$. The most accurate methods are the Talbot with $N'=20$ and TAME with $N'=33$ nodes and $r=100$, followed by the TAME method with $N'=12$ and $r=50$, and then by the Euler method with $N'=33$.

\section{Conclusions}
We developed a theoretical framework for the accuracy analysis of Abate--Whitt methods, based on a rational approximation problem. 

We focused on functions from the SE, the ME, and the LS classes, as well as on the first return time matrix $\Psi(t)$ for fluid queues. For each class, we provided theoretical bounds for the approximation error of an Abate--Whitt method, and a description how to choose an appropriate domain $\Omega$. We showed how to adapt the AAA algorithm to construct a TAME Abate--Whitt method for each domain $\Omega$. We provide precomputed parameters for the application in fluid queues. We hope that both the analysis and the new method can be of use for practitioners in queuing theory.

In this paper, we also discussed the numerical issues related to the computation of Abate--Whitt methods in floating-point arithmetic; we plan to focus our future research on optimizing further the selection of weights to avoid excessive growth. Another open problem is extending these results to more classes of functions and to numerical methods for the ILT outside the Abate--Whitt class.

\backmatter

\bmhead{Acknowledgments} The authors are thankful to  Bernhard Beckermann for literature pointers, to Michele Benzi for useful comments, to Andrea Macrì for help and discussion regarding Experiment E, and to Francesco Zigliotto for LaTeX help. 

\bmhead{Supplementary information}
The authors are members of INDAM (Istituto Nazionale di Alta Matematica) / GNCS.

FP has received financial support by the National Centre for HPC, Big Data and Quantum Computing–HPC, CUP B83C22002940006, funded by the European Union--NextGenerationEU, and by the Italian Ministry of University and Research (MUR) through the PRIN project 2022 ``MOLE: Manifold constrained Optimization and LEarning'',  code: 2022ZK5ME7 MUR D.D. financing decree n. 20428 of November 6th, 2024 (CUP B53C24006410006).

The authors have no competing interests to disclose.

\printbibliography

\end{document}

%% file: plot_magnitudeW.tex
\begin{tikzpicture}
\begin{semilogyaxis}[
    width = 0.47\textwidth,
    xlabel={$N$},
    ylabel={$\varepsilon$},
]
\addplot[
    color=redTAME, mark=\markTAME, mark size=\marksizeTAME
] table[x=n, y=error, col sep=comma] {dataset/magnitudeW/double.txt};
\addplot[
    color=blueTalbot, mark=\markTalbot, mark size=\marksizeTalbot
] table[x=n, y=error, col sep=comma,
mark=square*] {dataset/magnitudeW/vpa.txt};
\end{semilogyaxis}
\end{tikzpicture}
\begin{tikzpicture}
\begin{semilogyaxis}[
    width = 0.47\textwidth,
    xlabel={$N$},
    ylabel={$\max \abs{w_n}$},
    legend pos = north west,
]
\addplot[
    color=redTAME, mark=\markTAME, mark size=\marksizeTAME
] table[x=n, y=maxw, col sep=comma] {dataset/magnitudeW/double.txt};
\addplot[
    color=blueTalbot, mark=\markTalbot, mark size=\marksizeTalbot
] table[x=n, y=maxw, col sep=comma] {dataset/magnitudeW/vpa.txt};
\addlegendentry{\texttt{binary64}}
\addlegendentry{\texttt{vpa}}
\end{semilogyaxis}
\end{tikzpicture}

%% file: plot_optimal_parameters.tex

\pgfplotsset{
    mycolor2/.style={color={rgb,1:red,0.8;green,0.2;blue,0.2}},    
    mycolor3/.style={color={rgb,1:red,0.2;green,0.4;blue,0.8}},    
    mycolor4/.style={color={rgb,1:red,0.2;green,0.7;blue,0.2}},    
    mycolor5/.style={color={rgb,1:red,1.0;green,0.5;blue,0.0}},    
    mycolor6/.style={color={rgb,1:red,0.6;green,0.2;blue,0.8}},    
    mycolor7/.style={color={rgb,1:red,0.5;green,0.3;blue,0.1}},    
    mycolor8/.style={color={rgb,1:red,1.0;green,0.4;blue,0.7}},    
    mycolor9/.style={color={rgb,1:red,0.4;green,0.4;blue,0.4}},    
    mycolor10/.style={color={rgb,1:red,0.5;green,0.5;blue,0.0}},   
    mycolor11/.style={color={rgb,1:red,0.0;green,0.8;blue,0.8}},   
    mycolor12/.style={color={rgb,1:red,0.8;green,0.0;blue,0.8}},   
    mycolor13/.style={color={rgb,1:red,0.4;green,0.8;blue,0.0}},   
    mycolor14/.style={color={rgb,1:red,0.9;green,0.9;blue,0.0}},   
    mycolor15/.style={color={rgb,1:red,0.0;green,0.5;blue,0.5}},   
}

\pgfplotsset{
    mymarker3/.style={mark=x},                    
    mymarker4/.style={mark=o},                    
    mymarker5/.style={mark=+},                    
    mymarker6/.style={mark=square},               
    mymarker7/.style={mark=star},                 
    mymarker8/.style={mark=triangle},             
    mymarker9/.style={mark=asterisk},             
    mymarker10/.style={mark=diamond},             
    mymarker11/.style={mark=Mercedes star},       
    mymarker12/.style={mark=pentagon},            
    mymarker13/.style={mark=x},                   
    mymarker14/.style={mark=o},                   
    mymarker15/.style={mark=+},                   
    mymarker16/.style={mark=square},              
    mymarker17/.style={mark=star},                
    mymarker18/.style={mark=triangle},            
    mymarker19/.style={mark=asterisk},            
    mymarker20/.style={mark=diamond},             
    mymarker21/.style={mark=Mercedes star},       
    mymarker22/.style={mark=pentagon},            
    mymarker23/.style={mark=x},                   
    mymarker24/.style={mark=o},                   
    mymarker25/.style={mark=+},                   
}

\begin{tikzpicture}
\begin{semilogyaxis}[
    width = \textwidth,
    height=7.6cm,
    xlabel={$r$},
    ylabel={\currentylabel},
    legend style={
        at={(0.5,-0.1)},
        anchor=north,
        nodes={scale=0.5, transform shape}
    },
    legend columns=4
]

\foreach \i in {3,...,24}{
    \foreach \j in {2,...,12} {
        \IfFileExists{./dataset/optimalpar/plot-\i-\j.txt}{
            \edef\temp{%
                \noexpand\addplot[only marks, forget plot, mycolor\j, mymarker\i]
                table[x=lambda, y=\currentycolumn, col sep=comma] 
                {./dataset/optimalpar/plot-\i-\j.txt};%
            }%
            \temp
        }{
        }
    }
}

\addlegendimage{mark size=3pt,mycolor2,mymarker3}
\addlegendentry{$N'=2,N=3$}
\addlegendimage{mark size=3pt,mycolor2,mymarker4}
\addlegendentry{$N'=2,N=4$}
\addlegendimage{mark size=3pt,mycolor3,mymarker5}
\addlegendentry{$N'=3,N=5$}
\addlegendimage{mark size=3pt,mycolor3,mymarker6}
\addlegendentry{$N'=3,N=6$}
\addlegendimage{mark size=3pt,mycolor4,mymarker7}
\addlegendentry{$N'=4,N=7$}
\addlegendimage{mark size=3pt,mycolor4,mymarker8}
\addlegendentry{$N'=4,N=8$}
\addlegendimage{mark size=3pt,mycolor5,mymarker9}
\addlegendentry{$N'=5,N=9$}
\addlegendimage{mark size=3pt,mycolor5,mymarker10}
\addlegendentry{$N'=5,N=10$}
\addlegendimage{mark size=3pt,mycolor6,mymarker11}
\addlegendentry{$N'=6,N=11$}
\addlegendimage{mark size=3pt,mycolor6,mymarker12}
\addlegendentry{$N'=6,N=12$}
\addlegendimage{mark size=3pt,mycolor7,mymarker13}
\addlegendentry{$N'=7,N=13$}
\addlegendimage{mark size=3pt,mycolor7,mymarker14}
\addlegendentry{$N'=7,N=14$}
\addlegendimage{mark size=3pt,mycolor8,mymarker15}
\addlegendentry{$N'=8,N=15$}
\addlegendimage{mark size=3pt,mycolor8,mymarker16}
\addlegendentry{$N'=8,N=16$}
\addlegendimage{mark size=3pt,mycolor9,mymarker17}
\addlegendentry{$N'=9,N=17$}
\addlegendimage{mark size=3pt,mycolor9,mymarker18}
\addlegendentry{$N'=9,N=18$}
\addlegendimage{mark size=3pt,mycolor10,mymarker19}
\addlegendentry{$N'=10,N=19$}
\addlegendimage{mark size=3pt,mycolor10,mymarker20}
\addlegendentry{$N'=10,N=20$}
\addlegendimage{mark size=3pt,mycolor11,mymarker21}
\addlegendentry{$N'=11,N=21$}
\addlegendimage{mark size=3pt,mycolor11,mymarker22}
\addlegendentry{$N'=11,N=22$}
\addlegendimage{mark size=3pt,mycolor12,mymarker23}
\addlegendentry{$N'=12,N=23$}
\addlegendimage{mark size=3pt,mycolor12,mymarker24}
\addlegendentry{$N'=12,N=24$}



\end{semilogyaxis}
\end{tikzpicture}


%% file: plot_Z-beta.tex
\begin{tikzpicture}
\begin{axis}[
    xlabel={$\Re(z)$},
    ylabel={$\Im(z)$},
    width=9.8cm,
    height=9.8cm,
    ymin=-20,
    ymax=20,
    xmin=-25,
    xmax=18,
    axis lines = middle,
    xticklabel style={
        font=\small,
    },
    yticklabel style={
        font=\small,
        /pgf/number format/fixed,
    },
    legend pos = outer north east,
    legend cell align={left},
    legend style = {at={(0.3,0.88)}, anchor = north east,nodes={scale=0.8, transform shape}},
    mark options={solid},
]

\filldraw[color=blue!70!black, fill=blue!95!black,thick, on layer=axis background](axis cs:-4,0) circle (4);
node at (-10,0) {ciao};

\addplot [color=blue!70!black, mark=none, ultra thick] (0,0);

\addplot [color=cyanEuler, mark=\markEuler, mark size=3pt, dashed] table {dataset/nodes/nodes_euler_n35_conj.txt};
\addplot [color=greenGS, mark=\markGaver, mark size=2.7pt, dashed] table {dataset/nodes/nodes_gaver_n16.txt};
\addplot [color=blueTalbot, mark=\markTalbot, mark size = 2.3pt, dashed] table {dataset/nodes/nodes_talbot_n20_conj.txt};
\addplot [color=bordeauxCME, mark=\markCME, mark size = 3pt, dashed] table {dataset/nodes/nodes_cme_n60_conj.txt};
\addplot [color=orangeZakian, mark=\markZakian, mark size = 2.6pt, dashed] table {dataset/nodes/nodes_zakian_n4_conj.txt};
\addplot [color=redTAME, mark=\markTAME, mark size=2.1, dashed] table {dataset/nodes/nodes_tame_n5_r4_conj.txt};

\addlegendentry{domain $\Omega$}
\addlegendentry{Euler $N'=35$ }
\addlegendentry{Gaver $N'=16$ }
\addlegendentry{Talbot $N'=20$ }
\addlegendentry{CME $N'=60$ }
\addlegendentry{Zakian $N'=4$ }
\addlegendentry{TAME $N'=5$ }

\end{axis}
\end{tikzpicture}

%% file: plot_deltas.tex
\begin{tikzpicture}
\centering
\begin{axis}[
    xlabel={$y$},
    ylabel={$\deltaxx_N(y)$},
    width=12cm,
    height=5.5cm,
    ymin=-12,
    ymax=40,
    xmin=0,
    xmax=2,
    ylabel style={
        yshift=-6mm, 
        xshift=-3mm, 
    },
    xticklabel style={
        font=\small,
    },
    yticklabel style={
        font=\small,
        /pgf/number format/fixed,
    },
    legend pos = north east,
    legend cell align={left},
    legend style = {nodes={scale=0.8, transform shape}},
]


    \addplot [color=cyanEuler, mark=none] table {\pathtodir/delta_euler_n35.txt};
    \addlegendentry{Euler $N'=35$}

    \addplot [color=greenGS, mark=none] table {\pathtodir/delta_gaver_n16.txt};
    \addlegendentry{Gaver $N'=16$}


    \addplot [color=bordeauxCME, mark=none] table {\pathtodir/delta_cme_n60.txt};
    \addlegendentry{CME $N'=60$}


    \addplot [color=orangeZakian, mark=none] table {\pathtodir/delta_zakian_n04.txt};
    \addlegendentry{Zakian $N'=4$}

    \addplot [color=redTAME, mark=none] table {\pathtodir/delta_tame_n05_circle_r4.txt};
    \addlegendentry{TAME $N'=5$}

\end{axis}
\end{tikzpicture}

%% file: plot_experA.tex
\begin{tikzpicture}
\begin{semilogyaxis}[
    xlabel={reduced nodes $N'$},
    ylabel={\currentylabel},
    width=12cm,
    height=5.8cm,
    ymin=10^-15,
    ymax=8,
    xmin=0,
    xmax=50,
    xlabel style={
        yshift=2mm, 
    },
    xticklabel style={
        font=\small,
    },
    yticklabel style={
        font=\small,
        /pgf/number format/fixed,
    },
    legend pos = outer north east,
    legend cell align={left},
]

\addplot[mark=none, dashed, thin, black, samples=2, domain=0:50, forget plot] {1};
\addplot [color=cyanEuler, mark=\markEuler, mark size = \marksizeEuler] table {\pathtodir/euler.txt};
\addplot [color=greenGS, mark=\markGaver,  mark size = \marksizeGaver] table {\pathtodir/gaver.txt};
\addplot [color=blueTalbot, mark=\markTalbot, mark size = \marksizeTalbot] table {\pathtodir/talbot.txt};
\addplot [color=bordeauxCME, mark=\markCME, mark size = \marksizeCME] table {\pathtodir/cme.txt};
\addplot [color=orangeZakian, mark=\markZakian, mark size = \marksizeZakian] table {\pathtodir/zakian.txt};
\addplot [color=redTAME, mark=\markTAME, mark size = \marksizeTAME] table {\pathtodir/tame1.txt};

\addlegendentry{Euler}
\addlegendentry{Gaver}
\addlegendentry{Talbot}
\addlegendentry{CME}
\addlegendentry{Zakian}
\addlegendentry{TAME}

\end{semilogyaxis}
\end{tikzpicture}

%% file: plot_experA_bounds.tex
\begin{tikzpicture}
\begin{semilogyaxis}[
    xlabel={reduced nodes $N'$},
    width=7cm,
    height=4.4cm,
    ymin=10^-15,
    ymax=8,
    xmin=0,
    xmax=50,
    xlabel style={
        yshift=2.5mm, 
    },
    xticklabel style={
        font=\small,
    },
    yticklabel style={
        font=\small,
        /pgf/number format/fixed,
    },
    legend pos=\legendpos,
    legend cell align={left},
    legend style={font=\scriptsize,nodes={scale=0.73, transform shape}},
]
\definecolor{purpleBounds}{rgb}{0.406, 0, 0.566}
\addplot[mark=none, dashed, thin, black, samples=2, domain=0:50, forget plot] {1};
\addplot [color=\methodColor, mark=\methodMark, mark size = \methodMarkSize] table {\pathtodir/\methodError};
\addplot [color=black, mark=*, mark size = 1.3pt] table {\pathtodir/\methodRatAppr};
\addplot [color=purpleBounds, mark=square*, mark size = 1.1pt, dashed] table {\pathtodir/\methodMoments}; 

\addlegendentry{\methodName}
\addlegendentry{Rat. appr.}
\addlegendentry{Moments}

\end{semilogyaxis}
\end{tikzpicture}

%% file: plot_experB.tex
\begin{tikzpicture}
\begin{semilogyaxis}[
    xlabel={\currentxlabel},
    ylabel={},
    width=7cm,
    height=4.3cm,
    ymin=10^-15,
    ymax=8,
    xmin=0,
    xmax=28,
    xlabel style={
        yshift=2.5mm, 
    },
    xticklabel style={
        font=\small,
    },
    yticklabel style={
        font=\small,
        /pgf/number format/fixed,
    },
    legend pos = north east,
    legend cell align={left},
]

\addplot[mark=none, dashed, thin, black, samples=2, domain=0:30, forget plot] {1};
\addplot [color=redTAME, mark=\markTAME, mark size = \marksizeTAME] table {\pathtodir/tame1.txt};
\addplot [color=blueTalbot, mark=\markTalbot, mark size = \marksizeTalbot] table {\pathtodir/tame2.txt};
\addplot [color=bordeauxCME, mark=\markCME, mark size = \marksizeCME] table {\pathtodir/tame3.txt};
\addplot [color=greenGS, mark=\markGaver,  mark size = \marksizeGaver] table {\pathtodir/tame4.txt};
\addplot [color=orangeZakian, mark=\markZakian, mark size = \marksizeZakian] table {\pathtodir/tame5.txt};
\addplot [color=cyanEuler, mark=\markEuler, mark size = \marksizeEuler] table {\pathtodir/tame6.txt};



\end{semilogyaxis}
\end{tikzpicture}

%% file: plot_experB_withLegend.tex
\begin{tikzpicture}
\begin{semilogyaxis}[
    xlabel={\currentxlabel},
    ylabel={},
    width=7cm,
    height=4.3cm,
    ymin=10^-15,
    ymax=8,
    xmin=0,
    xmax=28,
    xlabel style={
        yshift=2.5mm, 
    },
    xticklabel style={
        font=\small,
    },
    yticklabel style={
        font=\small,
        /pgf/number format/fixed,
    },
    legend pos = outer north east,
    legend cell align={left},
    legend style={nodes={scale=0.8, transform shape}},
]

\addplot[mark=none, dashed, thin, black, samples=2, domain=0:50, forget plot] {1};
\addplot [color=redTAME, mark=\markTAME, mark size = \marksizeTAME] table {\pathtodir/tame1.txt};
\addplot [color=blueTalbot, mark=\markTalbot, mark size = \marksizeTalbot] table {\pathtodir/tame2.txt};
\addplot [color=bordeauxCME, mark=\markCME, mark size = \marksizeCME] table {\pathtodir/tame3.txt};
\addplot [color=greenGS, mark=\markGaver,  mark size = \marksizeGaver] table {\pathtodir/tame4.txt};
\addplot [color=orangeZakian, mark=\markZakian, mark size = \marksizeZakian] table {\pathtodir/tame5.txt};
\addplot [color=cyanEuler, mark=\markEuler, mark size = \marksizeEuler] table {\pathtodir/tame6.txt};

\addlegendentry{TAME optimal}
\addlegendentry{TAME $r=0.5$}
\addlegendentry{TAME $r=1$}
\addlegendentry{TAME $r=3$}
\addlegendentry{TAME $r=10$}
\addlegendentry{TAME $r=100$}
\end{semilogyaxis}
\end{tikzpicture}

%% file: plot_experC.tex
\begin{tikzpicture}
\begin{semilogyaxis}[
    xlabel={reduced nodes $N'$},
    ylabel={$\norm*{f(1)-f_N(1)}_\infty$},
    width=12cm,
    height=6cm,
    ymin=10^-15,ymax=8,
    xmin=0,xmax=31,
    xlabel style={
        yshift=2mm, 
    },
    xticklabel style={
        font=\small,
    },
    yticklabel style={
        font=\small,
        /pgf/number format/fixed,
    },
    legend pos = outer north east,
    legend cell align={left},
    legend style={nodes={scale=0.5, transform shape}},
]

\definecolor{purpleTAME}{rgb}{0.406, 0, 0.566}

\addplot[mark=none, dashed, thin, black, samples=2, domain=0:50, forget plot] {1};
\addplot [color=cyanEuler, mark=\markEuler, mark size = \marksizeEuler] table {\pathtodir/euler.txt};
\addplot [color=greenGS, mark=\markGaver,  mark size = \marksizeGaver] table {\pathtodir/gaver.txt};
\addplot [color=blueTalbot, mark=\markTalbot, mark size = \marksizeTalbot] table {\pathtodir/talbot.txt};
\addplot [color=bordeauxCME, mark=\markCME, mark size = \marksizeCME] table {\pathtodir/cme.txt};
\addplot [color=orangeZakian, mark=\markZakian, mark size = \marksizeZakian] table {\pathtodir/zakian.txt};
\addplot [color=redTAME, mark=\markTAME, mark size = \marksizeTAME] table {\pathtodir/tame1.txt};
\addplot [color=redTAME!50!white, mark=10-pointed star, mark size = 2.3pt] table {\pathtodir/tame2.txt};
\addplot [color=redTAME!70!black, mark=otimes, mark size = 2.2pt] table {\pathtodir/tame3.txt};
\addplot [color=purpleTAME, mark=square, mark size = 2pt] table {\pathtodir/tame4.txt};

\addlegendentry{Euler}
\addlegendentry{Gaver}
\addlegendentry{Talbot}
\addlegendentry{CME}
\addlegendentry{Zakian}
\addlegendentry{TAME Table~\ref{tab:optimal TAME}}
\addlegendentry{TAME Thm.~\ref{thm:Q_circle}}
\addlegendentry{TAME Thm.~\ref{thm:Q_rectangle_large}}
\addlegendentry{TAME Thm.~\ref{thm:Q_rectangle_small}}

\end{semilogyaxis}
\end{tikzpicture}

%% file: plot_experC_fov.tex
\begin{tikzpicture}
\begin{axis}[
   axis equal,
   enlargelimits={false, 0.01},
   width=0.65\textwidth,
   legend pos = {north east},
   legend cell align={left},
   legend style={nodes={scale=0.46, transform shape}},
]

\draw[gray, thin] (axis cs:\pgfkeysvalueof{/pgfplots/xmin},0) -- (axis cs:\pgfkeysvalueof{/pgfplots/xmax},0);
\draw[gray, thin] (axis cs:0,\pgfkeysvalueof{/pgfplots/ymin}) -- (axis cs:0,\pgfkeysvalueof{/pgfplots/ymax});

\draw[
    draw=yellow,
    fill=yellow!80, 
    fill opacity = 0.15, 
    thick, 
    radius={transformdirectionx(sqrt(15))}
] (axis cs:-1,0) circle;    
\addplot[domain=0:2*pi,samples=5, draw = none, forget plot] ({-1+sqrt(15)*cos(deg(x))}, {sqrt(15)*sin(deg(x))}); 
\addlegendimage{fill=yellow, draw=yellow, opacity=1, thick}
\addlegendentry{Theorem~\ref{thm:Q_circle}}

\addplot[
   draw=blue, 
   fill=cyan!60, 
   fill opacity=0.4, 
   thick,
] coordinates {(-3.0616,-1.8810) (-3.0616,1.8810) (1.4365,1.8810) (1.4365,-1.8810) (-3.0616,-1.8810) (-3.0616,1.8810)};
\addlegendentry{Theorem~\ref{thm:Q_rectangle_large}}

\addplot[
   draw=none, 
   fill=red!40, 
   fill opacity=0.6, 
   forget plot,
] coordinates {(-1.1503,-0.1964) (-1.1503,0.1964) ( 0.0074,0.1964) ( 0.0074,-0.1964) (-1.1503,-0.1964) (-1.1503,0.1964)};

\addplot[
   draw=green!80!black, 
   fill=green!50, 
   fill opacity=0.9, 
   thick,
   forget plot,
] 
table{
    0.0074   -0.0004
    0.0074   -0.0009
    0.0073   -0.0013
    0.0073   -0.0015
    0.0073   -0.0021
    0.0071   -0.0030
    0.0071   -0.0032
    0.0069   -0.0038
    0.0068   -0.0042
    0.0067   -0.0047
    0.0065   -0.0053
    0.0064   -0.0055
    0.0063   -0.0057
    0.0061   -0.0061
    0.0056   -0.0073
    0.0053   -0.0080
    0.0050   -0.0085
    0.0048   -0.0088
    0.0044   -0.0095
    0.0042   -0.0098
    0.0038   -0.0104
    0.0034   -0.0109
    0.0032   -0.0112
    0.0022   -0.0125
    0.0016   -0.0132
    0.0013   -0.0135
    0.0003   -0.0146
   -0.0005   -0.0154
   -0.0009   -0.0158
   -0.0018   -0.0167
   -0.0023   -0.0171
   -0.0028   -0.0176
   -0.0046   -0.0190
   -0.0059   -0.0201
   -0.0093   -0.0225
   -0.0103   -0.0232
   -0.0113   -0.0239
   -0.0125   -0.0246
   -0.0165   -0.0271
   -0.0198   -0.0290
   -0.0285   -0.0335
   -0.0379   -0.0378
   -0.0572   -0.0456
   -0.0639   -0.0482
   -0.0813   -0.0543
   -0.0925   -0.0580
   -0.1230   -0.0675
   -0.1438   -0.0735
   -0.2033   -0.0895
   -0.2451   -0.1000
   -0.2966   -0.1121
   -0.4231   -0.1395
   -0.4889   -0.1525
   -0.5992   -0.1722
   -0.6726   -0.1834
   -0.7200   -0.1895
   -0.7379   -0.1914
   -0.7667   -0.1939
   -0.7788   -0.1948
   -0.7898   -0.1954
   -0.8407   -0.1962
   -0.8474   -0.1960
   -0.8537   -0.1957
   -0.8801   -0.1936
   -0.9072   -0.1893
   -0.9105   -0.1886
   -0.9137   -0.1879
   -0.9197   -0.1865
   -0.9227   -0.1857
   -0.9255   -0.1849
   -0.9339   -0.1824
   -0.9366   -0.1815
   -0.9420   -0.1796
   -0.9474   -0.1775
   -0.9501   -0.1764
   -0.9529   -0.1752
   -0.9557   -0.1740
   -0.9586   -0.1727
   -0.9648   -0.1698
   -0.9681   -0.1682
   -0.9716   -0.1664
   -0.9753   -0.1644
   -0.9887   -0.1569
   -1.0151   -0.1406
   -1.0330   -0.1286
   -1.0526   -0.1147
   -1.0621   -0.1077
   -1.0709   -0.1010
   -1.0789   -0.0948
   -1.0859   -0.0892
   -1.1019   -0.0759
   -1.1058   -0.0725
   -1.1092   -0.0694
   -1.1122   -0.0666
   -1.1213   -0.0576
   -1.1231   -0.0558
   -1.1247   -0.0540
   -1.1262   -0.0524
   -1.1312   -0.0466
   -1.1322   -0.0452
   -1.1341   -0.0427
   -1.1350   -0.0415
   -1.1358   -0.0404
   -1.1366   -0.0393
   -1.1374   -0.0381
   -1.1400   -0.0340
   -1.1417   -0.0311
   -1.1422   -0.0301
   -1.1431   -0.0283
   -1.1444   -0.0257
   -1.1452   -0.0240
   -1.1455   -0.0231
   -1.1458   -0.0223
   -1.1468   -0.0199
   -1.1470   -0.0191
   -1.1473   -0.0183
   -1.1475   -0.0175
   -1.1478   -0.0167
   -1.1482   -0.0152
   -1.1484   -0.0144
   -1.1486   -0.0137
   -1.1489   -0.0122
   -1.1492   -0.0107
   -1.1494   -0.0100
   -1.1495   -0.0093
   -1.1496   -0.0085
   -1.1499   -0.0064
   -1.1500   -0.0057
   -1.1500   -0.0050
   -1.1502   -0.0028
   -1.1502   -0.0021
   -1.1502   -0.0014
   -1.1502    0.0014
   -1.1502    0.0028
   -1.1502    0.0035
   -1.1501    0.0042
   -1.1500    0.0050
   -1.1500    0.0057
   -1.1497    0.0078
   -1.1495    0.0093
   -1.1494    0.0100
   -1.1492    0.0107
   -1.1489    0.0122
   -1.1488    0.0129
   -1.1486    0.0137
   -1.1484    0.0144
   -1.1482    0.0152
   -1.1478    0.0167
   -1.1470    0.0191
   -1.1465    0.0207
   -1.1458    0.0223
   -1.1455    0.0231
   -1.1452    0.0240
   -1.1448    0.0248
   -1.1440    0.0265
   -1.1436    0.0274
   -1.1431    0.0283
   -1.1422    0.0301
   -1.1412    0.0320
   -1.1406    0.0330
   -1.1394    0.0350
   -1.1381    0.0371
   -1.1374    0.0381
   -1.1366    0.0393
   -1.1350    0.0415
   -1.1341    0.0427
   -1.1312    0.0466
   -1.1300    0.0479
   -1.1288    0.0494
   -1.1262    0.0524
   -1.1247    0.0540
   -1.1231    0.0558
   -1.1213    0.0576
   -1.1194    0.0596
   -1.1173    0.0617
   -1.1122    0.0666
   -1.1092    0.0694
   -1.0973    0.0798
   -1.0859    0.0892
   -1.0789    0.0948
   -1.0709    0.1010
   -1.0621    0.1077
   -1.0526    0.1147
   -1.0330    0.1286
   -1.0237    0.1349
   -1.0004    0.1499
   -0.9942    0.1537
   -0.9887    0.1569
   -0.9794    0.1622
   -0.9753    0.1644
   -0.9648    0.1698
   -0.9557    0.1740
   -0.9529    0.1752
   -0.9474    0.1775
   -0.9420    0.1796
   -0.9366    0.1815
   -0.9339    0.1824
   -0.9255    0.1849
   -0.9227    0.1857
   -0.9137    0.1879
   -0.9038    0.1900
   -0.9003    0.1907
   -0.8928    0.1919
   -0.8887    0.1925
   -0.8845    0.1930
   -0.8754    0.1941
   -0.8652    0.1950
   -0.8537    0.1957
   -0.8474    0.1960
   -0.8407    0.1962
   -0.8336    0.1963
   -0.8260    0.1964
   -0.8092    0.1961
   -0.7788    0.1948
   -0.7533    0.1928
   -0.7379    0.1914
   -0.5992    0.1722
   -0.5487    0.1635
   -0.4889    0.1525
   -0.2451    0.1000
   -0.1438    0.0735
   -0.1230    0.0675
   -0.0813    0.0543
   -0.0572    0.0456
   -0.0513    0.0434
   -0.0462    0.0413
   -0.0379    0.0378
   -0.0344    0.0362
   -0.0313    0.0348
   -0.0285    0.0335
   -0.0217    0.0300
   -0.0165    0.0271
   -0.0113    0.0239
   -0.0103    0.0232
   -0.0093    0.0225
   -0.0075    0.0213
   -0.0067    0.0207
   -0.0059    0.0201
   -0.0028    0.0176
   -0.0023    0.0171
   -0.0009    0.0158
   -0.0005    0.0154
   -0.0001    0.0150
    0.0007    0.0142
    0.0019    0.0128
    0.0022    0.0125
    0.0032    0.0112
    0.0036    0.0107
    0.0038    0.0104
    0.0042    0.0098
    0.0053    0.0080
    0.0056    0.0073
    0.0057    0.0071
    0.0059    0.0066
    0.0060    0.0064
    0.0063    0.0057
    0.0065    0.0053
    0.0065    0.0051
    0.0067    0.0047
    0.0068    0.0042
    0.0070    0.0036
    0.0070    0.0034
    0.0071    0.0032
    0.0072    0.0027
    0.0072    0.0025
    0.0073    0.0017
    0.0073    0.0015
    0.0073    0.0013
    0.0074    0.0009
    0.0074    0.0008
    0.0074    0.0004
    0.0074    0.0002
    0.0074    0.0000
};

\addplot[
   draw=red, 
   thick,
] coordinates {(-1.1503,-0.1964) (-1.1503,0.1964) ( 0.0074,0.1964) ( 0.0074,-0.1964) (-1.1503,-0.1964) (-1.1503,0.1964)};
\addlegendentry{Theorem~\ref{thm:Q_rectangle_small}}

\addlegendimage{fill=red, draw=green, opacity=1, thick}
\addlegendentry{$\mathbb{W}(Q)$}

\end{axis}
\end{tikzpicture}

%% file: plot_experD_triangle.tex
\begin{tikzpicture}
    \begin{axis}[
        xlabel={$t$},
        ylabel={$f_N(t)$},
        width=12cm, height=5cm,
        ymin=-0.5, ymax=1.5,
        xmin=0, xmax=\experDxmax,
        xlabel style={yshift=7.6mm, xshift=6.2cm,},
        ylabel style={yshift=-3mm,},  
        xticklabel style={font=\small,},
        yticklabel style={font=\small,/pgf/number format/fixed,},
        legend pos = outer north east,
        legend cell align={left},
        legend style={nodes={scale=0.8, transform shape}},
    ]
    \addplot [color=black, mark=none, line width = 0.4pt] table {\pathtodir/function_values.txt};
    \addlegendentry{$f(t)$}
    \addplot [color=blueTalbot, mark=none, densely dotted] table {\pathtodir/talbotcutoff_n20.txt};\addlegendentry{Talbot $N'=20$}
    \addplot [color=greenGS, mark=none] table {\pathtodir/gaver_n16.txt};\addlegendentry{Gaver $N'=16$}
    \addplot [color=orangeZakian, mark=none] table {\pathtodir/zakian_n04.txt};\addlegendentry{Zakian $N'=4$}
    \addplot [color=cyanEuler, mark=none] table {\pathtodir/euler_n33.txt};\addlegendentry{Euler $N'=33$}
    \addplot [color=redTAME, mark=none] table {\pathtodir/tame1.txt};\addlegendentry{TAME $N'=20$}
    \addplot [color=purpleCME, mark=none] table {\pathtodir/cme_n75.txt};\addlegendentry{CME $N'=75$}
    \end{axis}
\end{tikzpicture}

\begin{tikzpicture}
    \begin{axis}[
        xlabel={$t$},
        ylabel={$f_N(t)-f(t)$},
        width=12cm, height=5cm,
        ymin=-0.04, ymax=0.04,
        xmin=0, xmax=\experDxmax,
        xlabel style={yshift=7.6mm, xshift=6.2cm,},
        ylabel style={yshift=-3mm,},  
        xticklabel style={font=\small,},
        yticklabel style={font=\small,/pgf/number format/fixed,},
        legend pos = outer north east,
        legend cell align={left},
        legend style={nodes={scale=0.8, transform shape}},
    ]

    \addplot [color=cyanEuler, mark=none, dashed] table {\pathtodir/err_euler_n33.txt};\addlegendentry{Euler $N'=33$}
    \addplot [color=redTAME, mark=none] table {\pathtodir/err_tame1.txt};\addlegendentry{TAME $N'=20$}
    \addplot [color=bordeauxCME, mark=none, densely dotted] table {\pathtodir/err_cme_n20.txt};\addlegendentry{CME $N'=20$}
    \addplot [color=purpleCME, mark=none] table {\pathtodir/err_cme_n75.txt};\addlegendentry{CME $N'=75$}

    \end{axis}
\end{tikzpicture}

%% file: plot_experD_square.tex
\begin{tikzpicture}
    \begin{axis}[
        xlabel={$t$},
        ylabel={$f_N(t)$},
        width=12cm, height=5cm,
        ymin=-0.5, ymax=1.5,
        xmin=0, xmax=\experDxmax,
        xlabel style={yshift=7.6mm, xshift=6.2cm,},
        ylabel style={yshift=-3mm,},  
        xticklabel style={font=\small,},
        yticklabel style={font=\small,/pgf/number format/fixed,},
        legend pos = outer north east,
        legend cell align={left},
        legend style={nodes={scale=0.8, transform shape}},
    ]
    \addplot [color=black, mark=none, line width = 0.4pt] table {\pathtodir/function_values.txt};
    \addlegendentry{$f(t)$}
    \addplot [color=blueTalbot, mark=none, densely dotted] table {\pathtodir/talbotcutoff_n20.txt};\addlegendentry{Talbot $N'=20$}
    \addplot [color=greenGS, mark=none, loosely dashed] table {\pathtodir/gaver_n16.txt};\addlegendentry{Gaver $N'=16$}
    \addplot [color=orangeZakian, mark=none, dashed] table {\pathtodir/zakian_n04.txt};\addlegendentry{Zakian $N'=4$}
    \addplot [color=cyanEuler, mark=none] table {\pathtodir/euler_n33.txt};\addlegendentry{Euler $N'=33$}
    \addplot [color=redTAME, mark=none] table {\pathtodir/tame1.txt};\addlegendentry{TAME $N'=20$}
    \end{axis}
\end{tikzpicture}

\begin{tikzpicture}
    \begin{axis}[
        xlabel={$t$},
        ylabel={$f_N(t)$},
        width=12cm, height=5cm,
        ymin=-0.5, ymax=1.5,
        xmin=0, xmax=\experDxmax,
        xlabel style={yshift=7.6mm, xshift=6.2cm,},
        ylabel style={yshift=-3mm,},  
        xticklabel style={font=\small,},
        yticklabel style={font=\small,/pgf/number format/fixed,},
        legend pos = outer north east,
        legend cell align={left},
        legend style={nodes={scale=0.8, transform shape}},
    ]
    \addplot [color=black, mark=none, line width = 0.4pt] table {\pathtodir/function_values.txt};
    \addlegendentry{$f(t)$}
    \addplot [color=cyanEuler, mark=none] table {\pathtodir/euler_n33.txt};\addlegendentry{Euler $N'=33$}
    \addplot [color=redTAME, mark=none] table {\pathtodir/tame1.txt};\addlegendentry{TAME $N'=20$}
    \addplot [color=purpleCME, mark=none, thick] table {\pathtodir/cme_n75.txt};\addlegendentry{CME $N'=75$}
    \addplot [color=bordeauxCME, mark=none, densely dotted] table {\pathtodir/cme_n20.txt};\addlegendentry{CME $N'=20$}
    \end{axis}
\end{tikzpicture}

%% file: plot_experE.tex
\begin{tikzpicture}
    \begin{axis}[
        xlabel={$t$},
        ylabel={$C(t)$},
        width=12cm, height=5cm,
        ymin=-5, ymax=75,
        xmin=0, xmax=\experDxmax,
        xlabel style={yshift=7.6mm, xshift=6.1cm,},
        xticklabel style={font=\small,},
        yticklabel style={font=\small,/pgf/number format/fixed,},
        legend pos = outer north east,
    ]
    \addplot[mark=none, dashed, thin, black, samples=2, domain=0:\experDxmax, forget plot] {0};
    \addplot [color=black, mark=none, line width = 1pt] table {\pathtodir/function_values.txt};
    \addlegendentry{$C(t)$}
    \end{axis}
\end{tikzpicture}

\begin{tikzpicture}
    \begin{axis}[
        xlabel={$t$},
        ylabel={$\log (\left| C_N(t)-C(t)\right|) $},
        width=12cm, height=6cm,
        ymin=-15.6, ymax=-1,
        xmin=0, xmax=\experDxmax,
        xlabel style={yshift=7.6mm, xshift=6.1cm,},        
        xticklabel style={font=\small,},
        yticklabel style={font=\small,/pgf/number format/fixed,},
        legend pos = outer north east,
        legend cell align={left},
        legend style={nodes={scale=0.5, transform shape}},
    ]
    \addplot [color=blueTalbot, mark=none] table {\pathtodir/logerr_talbot_n20.txt};\addlegendentry{Talbot $N'=20$}
    \addplot [color=cyanEuler, mark=none] table {\pathtodir/logerr_euler_n33.txt};\addlegendentry{Euler $N'=33$}
    \addplot [color=bordeauxCME, mark=none] table {\pathtodir/logerr_cme_n75.txt};\addlegendentry{CME $N'=75$}
    \addplot [color=greenGS, mark=none] table {\pathtodir/logerr_gaver_n16.txt};\addlegendentry{Gaver $N'=16$}
    \addplot [color=orangeZakian, mark=none] table {\pathtodir/logerr_zakian_n04.txt};\addlegendentry{Zakian $N'=4$}
    \addplot [color=purpleTAME, mark=none] table {\pathtodir/logerr_tame2.txt};    \addlegendentry{TAME $N'=12, r=50$}
    \addplot [color=redTAME, mark=none] table {\pathtodir/logerr_tame1.txt};    \addlegendentry{TAME $N'=33, r=100$}
    \end{axis}
\end{tikzpicture}